\documentclass[12pt]{amsart}
\usepackage{epsfig,color}

\makeatother

\theoremstyle{definition}

\def\fnum{equation}
\newtheorem{Thm}[\fnum]{Theorem}
\newtheorem{Cor}[\fnum]{Corollary}

\newtheorem{Lem}[\fnum]{Lemma}

\newtheorem{Def}[\fnum]{Definition}
\newtheorem{Exa}[\fnum]{Example}
\newtheorem{Rem}[\fnum]{Remark}
\newtheorem{Pro}[\fnum]{Proposition}

\numberwithin{equation}{section}

\newcommand{\HH}{{\mathcal{H}}}
\newcommand{\HHH}{{\bar{\mathcal{H}}}}

\newcommand{\nn}{{\bf{n}}}
\newcommand{\Ric}{{\text{Ric}}}

\newcommand{\Low}{{\mathcal{L}}}
\newcommand{\Loww}{{\bar{\mathcal{L}}}}

\newcommand{\Hess}{{\text {Hess}}}
\def\ZZ{{\bold Z}}
\def\RR{{\bold R}}
\def\SS{{\bold S}}

\def\CC{{\bold C }}
\newcommand{\dv}{{\text {div}}}
\newcommand{\e}{{\text {e}}}

\newcommand{\Area}{{\text {Area}}}
\newcommand{\Length}{{\text {Length}}}

\newcommand{\cB}{{\mathcal{B}}}

\newcommand{\cL}{{\mathcal{L}}}
\newcommand{\cH}{{\mathcal{H}}}

\newcommand{\eqr}[1]{(\ref{#1})}

\newcommand{\elltwo}{{{J}}}

\begin{document}

\title[Schr\"odinger operators on manifolds with
cylindrical ends]{Three circles theorems for Schr\"odinger
operators on  cylindrical ends and geometric applications}

\author{Tobias H. Colding}%
\address{MIT\\
77 Massachusetts Avenue, Cambridge, MA 02139-4307\\
and Courant Institute of Mathematical Sciences\\
251 Mercer Street, New York, NY 10012.}
\author{Camillo De Lellis}
\address{Institut f\"ur Mathematik, Universit\"at Z\"urich\\
Winterthurerstr. 190\\ 8057 Z\"urich (CH)}
\author{William P. Minicozzi II}%
\address{Department of Mathematics\\
Johns Hopkins University\\
3400 N. Charles St.\\
Baltimore, MD 21218}

\thanks{The first and the third author
were partially supported by NSF Grants DMS 0104453, DMS 0606629
and DMS 0405695}


\email{colding@math.mit.edu, delellis@math.unizh.ch, and
minicozz@math.jhu.edu}

\maketitle

\begin{abstract}
We show that for a Schr\"odinger operator with   bounded potential
on a manifold with cylindrical ends the space of solutions which
grows at most exponentially at infinity is finite dimensional and,
for a dense set of potentials (or, equivalently for a surface, for
a fixed potential and a dense set of metrics), the constant
function zero is the only solution that vanishes at infinity.
Clearly, for general potentials there can be many solutions that
vanish at infinity.

These results follow from a three circles inequality (or log
convexity inequality) for the Sobolev norm of a solution $u$ to a
Schr\"odinger equation on a product $N\times [0,T]$, where $N$ is
a closed manifold with a certain spectral gap.  Examples of such
$N$'s are all (round) spheres $\SS^n$ for $n\geq 1$ and all Zoll
surfaces.

   Finally, we discuss some examples arising in geometry of such
manifolds and Schr\"odinger operators.
\end{abstract}

\section{Introduction}

Many problems in Geometric Analysis are about the space of
solutions of non-linear PDE's, like solutions of the Yang-Mills
equation, the Einstein equation, the Yamabe equation, the harmonic
map equation, the minimal surface equation, etc. For such problems
it is often of interest to estimate ``how many'' solutions there
are and be able to say something about their properties.
Infinitesimally, the space of nearby solutions to a given solution
solve a linear PDE,  which is often a Schr\"odinger equation.  For
this reason it is therefore very useful when one can say that the
space of solutions (with some constraints at infinity) to a
Schr\"odinger equation is finite dimensional and even more
significant when one can say that the trivial solution, that is,
the function that is identically zero, is the only such solution.
The first case corresponds to that the ``tangent space'' is finite
dimensional and the second case corresponds to that the space of
solutions is infinitesimally rigid.  We will return to some
specific examples later in the introduction after stating our main
results.

Let $M$ be a complete non-compact  $(n+1)$-dimensional Riemannian
manifold with  finitely many ends $E_1 , \dots , E_k$. Suppose
also that $M \setminus \bigcup_{j=1}^k E_j$ has compact closure
and each end is cylindrical.  By cylindrical we will mean
different things depending on whether $n=1$, in which case more
general ends will be allowed, or $n\geq 2$. For $n\geq 2$ we
assume that each end $E_i$ is isometric to a product of a closed
manifold $N_i$ and a half-line $[0,\infty)$, whereas, for $n=1$ we
assume only that each end is bi-Lipschitz to $\SS^1\times
[0,\infty)$ and has bounded geometry. Recall that a surface (or
manifold) has bounded geometry if its sectional curvature is
bounded above and below and the injectivity radius is bounded away
from zero.

We will consider Schr\"odinger operators $L = \Delta_{M} + V$ on
the manifold $M$ and on each cylindrical end use coordinates
 $(\theta , t )$.
Given a constant $\alpha$, let $H_{\alpha} (M)= H_{\alpha} (M, L)$
be the linear space of all solutions $u$ of $L u = 0$ that grow
slower than $\exp (\alpha r)$, where $r$ is the distance to a
fixed point. That is, for any fixed point $p$
\begin{equation}    \label{e:defhal}
    \limsup_{r\to \infty}
    \max_{\partial \cB_r(p)} \e^{-\alpha \, r} |u|  = 0 \, ,
\end{equation}
Where $\cB_r(p)$ is the intrinsic ball of radius $r$ and center
$p$. Note that   $H_0 (M)$ is the set of solutions that vanish at
infinity.

One of our main results is the next theorem about the solutions of
Schr\"odinger operators on manifolds with cylindrical ends, where
the cross-section of each end has a (infinite) sequence of
eigenvalues $\lambda_{m_i}$ for the Laplacian with
\begin{equation}    \label{e:spgaps}
\lambda_{m_i}-\lambda_{m_i-1}\to \infty \, .
\end{equation}
 This
last condition on the spectral gaps is satisfied on any round
sphere $\SS^n$ for $n\geq 1$.  On $\SS^n$, the eigenvalues occur
with multiplicity in clusters with the $m$-th cluster at $m^2 +
(n-1)\, m$.  The spectral gap condition is also satisfied on any
Zoll surface (normalized so the closed geodesics have length
$2\pi$).  The eigenvalues of a Zoll surface occur in clusters,
where the eigenvalues in the $m$-th cluster all lie in the
interval
\begin{equation}
 J_m = \left[ (m + \beta/4)^2 - K , (m + \beta/4)^2 + K \right]
\end{equation}
for constants $K$ and $\beta$; see
 Guillemin, \cite{Gu}, and Colin de Verdi\`ere, \cite{Cv}.  Notice
 that the gap between $J_m$ and $J_{m+1}$  grows linearly in $m$,
  as did the spectral gaps for $\SS^n$,
 thus giving the required spectral gap.{\footnote{Weyl's
 asymptotic formula gives for a general closed $n$-dimensional
 manifold that $\lambda_m \approx m^{n/2}$, so   \eqr{e:spgaps}
 does not   hold in general for $n \geq 2$.}}

\begin{Thm}     \label{t:main}
Let $M$ be a complete non-compact $(n+1)$-dimensional manifold
with finitely many cylindrical ends satisfying \eqr{e:spgaps}.
\begin{enumerate}
\item[(1)] If $V$ is a $C^{0,1}$ bounded{\footnote{A function $f$
is in $C^{0,1}$ if it is both bounded and Lipschitz.  The
$C^{0,1}$ norm  is
$$
    ||f||_{C^{0,1}} = \sup_M |f| + \sup_{x\ne y \in M} \, \frac{|f(x)
    - f(y)|}{|x-y|} \, .
$$}}
   function{\footnote{We will prove that both parts (1) and (2) of the
theorem also hold for bounded potentials $V$, whenever the
cross-section of each end is a round $\SS^n$, $n\geq 1$,   or a
Zoll surface.}} (potential) on $M$, then $H_{\alpha} (M,
\Delta_{M} +V)$ is finite dimensional for every $\alpha$; the
bound for   $\dim H_{\alpha}$ depends only on $M$, $\alpha$, and
$||V||_{C^{0,1}}$. \item[(2)] For a dense set of $C^{0,1}$ bounded
potentials $H_{0}$ contains only the constant function zero; for a
surface this is equivalent to that, for a fixed potential, there
is a dense set of metrics (with finitely many cylindrical ends)
where $H_{0}= \{ 0 \}$.
\end{enumerate}
\end{Thm}

 For easy reference, we state also this theorem in the special
 case of surfaces.

\begin{Thm}
Let $M$ be a complete non-compact surface with finitely many
cylindrical ends.
\begin{enumerate}
\item[(1)] If $V$ is a bounded
   function (potential) on $M$, then $H_{\alpha} (M,
\Delta_{M} +V)$ is finite dimensional for every $\alpha$; the
bound for   $\dim H_{\alpha}$ depends only on $M$, $\alpha$, and
$||V||_{L^{\infty}}$. \item[(2)] For a dense set of bounded
potentials, $H_{0}$ contains only the constant function zero; or
equivalently, for a fixed potential, there is a dense set of
metrics (with finitely many cylindrical ends) where $H_{0}= \{ 0
\}$.
\end{enumerate}
\end{Thm}

Even the special case of our theorem where $M=\SS^1\times \RR$ is
a flat cylinder is of interest.  In that case we can define spaces
$H_+$ and $H_-$ of solutions to the Schr\"odinger equation where
$H_+$ are the solutions that vanish at $+\infty$ and $H_-$ the
space that vanishes at $-\infty$ and thus $H_0$ is the
interesection of the two.  In this case both $H_+$ and $H_-$ can
be infinite dimensional, as can be seen when  $V\equiv 0$ by
considering separation of variable solutions:
\begin{align}
     \{ \e^{kt} \,
\cos ( k\theta) {\text{ and }} \e^{kt} \,
\sin ( k\theta) \, | \, k\in \ZZ , k < 0 \} &\subset H_+ \, , \\
 \{ \e^{kt} \,
\cos ( k\theta) {\text{ and }} \e^{kt} \, \sin ( k\theta) \, | \,
k\in \ZZ , k > 0 \} &\subset H_- \, .
\end{align}
In particular, one can easily construct (non-generic) compactly
supported potentials $V$ on the flat cylinder $\SS^1\times \RR$
where $H_0$ is non-trivial by patching together exponentially
decaying solutions on each end.

One of the key ingredients  in the proof of Theorem \ref{t:main}
is a three circles inequality (or log convexity inequality) for
the Sobolev norm of a solution $u$ to a Schr\"odinger equation on
a product $N\times [0,T]$, where $N$ satisfies \eqr{e:spgaps}.  We
will state the first version of the three circles theorem next
when $N$ is a sphere or a Zoll surface and the dependence of the
constants is cleanest; see Theorem \ref{t:ests0} below for the
statement for a general $N$ satisfying \eqr{e:spgaps}.

\begin{Thm}     \label{t:ests}
Let $N = \SS^n$ for any $n \geq 1$ or a Zoll surface. There exists
a constant $C
> 0$ depending on  $N$ and  $||V||_{C^{0,1}}$ so that if $u$ is a
solution to the Schr\"odinger equation $\Delta u+V\,u=0$ on
$N\times [0,T]$ and $\alpha$ satisfies
\begin{equation}    \label{e:defalphafirst0}
\alpha \geq \frac{1}{T}\left[\log \frac{I(T)}{I(0)}\right]\, ,
\end{equation}
then $u$'s $W^{1,2}$ norm at $0<t<T$ satisfies the following three
circles type inequality (logarithmic convexity type inequality)
\begin{equation}    \label{e:eo8}
    \log I(t)\leq   C   +    (C + |\alpha|) \, t + \log I(0)
    \, .
\end{equation}
Here
\begin{equation}
I(s)=\int_{N\times \{s\} } \left( u^2 + |\nabla u|^2 \right)
\,d\theta\, .
\end{equation}
\end{Thm}

Our argument actually gives a stronger bound than we record in
Theorem \ref{t:ests0}, but we have tailored the statement to fit
our geometric applications.

 Even if the potential is merely bounded, and not
 Lipschitz, we get the following estimate:

 \begin{Thm}     \label{t:estsc0}
  Let $N = \SS^n$ for any $n \geq 1$ or a Zoll surface.
There exists a constant $C
> 0$ depending on  $N$ and    $||V||_{L^{\infty}}$ so that if $u$ is a
solution to the Schr\"odinger equation $\Delta u+V\,u=0$ on
$N\times [0,T]$ and $\alpha$ satisfies
\begin{equation}    \label{e:defalphafirst1}
\alpha \geq \frac{1}{T}\left[\log \frac{\int_{N \times \{ T \} }
u^2 }{\int_{N \times \{ 0 \} } u^2 }\right] \, ,
\end{equation}
then $u$'s $L^{2}$ norm at $0<t<T$ satisfies
\begin{equation}    \label{e:eo11}
    \log \left( \int_{ N \times \{ t \} }  u^2 \, d\theta \right)
        \leq C + (C + |\alpha|) \, t + \log I(0)\, .
\end{equation}
\end{Thm}

One of the main reasons why such estimates are useful is that it
shows that if a solution grows/decays initially with at least a
certain rate (the constant $C$ in \eqr{e:eo8}  and \eqr{e:eo11}
gives a threshold), then it will keep growing/decaying
indefinitely.

As an immediate corollary of the general version of Theorem
\ref{t:ests} where $N$ is only assumed to satisfy \eqr{e:spgaps},
i.e., Theorem \ref{t:ests0}, (and Schauder estimates) we get the
following:

\begin{Cor}     \label{c:ests}
Let $N$ be a closed $n$-dimensional manifold satisfying
\eqr{e:spgaps}. Given $\alpha \in \RR$, there exists a constant
$\nu > 0$ depending on $\alpha$, the $C^{0,1}$ norm of $V$, and
$N$ so that if $u \in H_{\alpha}(N \times [0,\infty))$, then its
$W^{1,2}$ norm grows at most exponentially with the estimate
\begin{equation}
    \int_{N\times \{t\} } \left( u^2 + |\nabla u|^2 \right)  \,d\theta \leq
    \nu \, \e^{\nu \, t} \, \int_{N\times \{0\} } \left( u^2 + |\nabla u|^2
\right)
    \,d\theta  \, .
\end{equation}
\end{Cor}

\begin{Rem}
The corollary also holds for bounded potentials $V$ whenever $N$
is an $n$-dimensional sphere or a Zoll surface; in this case, we
apply Theorem \ref{t:estsc0}.
\end{Rem}

\subsection{Examples from geometry}
Let $\Sigma\subset M^3$ be a  smooth surface (possibly with
boundary) in a complete Riemannian $3$-manifold $M$ and with
orientable normal bundle. Given a function $\phi$ in the space
$C^{\infty}_0(\Sigma)$ of infinitely differentiable (i.e.,
smooth), compactly supported functions on $\Sigma$, consider the
one-parameter variation
\begin{equation}
\Sigma_{t,\phi}=\{x+\exp_x(t\,\phi (x)\,\nn_{\Sigma}(x)) | x\in
\Sigma\}\, .
\end{equation}
Here $\nn_{\Sigma}$ is the unit normal to $\Sigma$ and $\exp$ is
the exponential map on $M$.{\footnote{For instance, if $M=\RR^3$,
then $\exp_x(v)=x+v$.}} The so-called first variation formula of
area is the equation (integration is with respect to the area of
$\Sigma$)
\begin{equation}  \label{e:frstvar}
\left.\frac{d}{dt} \right|_{t=0}\Area (\Sigma_{t,\phi})
=\int_{\Sigma}\phi\,H\, ,
\end{equation}
where the {\emph{mean curvature}} $H$   of $\Sigma$ is the sum of
the principal curvatures $\kappa_1$, $\kappa_2$.{\footnote{When
$\Sigma$ is non-compact, $\Sigma_{t,\phi}$ in \eqr{e:frstvar} is
replaced by $\Gamma_{t,\phi}$, where $\Gamma$ is any compact set
containing the support of $\phi$.}} The surface $\Sigma$ is said
to be a {\it minimal} surface (or just minimal) if
\begin{equation}
\left.\frac{d}{dt} \right|_{t=0}\Area (\Sigma_{t,\phi})=0
\,\,\,\,\,\,\,\,\,\,\,\text{ for all } \phi\in
C^{\infty}_0(\Sigma)
\end{equation}
 or, equivalently by \eqr{e:frstvar}, if the
mean curvature $H$ is identically zero.  Thus $\Sigma$ is minimal
if and only if it is a critical point for the area functional.

Since a critical point is not necessarily a minimum the term
``minimal'' is misleading, but it is time-honored. A computation
shows that if $\Sigma$ is minimal, then
\begin{equation}
\left. \frac{d^2}{dt^2} \right|_{t=0}\Area (\Sigma_{t,\phi})
=-\int_{\Sigma}\phi\,L_{\Sigma}\phi  \, ,
\end{equation}
 where $L_{\Sigma}\phi
=\Delta_{\Sigma}\phi+|A|^2\phi+\Ric_M(\nn_{\Sigma},\nn_{\Sigma})\phi$
is the second variational (or Jacobi) operator. Here
$\Delta_{\Sigma}$ is the Laplacian on $\Sigma$,
$\Ric_M(\nn_{\Sigma},\nn_{\Sigma})$ is the Ricci curvature of $M$
in the direction of the unit normal to $\Sigma$, and $A$ is the
second fundamental form of $\Sigma$.  So $A$ is the covariant
derivative of the unit normal of $\Sigma$ and
$|A|^2=\kappa_1^2+\kappa_2^2$.

 For us, the key is that {\emph{the second variational operator is a
Schr\"odinger operator}} with potential
$V=|A|^2+\Ric(\nn_{\Sigma},\nn_{\Sigma})$.

A useful example to keep in mind is that of the catenoid. The
catenoid is the complete embedded minimal surface in $\RR^3$ that
is given by conformally embedding the flat $2$-dimensional
cylinder into $\RR^3$ by
\begin{equation}
    (\theta ,t) \to (-\cosh t \, \sin \theta , \cosh t \cos \theta , t ) \, .
\end{equation}
A calculation shows that pulling back the second variational
operator to the flat cylinder gives a rotationally symmetric
Schr\"odinger operator with potential
\begin{equation}
   V(\theta,t)= V(t) = 2 \, \cosh^{-2} (t) \, .
\end{equation}
Similarly, each  of the singly-periodic minimal surfaces known as
the Riemann examples is conformal to a flat cylinder with a
periodic set of punctures. Pulling back the second variational
operator to the flat cylinder gives a Schr\"odinger operator with
bounded potential.

A minimal surface $\Sigma$ is said to be stable if
\begin{equation}
\left. \frac{d^2}{dt^2} \right|_{t=0}\Area (\Sigma_{t,\phi})\geq 0
\,\,\,\,\,\,\,\,\,\,\,\text{ for all } \phi\in
C^{\infty}_0(\Sigma)\, .
\end{equation}
The {\it Morse index} of $\Sigma$ is the index of the critical
point $\Sigma$ for the area functional, that is, the number of
negative eigenvalues (counted with multiplicity) of the second
derivative of area; i.e., the number of negative eigenvalues of
$L$.{\footnote{By convention,  an eigenfunction $\phi$ with
eigenvalue $\lambda$ of $L$ is a solution of $L\,\phi+\lambda\,
\phi=0$.}}  Thus $\Sigma$ is {\it stable} if the index is zero. If
$\lambda=0$, then $\phi$ is said to be a Jacobi field.

\begin{figure}
\begin{center}
\input{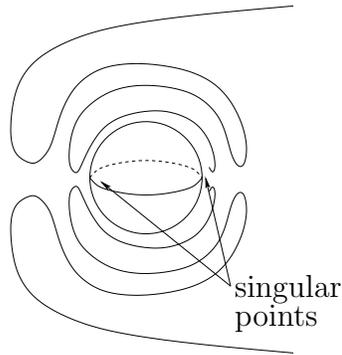}
\caption{One of the two possible singular laminations in half of a
neighborhood of a strictly stable $2$-sphere.  There are two
leaves.  Namely, the strictly stable $2$-sphere and half of a
cylinder. The cylinder accumulates towards the $2$-sphere through
catenoid type necks.  In fact, the lamination has two singular
points over which the necks accumulate.} \label{f:f1}
\end{center}
\end{figure}

\begin{figure}
\begin{center}
\input{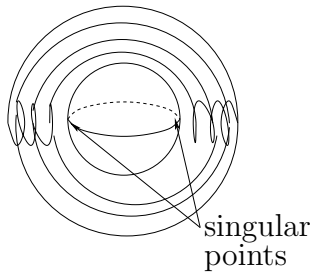}
\caption{One of the two possible singular laminations in half of a
neighborhood of a strictly stable $2$-sphere.  There are two
leaves.  Namely, the strictly stable $2$-sphere and half of a
cylinder. The cylinder accumulates towards the $2$-sphere and is
obtained by gluing together two oppositely oriented double spiral
staircases. Each   double spiral staircase winds tighter and
tighter as it approaches the $2$-sphere and, thus, never actually
reaches the $2$-sphere.} \label{f:f2}
\end{center}
\end{figure}

Suppose that $M^3$ is a fixed closed $3$-manifold with a
bumpy{\footnote{Bumpy means that no closed minimal surface
$\Sigma$ has $0$ as an eigenvalue of $L_{\Sigma}$ and the space of
such metrics is of Baire category by a result of B. White,
\cite{W1}.}} metric with positive scalar curvature and let
$\Sigma_i$ be a sequence without repeats, i.e., with $\Sigma_i\ne
\Sigma_j$ for $i\ne j$, of embedded minimal surfaces of a given
fixed genus. After possibly passing to a subsequence one expects
that it converges to a singular lamination{\footnote{A lamination
is a foliation except for that it is not assumed to foliate the
entire manifold.}} that looks like one of the two illustrated
below:

\begin{quote}
One expects that any singular limit lamination has only finitely
many leaves.  Each closed leaf is a strictly stable $2$-sphere.
Each non-compact leaf has only finitely many ends and each end
accumulates around exactly one of the closed leaves.  The
accumulation looks almost exactly as in either Figure \ref{f:f1}
or Figure \ref{f:f2}.
\end{quote}

Indeed the lamination in Figure \ref{f:f1} can happen as a limit
of fixed genus embedded minimal surfaces in a $3$-manifold, see
\cite{CD} (even in a $3$-manifold with positive scalar curvature);
cf. also with B. White, \cite{W2} and M. Calle and D. Lee,
\cite{CaL}.

  For us, the key is that (see Section \ref{s:fromgeom}):

  \vskip1mm \noindent
{\emph{Each non-compact leaf is conformally a
Riemann surface with finitely many cylindrical ends}},
 \vskip1mm

and under the conformal change,
 \vskip1mm \noindent
{\emph{the second variational operator becomes a Schr\"odinger
operator with Lipschitz bounded potential}.}
 \vskip1mm

One would like to understand the moduli space of such non-compact
minimal surfaces.  Infinitesimally, the space of nearby
non-compact minimal surfaces with finitely many ends, each as
Figure \ref{f:f1} or \ref{f:f2}, are solutions of the second
variational equation on the initial surface.  Thus, we are led to
analyze the solutions of this Schr\"odinger equation.

\subsection{Schr\"odinger operators on $\RR^{n+1}$}

 Theorem \ref{t:estsc0} implies a three circles
inequality, and a corresponding strong unique continuation
theorem, for a Euclidean operator
\begin{equation}
    L = \Delta_{\RR^{n+1}}   - (n-1)\, |x|^{-1} \,
\partial_{|x|} + V (x)  \, ,
\end{equation}
where $\partial_{|x|}$ is the radial derivative and the potential
$V(x)$
  satisfies
\begin{equation}
 \label{e:dec1}
    |V(x)| \leq C \, |x|^{-2}   \, .
\end{equation}
This unique continuation  does not follow from the well-known
sharp result for potentials $V \in L^{ \frac{n+1}{2} }
(\RR^{n+1})$ of Jerison and Kenig, \cite{JK}.  It also does not
follow from the unique continuation result of Garofalo and Lin,
\cite{GL}, which holds when $|x|^2 \, |V(x)|$ goes to zero at a
definite rate. To our knowledge, the sharpest unique continuation
results for Euclidean operators of this general form are given in
Pan and Wolff, \cite{PW}. In that paper, they consider operators
$\Delta_{\RR^{n+1}}   + W(x) \cdot \nabla_{\RR^{n+1}} + V (x)$,
where $V$ satisfies \eqr{e:dec1} for some constant and $W$
satisfies $|x| \, |W(x)| \leq C_0$ for a fixed small constant
$C_0$.

To see why Theorem \ref{t:estsc0} applies to the operator $L$, it
will be convenient to work in ``exponential polar coordinates'' $(
\theta = x/|x| , \, t = \log |x|) \in \SS^n \times \RR $.  In
these coordinates, the chain rule gives
\begin{align}
    \partial_{|x|} &= \e^{-t} \,
\partial_t \, , \\
    \partial_{|x|}^2 &= \e^{-2t} \,
\left( \partial_t^2 - \partial_t \right) \, .
\end{align}
Using this, we can rewrite the
   Euclidean Laplacian
$\Delta_{\RR^{n+1}}$   as
\begin{equation}
    \Delta_{\RR^{n+1}} =   \partial_{|x|}^2 + \frac{n}{|x|} \, \partial_{|x|}  +
    |x|^{-2} \, \Delta_{\SS^n} =
     \e^{-2t} \, \Delta_{\SS^n \times \RR} + \e^{-2t} \, (n-1) \, \partial_t  \,
     .
\end{equation}
 Therefore, the
Euclidean operator $L$ can be written
\begin{equation}
    \e^{2t} \, L  =  \Delta_{\SS^n \times \RR}      +
    \e^{2t} \, V(\e^t \, \theta)
     \, .
\end{equation}
In particular, if $V$   satisfies \eqr{e:dec1}, then the operator
$\e^{2t} \, L $ can be written as $\Delta_{\SS^n \times \RR} +
\tilde{V}$, where the potential $\tilde{V}$  is bounded. It
follows that Theorem \ref{t:estsc0} applies to an operator $L$
satisfying \eqr{e:dec1}.

\subsection{Outline of the paper}
  In Section
\ref{s:fourier},  on the   half-cylinder $N \times [0,\infty)$
with coordinates $(\theta , t)$, we introduce notation for the
  Fourier coefficients (or spectral projections) of a function
$f(\theta, t)$ on each cross-section $t=$ constant.

(In Appendix \ref{s:sone}, we specialize  to the case of a
cylinder $N \times \RR$ and a rotationally symmetric potential
$V(\theta,t) = V(t)$. This is meant only to explain some of the
ideas in a simple case and the results will not be used elsewhere.
Given a solution $u$ of the Schr\"odinger equation, an easy
calculation shows that the Fourier coefficients of $u$ satisfy an
ODE as a function of $t$. It follows from a Riccati comparison
argument, that any sufficiently high Fourier coefficient of $u$
grows exponentially at either plus infinity or minus infinity. In
particular, if the solution $u$ vanishes at both plus and minus
infinity, then all sufficiently high Fourier coefficients vanish.
It follows from this that the space $H_0$ is finite dimensional.
Similarly for $H_{\alpha}$ when $\alpha > 0$.)

In Section \ref{s:sthree}, we prove the three circles theorem for
Lipschitz  potentials, i.e., Theorem \ref{t:ests}.  Unlike the
case of   rotationally symmetric potentials, the individual
Fourier coefficients will no longer satisfy a useful ODE, but we
will still be able to show that the simultaneous projection of a
solution $u$ onto all sufficiently large Fourier eigenspaces
satisfies a useful differential inequality. To give a feel for the
proof, we will now outline the argument. For each $t\in [0,T]$,
let $[u]_j(t)$ be the $j$-th Fourier coefficient of a solution $u$
restricted to the $t$-th slice. Define functions of $t$ by
$\cL_m=\sum_{j=0}^{m-1}\left[([u]_j')^2+(1+\lambda_j)[u]_j^2\right]$
and
$\cH_m=\sum_{j=m}^{\infty}\left[([u]_j')^2+(1+\lambda_j)[u]_j^2\right]$
and note that the sum of the two is  the Sobolev norm.   A
computation   shows that they satisfy the two differential
inequalities:  $\cH_m''\geq (4\lambda_m-C)\,\cH_m-C\cL_m$ and
$\cL_m''\leq (4\lambda_{m-1}+C)\,\cL_m+C\cH_m$, for some constant
$C$ depending only on the Lipschitz norm of the potential and in
particular not on $m$.  Subtracting the second inequality from the
first and using the spectral gap yields that $[\cH_m-\cL_m]''\geq
(4\lambda_{m-1}+2\kappa)\, [\cH_m-\cL_m]$ for some positive
constant $\kappa$ and $m$ sufficiently large.  We then use this
differential inequality and the maximum principle applied to the
function $f(t)=\e^{-\alpha t} [\cH_m-\cL_m]$, where $\alpha$ is
the logarithmic growth rate of the Sobolev norm from $t=0$ to
$t=T$, to conclude that $\cH_m (t)$ is bounded in terms of
$\e^{\alpha t}\,I(0)+\cL_m(t)$.  Inserting this back into the
first order differential inequality that $\cL_m$ satisfies easily
gives a bound for $\cL_m(t)$ (and hence for $\cH_m(t)$ and $I(t)$)
in terms of $\e^{\alpha \, t}\,I(0)$.  Unravelling it all yields
the desired three circles inequality, i.e., Theorem \ref{t:ests}.
In Section \ref{s:ne}, we prove a three circles inequality when
the potential $V$ is bounded, i.e., Theorem \ref{t:estsc0}.

 Using the results of Section \ref{s:ne},
we will show in Section \ref{s:sfour} that the space $H_{\alpha}$
is finite dimensional  on a manifold with finitely many ends, each
of which is isometric to a half-cylinder. In Section
\ref{s:dense}, we show that the space $H_0$ is zero dimensional
for a dense set of
  potentials. Subsection \ref{s:open} gives an example where
the set of potentials with $H_0=\{0\}$ is not open.

In Section \ref{s:ssix}, we prove a uniformization theorem that
allows us to reduce the general case of surfaces with cylindrical
ends to the case where the ends are isometric to  flat
half-cylinders. Together with the results of Sections \ref{s:ne},
\ref{s:sfour}, and \ref{s:dense}, this proves the main theorem.

\section{Examples from geometry}        \label{s:fromgeom}

In this section, we will show that for each non-compact leaf of
the singular minimal lamination constructed  in \cite{CD} (see
  Figure \ref{f:f1}) our main results, Theorem  \ref{t:main} and
Theorem \ref{t:ests}, apply. Namely, we show the following
proposition.

\begin{Pro}     \label{p:cd}
Each non-compact leaf of the singular minimal lamination
constructed  in \cite{CD} is conformally a Riemann surface with
finitely many cylindrical ends and, after this conformal change,
the second variational operator becomes a Schr\"odinger operator
with  bounded potential.  In fact, the conformal change of metric
that we give below will directly make each end isometric to a flat
half-cylinder.
\end{Pro}

Let $M^3$ be a closed $3$-manifold with a Riemannian metric $g$
and $\cL$ a minimal lamination consisting of finitely many leaves,
as constructed in \cite{CD}.  Each compact leaf is a strictly
stable $2$-sphere and each non-compact leaf has only finitely many
ends, each end, a half infinite cylinder spiralling into one of
the strictly stable $2$-spheres as in Figure \ref{f:f1}.      To
prove the proposition, it is enough to show that we can
conformally change the metric on each end $\Sigma$  to make it a
flat cylinder and then show that, in this conformally changed
metric, the second variational operator becomes a Schr\"odinger
operator with bounded potential.

In this example, we can parametrize a neighborhood of the strictly
stable $2$-sphere by $\SS^2\times (-\varepsilon, \varepsilon)$ and
on $\SS^2$ use spherical coordinates $(\phi,\theta)$; $r\in
(-\varepsilon,\varepsilon)$ denotes the (signed) distance to the
strictly stable $2$-sphere.  In these coordinates the metric $g$
takes the form
\begin{equation}\label{e:metric}
dr^2 + \mu^2 (r) (d\phi^2 + \sin^2 \phi\, d\theta^2)\,
\end{equation}
(see equation (2) in \cite{CD}). Moreover, $\mu$ is a smooth
function with $\mu (0)=1$, $\mu' (0) = 0$ and $\mu'' > 0$.

\begin{figure}
\begin{center}
\input{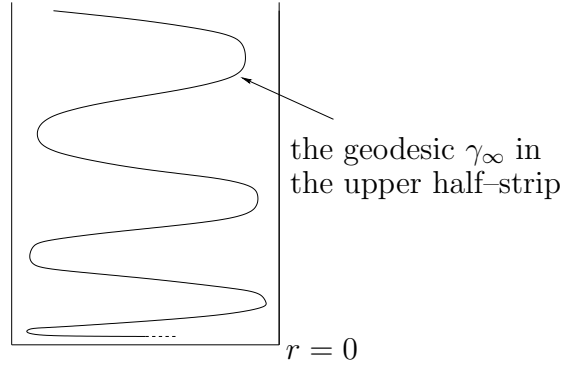}
\caption{The projection of the half-infinite cylinder $\Sigma$ in
$M$ is an infinite geodesic $\gamma_{\infty}$ in the upper
half-strip with the degenerate metric \eqr{e:degmet}.}
\label{f:f4}
\end{center}
\end{figure}

The minimal half-cylinder $\Sigma$ is $\SS^1$-invariant, i.e., it
is the preimage of a curve $\gamma_{\infty}$ on the strip $[0,
\pi]\times (-\varepsilon, \varepsilon)$ under the projection map
\begin{equation}    \label{e:pi}
(\phi, \theta, r) \;\mapsto\; (\phi, r)\, .
\end{equation}
As first remarked by Hsiang and Lawson in \cite{HsLa} (cf. with
section 2 of \cite{CD}), since $\Sigma$ is a critical point for
the area functional, $\gamma_{\infty}$ is a critical point for the
functional
\begin{equation}
  F (\gamma_{\infty})=\int_{\gamma_{\infty}}\, \textrm{length} ( \SS^1 \times \{\gamma_{\infty}
  (t)\})=\int_{\gamma_{\infty}}\, 2\pi\, \mu (r(t))\, \sin (\phi (t))\, .
\end{equation}
Therefore, $\gamma_{\infty}$ is an infinite geodesic for the
degenerate metric
\begin{equation}    \label{e:degmet}
\mu^2 (r) \sin^2 \phi\, (dr^2 + \mu^2 (r) d\phi^2)\, ,
\end{equation}
accumulating towards the geodesic segment $\{r=0\}$; see Figure
\ref{f:f4}.

If we assume that $t\mapsto (\phi (t), r(t))$ is the
parameterization of $\gamma_{\infty}$ by arc-length ($t>0$) in the
degenerate metric \eqr{e:degmet}, then
\begin{equation}\label{e:arclength}
\left(\frac{dr}{dt}\right)^2 + \mu^2 (r)
\left(\frac{d\phi}{dt}\right)^2 \;=\; \mu^{-2} (r) \sin^{-2} \phi\, .
\end{equation}
Therefore, if we parameterize $\Sigma$ by $(t, \theta)\mapsto
(\phi (t), r (t), \theta)$, the induced metric on $\Sigma$ is
\begin{eqnarray}
d\sigma^2 &=& \left[\left(\frac{dr}{dt}\right)^2 + \mu^2 (r)
\left(\frac{d\phi}{dt}\right)^2\right] dt^2 + \mu^2 (r) \sin^2 \phi\,
d\theta^2\nonumber\\
&=& \mu^{-2} (r) \sin^{-2} \phi\, dt^2 + \mu^2 (r) \sin^2
\phi\, d\theta^2\, .
\end{eqnarray}
Let $\tau$ be a new parameterization of $\gamma_{\infty}$, so that
\begin{equation}\label{e:reparam}
\frac{dt}{d\tau} = \mu^2 (r(t)) \sin^2 (\phi (t))\, .
\end{equation}
It follows that in the coordinates $(\tau, \theta)$ the metric on
$\Sigma$ takes the form
\begin{equation}\label{e:conformal}
\mu^2 (r(\tau)) \sin^2 \phi (\tau) \bigl(d\tau^2 + d\theta^2)\, ,
\end{equation}
i.e., $(\tau, \theta)$ is a conformal parameterization with
conformal factor $h = \mu (r) \sin \phi$.

To complete the proof of Proposition \ref{p:cd}, it only remains
to show that the second variational operator $L=\Delta_{d\sigma^2}
+ (|A|^2 + \Ric_M (\nn_\Sigma, \nn_\Sigma))$ on $\Sigma$ has the
same kernel as a Schr\"odinger operator $\tilde{L}$  with bounded
potential in the conformally changed metric $ds^2 = h^{-2}
d\sigma^2$.  We will do this in the next lemma for the operator
$\tilde{L} = h^2 \, L$.

\begin{Lem} \label{l:rescaledJacobi}
In the conformally changed metric $ds^2 = h^{-2} d\sigma^2$ (i.e.,
the flat metric on the half cylinder), the operator $\tilde{L} =
h^2 \, L $
 is a Schr\"odinger
operator with bounded potential.
\end{Lem}

\begin{proof} Since $ds^2 = h^{-2} d\sigma^2$,
we have $\Delta_{ds^2} = h^2 \Delta_{d\sigma^2}$ and, therefore,
$\tilde{L} = \Delta_{ds^2} + h^2 (|A|^2 + \Ric_M (\nn_\Sigma,
\nn_\Sigma))$ is a Schr\"odinger operator in the metric $ds^2$;
cf. \eqr{e:confH0}. Since both $\Ric_M (\nn_\Sigma, \nn_\Sigma)$
and $h$ are  bounded, to prove the proposition, it suffices to
show that $h^2 |A|^2$ is bounded.

In what follows, we will denote by $\dot{r}$ and $\dot{\phi}$ the
derivatives $\textstyle{\frac{dr}{d\tau}}$ and
$\textstyle{\frac{d\phi}{d\tau}}$, respectively. According to
\eqref{e:arclength} and \eqref{e:reparam} we have
\begin{equation}\label{e:speed}
(\dot{r})^2 + \mu^2 (r) (\dot{\phi})^2 \;=\; \mu^2 (r) \sin^2 \phi = h^2\,
.
\end{equation}
Set
\begin{equation}
A_{\theta \theta} \;=\; - g \left(\nn_\Sigma ,
\nabla_{\partial_\theta}
\partial_\theta\right) \, ,
\qquad A_{\tau\theta} \;=\; - g \left(\nn_\Sigma ,
\nabla_{\partial_\tau}
\partial_\theta\right) \, ,
\qquad A_{\tau \tau} \;=\; - g \left(\nn_\Sigma ,
\nabla_{\partial_\tau}
\partial_\tau\right)\, .
\end{equation}
By minimality, $A_{\theta\theta} = - A_{\tau \tau}$, and hence
\begin{equation}
h^2 |A|^2 \;=\; h^2 \left[ h^{-4} (A_{\theta\theta}^2 +
A_{\tau\tau}^2 + 2 A_{\tau\theta}^2)\right] \;=\; 2 \, h^{-2} \,
\left[ A_{\theta\theta}^2 +   A_{\tau \theta}^2 \right] \, .
\end{equation}
It can be readily checked that the normal $\nn= \nn_\Sigma$ is
given by
\begin{align}\label{e:normal}
\nn \;&=\; (\dot{r} \, \partial_\phi - \mu^2 (r)\, \dot{\phi} \,
\partial_r )/ [(\mu^2 (r)\, (\dot{r})^2 + \mu^4 (r)\, (\dot{\phi})^2)^{1/2}]  \\ \,
\;&=\; \mu^{-2} (r)\,  \sin^{-1} \, \phi\, (\dot{r}\,
\partial_\phi - \mu^2 (r)\,  \dot{\phi}\,  \partial_r)\, . \notag
\end{align}
Moreover, since also $\partial_{\tau}$ lies in the linear span of
$\partial_{r}$ and $\partial_{\phi}$ and the level sets $\theta$
equal constant are totally geodesic in the metric $g$, it follows
easily that $A_{\tau \theta} = 0$.  Finally,
\begin{eqnarray}
A_{\theta\theta} &=& - \mu^{-2} (r) \sin^{-1} \phi \left[ \dot{r} g
\left(\partial_\phi, \nabla_{\partial_\theta}
\partial_\theta\right)
- \mu^2 (r) \dot{\phi} g
\left(\partial_r, \nabla_{\partial_\theta} \partial_\theta\right)\right] \nonumber\\
&=& - \mu^{-2} (r) \sin^{-1} \phi \left[ -
\frac{\dot{r}}{2}\partial_\phi (g (\partial_\theta,
\partial_\theta))
+ \frac{\mu^2 (r) \dot{\phi}}{2}
\partial_r (g(\partial_\theta, \partial_\theta))\right]\nonumber\\
&=& - \mu^{-2} (r) \sin^{-1} \phi \left( - \mu^2 (r) \dot{r} \sin \phi \cos
\phi + \mu^3 (r) \dot{\phi} \mu' (r) \sin^2 \phi\right)\nonumber\\
&=& \dot{r} \cos
\phi - \mu (r) \mu' (r) \dot{\phi} \sin \phi \label{e:thetatheta}\, ,
\end{eqnarray}
and
\begin{align}
h^4 |A|^2 &\;=\; 2 ( \dot{r}\cos \phi - \mu (r) \mu' (r) \dot{\phi} \sin
\phi  )^2 \leq 4 \left[ ( \dot{r})^2 \cos^2 \phi +
 h^2 (\mu' (r))^2 (\dot{\phi})^2
\right]  \notag \\
& \leq 4 h^2 \, \left[ 1 +  (\mu' (r))^2  \right] \, , \label{e:partA}
\end{align}
where the last inequality follows from \eqref{e:speed}. The
desired bound on $h^2 |A|^2$ now follows.
\end{proof}

Next, consider the Jacobi
fields generated by sequences of spiralling cylinders
$\{\Sigma_n\}$ of the form above. Then these Jacobi fields grow at most exponentially
in $\tau$.

\begin{Def}\label{d:setofcylinders}
We let $M$ be the Riemannian manifold $\SS^2\times (-\varepsilon,
\varepsilon)$ with the metric $g$ of \eqref{e:metric}. Any
isometry $\Phi$ of the standard $\SS^2$ can be extended to an
isometry of $M$ in an obvious way, i.e. by mapping $(z,r) \in
\SS^2 \times (-\varepsilon, \varepsilon)$ to $(\Phi (z), r)$. We
denote by $G$ the set of such isometries. Finally, we denote by
$\mathcal{S}$ the set of minimal $\SS^1$--invariant cylinders
spiraling into $\SS^2\times \{0\}$. That is, $\Gamma$ is an
element of $\mathcal{S}$ if and only if there exist a minimal
cylinder $\Sigma$ and a $\Phi\in G$ such that $\Gamma = \Phi
(\Sigma)$ and $\Sigma$ is the lifting of a curve $\gamma_\infty$
under the projection map \eqref{e:pi}.
\end{Def}

Loosely speaking, the set of Jacobi fields generated by sequences
of elements of $\mathcal{S}$ gives the tangent space to
$\mathcal{S}$.  More precisely, let $\{\Sigma_k\}$ be a sequence
of elements of $\mathcal{S}$ that converges to   $\Sigma\in
\mathcal{S}$. Consider a sequence of increasing compact domains
$\Omega_0\subset \Omega_1\subset \ldots \subset \Sigma$ exhausting
$\Sigma$. For each $i$ we select $\varepsilon_i$ sufficiently
small and we consider the portion $T_i$ of the
$\varepsilon_i$--tubular neighborhood of $\Sigma$ which is ``lying
above'' $\Omega_i$, that is
\begin{equation}\label{e:Ti}
T_i\;=\;\bigl\{ x + \exp_x (s \nn_\Sigma (x))\, |\, x\in \Omega_i, s\in
(-\varepsilon_i, \varepsilon_i)\bigr\}\, .
\end{equation}
Let $i$ be given. By the standard regularity theory for minimal surfaces,
for $k$ large enough $\Sigma_k \cap T_i$ is a graph over $\Omega_i$, i.e.
\begin{equation}\label{e:graph}
\Sigma_k\cap T_i \;=\;
\bigl\{x + \exp_x (u_k (x) \nn_\Sigma (x))\, |\, x\in \Omega_i\}
\end{equation}
for some smooth function $u_k$.

We normalize $u_k$ to $f_k = u_k/\|u_k\|_{L^2 (\Omega_0)}$. Then,
a subsequence  converges to a nontrivial smooth function $f$ on
$\Sigma$ solving $\tilde{L} f =0$, where $\tilde{L}$ is the
operator of Lemma \ref{l:rescaledJacobi}. We denote by
$T_\Sigma\mathcal{S}$ the space of functions $c f$, where $f$ is
generated with the procedure above and $c$ is a real number.

\begin{Lem}\label{l:Halpha-bound}
There exists a constant $\alpha$ such that the following holds.
Consider any $\Sigma\in \mathcal{S}$ with the rescaled flat metric
$ds^2$ as in Lemma \ref{l:rescaledJacobi}. Then $T_\Sigma
\mathcal{S}\subset H_\alpha (\Sigma)$ for some $\alpha \geq 0$.
\end{Lem}
\begin{proof} Without loss of generality we can assume that $\Sigma$ is the lifting of
a curve $\gamma_\infty$ through the projection \eqref{e:pi}. Therefore,
we use on $\Sigma$ the coordinates $(\theta, \tau)$ introduced in
Lemma \ref{l:rescaledJacobi}.

Let $\mathcal{G}$ be the Lie Algebra generating $G$ and define the
linear space $V= \{g (X, \nn_\Sigma)|X\in \mathcal{G}\}$. Clearly,
$V$ is a space of bounded smooth functions on $\Sigma$. Moreover,
$V$ gives the Jacobi fields generated by minimal surfaces of the
form $\{\Phi_n (\Sigma)\}$ for sequences $\{\Phi_n\}\subset G$
converging to the identity. Therefore, any element $f\in T_\Sigma
\mathcal{S}$ can be written as $v+w$, where $v$ belongs to $V$ and
$w$ is a function of $T_\Sigma \mathcal{S}$ independent of the
variable $\theta$. We sketch a proof of this fact for the reader's
convenience. Let $f$ be a nontrivial element of $T_\Sigma
\mathcal{S}$  generated by a sequence of $\SS^1$--invariant
minimal cylinders $\Sigma_k$ as above. Then $\Sigma_k = \Phi_k
(\Gamma_k)$, where
\begin{itemize}
\item $\{\Phi_k\}$ is a sequence of isometries converging to the identity;
\item $\Gamma_k$ are liftings of curves $\gamma_k$ through the
projection \eqref{e:pi}.
\end{itemize}
Let $i$ be a given natural number.
For $k$ sufficiently large, $\Sigma_k\cap T_i$ has the form
\begin{equation}\label{e:graph2}
\Sigma_k\cap T_i \;=\;
\bigl\{x + \exp_x (u_k (x) \nn_\Sigma (x))\, |\, x\in \Omega_i\}
\end{equation}
and $u_k/\|u_k\|_{L^2 (\Omega_0)}$ converges to $f$.

On the other hand, by the standard theory of minimal surfaces,
the Hausdorff distance between $\Gamma_k\cap T_i$
and $\Sigma\cap T_i$ and $\Phi_k (\Gamma_k)\cap T_i$ and $\Gamma_k\cap T_i$
converge to $0$.
Hence, for $k$ sufficiently large, $\Gamma_k\cap T_i$ is a graph over
$\Omega_i$ and $\Phi (\Gamma_k)\cap T_i$ is a graph over $\Gamma_k$.
Thus we can find functions $v_k$ and $w_k$ such that
\begin{equation}\label{e:graph3}
T_i\cap \Gamma_k \;=\; \bigl\{ x + \exp_x (w_k \nn_\Sigma (x))\, |\,
x\in \Omega_i\bigr\}
\end{equation}
\begin{equation}\label{e:graph4}
T_i\cap \Phi_k (\Gamma_k)\;=\; \bigl\{x+ \exp_x (w_k \nn_\Sigma (x))
+ \exp_{x+\exp_x (w_k \nn_\Sigma (x))} (v_k \nn_{\Gamma_k} (x))\bigr\}\, .
\end{equation}
Note that $w_k$ is a function independent of $\theta$. Moreover, up to subsequences
we can assume that $w_k/\|w_k\|_{L^2 (\Omega_0)}$ converges to a function $w$.
Such a $w$ belongs to $T_\Sigma \mathcal{S}$ and depends only on the variable $\tau$.
Finally, up to subsequences, we can assume that $v_k/\|v_k\|_{L^2 (\Omega_0)}$
converges to an element $v$ of $V$.

By the theory of minimal surfaces, the Hausdorff distances between
$\Gamma_k\cap T_i$ and $\Sigma\cap T_i$ and $\Phi_k (\Gamma_k)\cap
T_i$ and $\Gamma_k\cap T_i$ are controlled by $\|u_k\|_{L^2
(\Omega_0)}$. Moreover, $u_k = w_k + v_k + o (\|u_k\|_{L^2
(\Omega_0)})$. Since $f$ is the limit of $u_k/\|u_k\|_{L^2
(\Omega_0)}$, $f$ must be a linear combination of $v$ and $w$.

Having shown the desired decomposition for any element of
$T_\Sigma \mathcal{S}$, since $V$ is a space of bounded functions,
it suffices to show the existence of $\alpha\geq 0$ such that
every function $f\in T_\Sigma \mathcal{S}$ independent of $\theta$
belongs to $H_\alpha (\Sigma)$. For any such $f$ we have, by Lemma
\ref{l:rescaledJacobi}, $f'' (\tau) = -V (\tau) f (\tau)$. Since
$V$ is bounded, this gives the inequality
\begin{equation}
|f''|\;\leq\; \|V\|_\infty |f| \;=\; a |f|\, .
\end{equation}
Consider the nonnegative locally Lipschitz function $g (\tau) =
|f' (\tau)| + |f(\tau)|$ and set $\alpha = \max \{a,1\}$. Then
\begin{equation}
g'\;\leq\; |f''| + |f'| \leq a |f| + |f'| \leq \alpha g\, .
\end{equation}
Hence, from Gronwall's inequality, we get $|f(\tau)|\leq g (\tau)
\leq g(0) e^{\alpha \tau}$ for $\tau \geq 0$, which is the desired
bound.
\end{proof}

\section{The spectral projection on a closed manifold} \label{s:fourier}

Suppose now that $N^n$ is an $n$-dimensional closed Riemannian
manifold and $\Delta_N$ is the Laplacian on $N$.  We will
generally use $\theta$ as a parameter on $N$.  Fix an
$L^2(N)$-orthonormal basis of $\Delta_N$ eigenfunctions $\phi_0 ,
\phi_1 , \dots$ with eigenvalues $0 = \lambda_0 < \lambda_1 \leq
\dots$, so that
\begin{equation}
    \Delta_N \phi_j = - \lambda_j \, \phi_j \, .
\end{equation}
Given an arbitrary $L^2$ function $f$ on $N$, we will let $[f]_j$
denote the inner product of $f$ with $\phi_j$
\begin{equation}
    [f]_j = \int_N f(\theta) \, \phi_j (\theta) \, d\theta \, .
\end{equation}
In analogy to the special case where $N=\SS^1$ (see below), we
will often refer to this as the $j$-th Fourier coefficient, or
$j$-th spectral projection. It follows that
\begin{equation}
    f(\theta) = \sum_{j=0}^{\infty} [f]_j  \, \phi_j (\theta)   \, .
\end{equation}

It will often be important to understand how the Fourier
coefficients of a function $f(\theta , t)$  on the half-cylinder
$N \times [0,\infty)$  vary as a function of $t$. To keep track of
these coefficients,  we  define  $[f]_j(t)$ by
\begin{equation}
    [f]_j(t) = \int_N f(\theta, t) \, \phi_j (\theta) \, d\theta \, .
\end{equation}

\subsection{The Fourier coefficients on a half-cylinder}

The simplest example of   spectral projection is when $N$ is the
unit circle $\SS^1$ with the standard orthonormal basis of
eigenfunctions
\begin{equation}
    \phi_0 =  \frac{1}{\sqrt{2\pi}} \, , \,   \left\{ \phi_{2k+1} = \frac{1}{\sqrt{\pi}}
\, \sin (k\theta)  \right\}_{k \geq 0} \, , \, \left\{ \phi_{2k} =
\frac{1}{\sqrt{\pi}} \, \cos (k\theta) \right\}_{k \geq 1} \, ,
\end{equation}
with eigenvalues $\lambda_0 = 0$ and $\lambda_{2k+1}= \lambda_{2k}
= k^2$.  In this case, the $[f]_{j}$'s are the Fourier
coefficients of the function $f$.

\section{General Lipschitz bounded potentials: The three circles inequality}
\label{s:sthree}

Throughout this section, $u$ will be a solution of
\begin{equation}    \label{e:equn}
    \Delta u = - Vu
\end{equation}
 on a product $N \times [0,t]$, where the potential $V$ will be
Lipschitz, but is no longer assumed to be rotationally symmetric.

The results of Appendix \ref{s:sone} in the rotationally symmetric
case where $V=V(t)$ were stated on an entire cylinder, but the
corresponding results for the half-cylinder
 motivate the general results of this section.  Namely, the ODE
   \eqr{e:uk} for the Fourier coefficients of $u$ as a function of
   $t$ implies that the $j$-th Fourier coefficient
    must either grow or decay exponentially if $\lambda_j > \sup V$.
    This same analysis holds
   even on a half-cylinder when $V$ is
   rotationally symmetric.  We will prove similar results in this
   section for a general bounded potential $V=V(\theta,t)$,
but things are more complicated  since multiplication by
$V(\theta,t)$ does not preserve  the   eigenspaces of $\Delta_N$
(i.e., $\phi_j (\theta)$).{\footnote{The reason that the ODE
\eqr{e:uk} is so simple
   is that the $j$-th Fourier coefficient of $V(t) u(\theta,t)$
is just $V(t)$ times
   the $j$-th Fourier coefficient of $u(\theta,t)$.}}

The main result of this section is Theorem \ref{t:ests0} below
that shows a three circles inequality for the Sobolev norm of a
solution of a Schr\"odinger equation on a product $N\times [0,T]$
where $N$ has the required spectral gaps.  This will give Theorem
\ref{t:ests} in the special case where $N$ is a round sphere or a
Zoll surface.  See the upshot to Section \ref{s:sthree} in the
introduction for an overview of the proof.

\subsection{The Fourier coefficients of $u$}
As in the rotationally symmetric case, it will be important to
understand how the Fourier coefficients $ [u]_{j}(t)$   and its
derivatives grow or decay as a function of $t$.

The next lemma gives the ODE's that govern how the Fourier
coefficients  $ [u]_{j}(t)$  grow or decay as   functions of
$t$; cf. the similar ODE
 \eqr{e:uk} in the rotationally symmetric case.

\begin{Lem}     \label{l:diff}
The Fourier coefficients $  [u]_{j}(t)$ satisfy
\begin{align}
      [u]_{j}'(t)  &= [u_t]_{j}(t) \, , \\
    [u]_j''(t) &=  \lambda_j \, [u]_j (t) - [Vu]_j (t)
    \, ,
    \label{e:dig} \\
    [u]_j'''(t) &= \lambda_j \, [u]_j' (t) - [\partial_t (Vu)]_{j} (t) \,
    .
\end{align}
\end{Lem}

\begin{proof}
 Differentiating $[u]_j(t)$ immediately gives the first claim.  To
 get the second claim, first differentiate again to get
\begin{equation}
    [u]_j''(t) =   \int_N u_{tt} (\theta,t) \, \phi_j  (\theta) \,
    d\theta \, .
\end{equation}
Next, bring in the equation $u_{tt} = - \Delta_N u - V u $
 and integrate by parts twice to get
\begin{align}
    [u]_j''(t) &= -   \int_{N \times \{ t \} }  \phi_j   \, \Delta_N u    \,   \,
    d\theta -   \int_{N \times \{ t \} }  V  \, u    \, \phi_j   \,
    d\theta    \notag \\
    &= \lambda_j [u]_j(t) - [Vu]_j (t)
\, . \label{e:pluig}
\end{align}
Differentiating again gives
\begin{equation}
    [u]_j'''(t) =
\lambda_j \, [u]_j' (t) - [\partial_t (Vu)]_{j} (t)   \, .
\end{equation}
\end{proof}

As mentioned above,  \eqr{e:dig} implies exponential growth (or
decay) of $[u]_j(t)$ when $V$ is rotationally symmetric and
$\lambda_j > \sup V$. However, this is not the case for a general
bounded $V$ since the ``error term'' $[Vu]_{j} (t)$ need not be
bounded by $[u]_j(t)$.
 We will  get around this in the next subsection by considering all
 of the $[u]_j$'s above a fixed value at the same time.
To get a well-defined quantity when we do this, we will have to
sum the squares of the $[u]_j$'s.  Unfortunately,  the quantity
$[u]_j^2$ does not satisfy as nice of an ODE, so we will have to
consider a slightly different quantity.
 To see why,
observe that when $V=0$,  then
\begin{equation}    \label{e:last3}
      \partial_t^2 \left[  ([u]_j')^2 + \lambda_j \, [u]_j^2 \right]
=  4 \lambda_j \,
      \left[  ([u]_j')^2 + \lambda_j \, [u]_j^2 \right]  \, .
\end{equation}
Equation \eqr{e:last3} suggests   looking at the quantity
$([u]_j')^2 + \lambda_j \, [u]_j^2$, but it will be more
convenient to look at the slightly different quantity $([u]_j')^2
+ (1+\lambda_j) \, [u]_j^2$. This is because $([u]_j')^2 +
\lambda_j \, [u]_j^2$ is a piece of the $L^2$ norm of   $\nabla
u$, but $([u]_j')^2 + (1+\lambda_j) \, [u]_j^2$ also includes part
of the $L^2$ norm of $u$ and, hence, corresponds to the full
$W^{1,2}$ norm of $u$; see equation \eqr{e:standardf} below.

\begin{Lem}     \label{l:diff2}
The quantity $\left[  ([u]_j')^2 + (1+\lambda_j) \, [u]_j^2
\right]$ satisfies the ODE's
\begin{align}
\partial_t \left[  ([u]_j')^2 +
(1+\lambda_j) \, [u]_j^2 \right] &=  (4\, \lambda_j + 2) \, [u]_j
\, [u]_j' - 2   [u]_j'\, [Vu]_j \, , \label{e:last4a}
         \\
      \partial_t^2 \left[  ([u]_j')^2 +
(1+\lambda_j) \, [u]_j^2 \right] &=  (4 \lambda_j + 2) \,
      \left[  ([u]_j')^2 + (1+\lambda_j) \, [u]_j^2 \right] - (4\lambda_j + 2) \, [u]_j^2   \label{e:last4} \\
    &\quad \quad - (6 \lambda_j +2)\, [u]_j [Vu]_{j} + 2 \, [Vu]_{j}^2   - 2\, [u]_j' \,
    [\partial_t (Vu)]_j   \, .
\notag
\end{align}
\end{Lem}

\begin{proof}
Using Lemma \ref{l:diff}, we get
\begin{align}
    \frac{1}{2} \, ([u]_j^2)'  &= [u]_j   \, [u]_j'  \, ,  \label{e:last0a} \\
    \frac{1}{2} \, ([u]_j^2)''  &=   \lambda_j \, [u]_j^2 +
    ([u]_j')^2 - [u]_j \, [Vu]_j \label{e:last0}    \, .
\end{align}
  Similarly, differentiating $([u]_j')^2$
gives
\begin{align}    \label{e:last1}
    \frac{1}{2} \, \partial_t ([u]_j')^2 &=   \lambda_j \,
    [u]_j \, [u]_j'  - [u]_j' \, [Vu]_j   \, ,  \\
     \frac{1}{2} \, \partial_t^2 ([u]_j')^2 &=
    \left(\lambda_j [u]_j - [Vu]_j
    \right)^2   + [u]_j' \left( \lambda_j \, [u]_j' - [\partial_t (Vu)]_j   \right)
     \, . \label{e:last2}
\end{align}
The lemma follows by combining \eqr{e:last0a} with \eqr{e:last1}
and then
 \eqr{e:last0} with
\eqr{e:last2}.
\end{proof}

The terms on the last  line  of \eqr{e:last4} are the error terms
that vanish when $V=0$.

\subsection{Projecting onto high frequencies}

In contrast to the rotationally symmetric case, the ODE's
 in the previous subsection do  not   imply exponential
 growth or decay of the individual Fourier coefficients.  This is
 because the error terms involve the Fourier coefficients of $Vu$
 and cannot be absorbed.  To get around this, we will instead
 consider simultaneously  all of the Fourier coefficients from some
 point on.
To be precise, we
 fix a large non-negative integer $m$ and let $\HH_{m}(t)$ be the ``high
 frequency''
 part of the norm
 of $u(t,\theta)$   given by
\begin{equation}
    \HH_m (t) =     \sum_{j= m}^{\infty}
        \left( ([u]_j')^2(t) + (1+\lambda_j) \, [u]_j^2(t) \right)  \, .
\end{equation}
Likewise, let $\Low_{m}(t)$ be the left over ``low frequency''
 part
\begin{equation}
    \Low_{m}(t)  =     \sum_{j=0}^{m-1} \left( ([u]_j')^2(t) + (1+\lambda_j) \, [u]_j^2(t) \right)   \, .
\end{equation}
Note that  $\HH_{m}(t)$ is  the contribution on the slice $N\times
\{ t \}$ to the square of the $W^{1,2}(N\times [0,T])$ norm of the
$L^2(N)$-projection of the function $u$ to the eigenspaces from
$m$ to $\infty$. Likewise, $\Low_{m}(t)$ is the $N\times \{ t \}$
part of the square of the $W^{1,2}$ norm of the
$L^2(N)$-projection of the function $u$ to the eigenspaces below
$m$.

The next lemma gives the key differential inequalities for
$\HH_{m}(t)$ and $\Low_{ m}(t)$.

\begin{Lem} \label{l:maindiff}
\begin{align}
    \HH_{m}'' (t) &\geq   (4 \lambda_m - 6) \, \HH_{m}(t) - 3
      \, \int_{N\times \{ t \} } \left[ (Vu)^2 + |\nabla (Vu)|^2 \right] \, d\theta    \,
    , \\
     \Low_{m}''(t)  &\leq (4 \lambda_{m-1} + 6) \, \Low_{m}(t) + 5
      \, \int_{N\times \{ t \} } \left[ (Vu)^2 + |\nabla (Vu)|^2 \right] \, d\theta
      \, .
\end{align}
\end{Lem}

\begin{proof}
We will first prove the bound for  $\HH_{m}''(t)$ and then argue
similarly for $\Low_{m}''(t)$. Applying Lemma \ref{l:diff2} and
then summing over $j$ gives
\begin{align}
      \HH_{ m}''  &=   \sum_{j= m}^{\infty} (4\lambda_j+2) \,
      \left[  ([u]_j')^2 + (1+\lambda_j) \, [u]_j^2 \right]  -
\sum_{j= m}^{\infty}
      \left[    (4\lambda_j+2) \, [u]_j^2 \right]
      \notag \\
    &\quad \quad -   \sum_{j= m}^{\infty}
\left[ (6 \, \lambda_j +2) [u]_j [Vu]_{j} - 2 \, [Vu]_{j}^2 + 2\,
[u]_j' \,
    [Vu]_{j}' \right] \, .
      \label{e:maind1}
\end{align}
The first sum on the first line is at least $(4\lambda_m+2)\,
\HH_m$, while the second is at least $-4 \, \HH_m$. We will now
handle each of the three ``error terms'' in the second line.
  First, the Cauchy-Schwarz inequality  gives
\begin{align}     \label{e:maind2}
      2\, \sum_{j= m}^{\infty}
        (1+\lambda_j)\, \left| [u]_j \, [Vu]_{j} \right|  &\leq
          \sum_{j= m}^{\infty} (1+\lambda_j) \, \left[
          [u]_j^2     +
        [Vu]_{j}^2     \right] \notag \\
        &\leq \HH_{ m} + \int_{N\times \{ t \} } \left[ (Vu)^2 +
|\nabla_{N}(Vu)|^2 \right] \, d\theta
        \, ,
\end{align}
where the second inequality used the standard relation between
  the Fourier coefficients of a function on $N$ and those of its
  derivative.  The second error term is clearly non-negative.
  For the last error term, we  again  use  the Cauchy-Schwarz
inequality   to get
\begin{equation}     \label{e:maind4}
      2 \, \sum_{j= m}^{\infty} \,  \left|
        [u]_j' \,
    [Vu]_{j}' \right| \leq
        \sum_{j= m}^{\infty} \,  \left[
          ([u]_j')^2    +
        ([Vu]_{j}')^2     \right] \leq \HH_{ m} +  \int_{N\times \{ t \} } (\partial_{t}(Vu))^2 \, d\theta
        \, .
\end{equation}
Substituting the bounds \eqr{e:maind2} and \eqr{e:maind4} into
\eqr{e:maind1}  gives
\begin{equation}     \label{e:maind5}
      \HH_{ m}''  \geq (4\lambda_m -6) \, \HH_{  m}  - 3   \int_{N\times \{ t \} } \left[ (Vu)^2 + |\nabla_N(Vu)|^2 \right] \,
      d\theta    -    \int_{N\times \{ t \} } (\partial_{t}(Vu))^2
      d\theta \, ,
\end{equation}
 giving the bound for $\HH_{  m}''$.

The bound for $\Low_{ m}''(t)$ follows similarly, except that the
second term on the first line of \eqr{e:maind1} now has the right
sign and the term corresponding to second error term for $\HH_m$
now has the wrong sign.  We bound this term by
\begin{equation}     \label{e:maind3}
        2 \sum_{j= 0}^{m-1}
         [Vu]_{j}^2 \leq  2 \,  \int_{N\times \{ t \} } (Vu)^2 \, d\theta
        \, .
\end{equation}
\end{proof}

\begin{Cor}     \label{c:corr}
There is a constant $C$   depending only on $||V||_{C^{0,1}}$ (but
{\underline{not}} on $m$) so that
\begin{align}
    \HH_{m}''  &\geq   (4 \lambda_m - C) \, \HH_{m} - C \, \Low_m \, ,
    \\
     \Low_{m}''  &\leq (4 \lambda_{m-1} + C) \, \Low_{ m} + C \, \HH_{m}
      \, .
\end{align}
\end{Cor}

\begin{proof}
Integrating by parts on the closed manifold $N$ and using that
$\nabla = \nabla_N + \partial_t$ gives
\begin{equation}    \label{e:standardf}
 \int_{ N \times \{ t\} } \left[ u^2 + |\nabla u|^2 \right] \, d\theta  =     \sum_{j=0}^{\infty}
\left[ (1+\lambda_j)[u]_j^2   +   ([u]_j')^2 \right] =   \Low_{m}
+ \HH_{m}   \, .
\end{equation}
It is easy to see that there is a constant $c$ depending on
$||V||_{C^{0,1}}$ so that
\begin{equation}    \label{e:bdvuterms}
\int_{ N \times \{ t\} } \left[ (Vu)^2 + |\nabla (Vu)|^2 \right]
\, d\theta  \leq c \, \int_{ N \times \{ t\} } \left[ u^2 +
|\nabla u|^2 \right] \, d\theta =  c \, \left( \Low_{m} + \HH_{m}
\right) \, ,
\end{equation}
where the equality used \eqr{e:standardf}. The corollary follows
from using this bound on the error terms in Lemma
\ref{l:maindiff}.
\end{proof}

\subsection{Taking advantage of gaps in the spectrum}

The next proposition proves a differential inequality for an
integer $m$ where $\lambda_{m} - \lambda_{m-1}$ is large.

\begin{Pro}     \label{l:growth}
There exists a constant $\kappa > 0$ depending on
$||V||_{C^{0,1}}$ so that if   $m$ is an integer with $\lambda_m -
\lambda_{m-1} > \kappa$, then
\begin{equation} \label{e:adv}
    \left( \HH_{m} -  \Low_{m} \right) ''  \geq   (4 \lambda_{m-1} + 2\kappa)
\, \left( \HH_{m} -  \Low_{m} \right)
      \, .
\end{equation}
\end{Pro}

\begin{proof}
To see this, apply Corollary \ref{c:corr}   to get $C$ depending
only on $||V||_{C^{0,1}}$ so that
\begin{align}
    \left( \HH_{m} -  \Low_{m} \right) ''  &\geq
(4 \lambda_m - C) \, \HH_{m} - C \, \Low_m
    - (4 \lambda_{m-1} + C) \, \Low_{ m} - C \, \HH_{m}  \notag \\
    &= (4 \lambda_{m-1} + 2C) \, \left( \HH_{m} -  \Low_{m} \right) +
    4\, (\lambda_m - \lambda_{m-1} - C) \, \HH_{m}
      \, .
\end{align}
\end{proof}

\subsection{The three circles inequality}

We will next use Proposition \ref{l:growth} to prove the three
circles inequality. In fact, we will prove a more general
inequality than the one stated in Theorem \ref{t:ests}.  To state
this, let $N$ be any closed $n$-dimensional Riemannian manifold
satisfying \eqr{e:spgaps} and set
\begin{equation}
I(s)=\int_{N\times \{s\} } \left( u^2 + |\nabla u|^2 \right)
\,d\theta\, .
\end{equation}

\begin{Thm}     \label{t:ests0}
There exists a constant $C > 0$ depending on $||V||_{C^{0,1}}$
 so that if $\alpha$ satisfies
\begin{equation}    \label{e:defalphafirst}
\alpha \geq \frac{1}{T}\left[\log \frac{I(T)}{I(0)}\right]\, ,
\end{equation}
then
\begin{equation}    \label{e:notso}
    \log I(t)\leq   C   +    (c_3 + C + |\alpha|) \, t + \log I(0)
    \, ,
\end{equation}
where the constant $c_3$ is given by{\footnote{The only place
where we use    the spectral gaps given by \eqr{e:spgaps} is to
    get an $m$ satisfying \eqr{e:c3}.}}
\begin{equation}        \label{e:c3}
    c_3 = \min_m \left\{  2\lambda_{m-1}^{1/2} - |\alpha| \, \big| \,
\lambda_m - \lambda_{m-1} > C {\text{ and
    }} 2\lambda_{m-1}^{1/2}   > |\alpha|
    \right\} \, .
\end{equation}
\end{Thm}

Before getting to the proof of Theorem \ref{t:ests0}, we will make
a few remarks.  First, when we have equality in
\eqr{e:defalphafirst}, then Theorem \ref{t:ests0} also applies to
the reflected function $\bar{u}(t) = u(T-t)$ with $-\alpha$ in
place of $\alpha$.  Next, observe that \eqr{e:notso} simplifies
considerably when
\begin{equation}    \label{e:defalphafirst2}
\alpha = \frac{1}{T}\left[\log \frac{I(T)}{I(0)}\right] \geq 0 \,
.
\end{equation}
Namely, when \eqr{e:defalphafirst2} holds, then we get
\begin{equation}    \label{e:alphap}
    \log I(t)\leq   C   +
(c_3 + C) \, t  +\frac{t}{T}\log I(T) +\frac{T-t}{T}\log I(0)\, .
\end{equation}

\begin{proof}
(of Theorem \ref{t:ests0}).
 We will first use  the spectral gap to
bound $\HH_m (t)$ in terms of $\Low_m (t)$ and $\HH_m (0)$ for
some fixed $m$. The key for this is that Proposition
\ref{l:growth} gives a constant $\kappa > 0$ depending on
$||V||_{C^{0,1}}$ so that if $m$ is an integer with
\begin{equation}    \label{e:kappa}
\lambda_m - \lambda_{m-1} \geq \kappa \, ,
\end{equation}
then we have
\begin{equation} \label{e:advya}
    \left( \HH_{m} -  \Low_{m} \right) ''  \geq   (4 \lambda_{m-1} + 2\kappa)
\, \left( \HH_{m} -  \Low_{m} \right)
      \, .
\end{equation}
  Fix some $m$ so that \eqr{e:kappa} holds and
\begin{equation}    \label{e:otherm}
    4 \lambda_{m-1} > \alpha^2 \, .
\end{equation}

On the interval $[0,T]$, we define a function $f$ by
\begin{equation}
f(t)=\e^{-\alpha t} \left( \HH_{m} -  \Low_{m} \right)(t)\, ,
\end{equation}
then
\begin{align}
f(0)&= \left( \HH_{m} -  \Low_{m} \right)(0)\, ,\\
f(T)&\leq   \frac{(\HH_{m} +  \Low_{m})(0)}{ (\HH_{m} +
\Low_{m})(T)} \, \left( \HH_{m} -  \Low_{m} \right)(T) \leq
(\HH_{m} + \Low_{m})(0)
    \, ,\\
f'&=\e^{-\alpha t}\left[\left( \HH_{m} -  \Low_{m} \right)'
-\alpha \left( \HH_{m} -  \Low_{m} \right)\right]\, ,\label{e:f'}\\
f''&=\e^{-\alpha t}\left[\left( \HH_{m} -  \Low_{m} \right)''
-2\alpha \left( \HH_{m} -  \Low_{m} \right)' +\alpha^2 \left(
\HH_{m} -  \Low_{m} \right)\right]\label{e:f''}\, .
\end{align}
By the maximum principle, at an interior maximum $t_0\in (0,T)$
for $f$, $f'(t_0)=0$ and $f''(t_0)\leq 0$.   Hence, by \eqr{e:f'}
and \eqr{e:f''}
\begin{equation}    \label{e:450}
\left( \HH_{m} -  \Low_{m} \right)''(t_0) \leq \alpha^2 \left(
\HH_{m} -  \Low_{m} \right)(t_0)\, .
\end{equation}
However, this contradicts \eqr{e:advya} and \eqr{e:otherm} if
$f(t_0) > 0$, so we conclude that
 $f$ does not have a positive interior maximum.
Therefore, for all $t\in [0,T]$, we have that
\begin{equation}
    f(t)\leq \max \{0 , f(0) , f(T) \} \leq \left( \HH_{m} +  \Low_{m}
    \right)(0) = I(0) \, .
\end{equation}
 This implies that
$\left( \HH_{m} -  \Low_{m} \right)(t) \leq \e^{\alpha t} \, I(0)$
and, hence,
\begin{equation} \label{e:fupp1}
\HH_{m} (t) \leq  \Low_{m} (t) + \e^{\alpha t} \, I (0)\, .
\end{equation}

To complete the proof, we will substitute \eqr{e:fupp1} into a
differential inequality for $\Low_m (t)$ and use this to prove an
exponential upper bound for $\Low_m (t)$. To get the differential
inequality,  recall that \eqr{e:last4a} in Lemma \ref{l:diff2}
gives
\begin{align}
\left| \, \partial_t \left[  ([u]_j')^2 + (1+\lambda_j) \, [u]_j^2
\right] \, \right| &=  \left| \, (4\,
\lambda_j + 2) \, [u]_j \, [u]_j' - 2 [Vu]_{j} \, [u]_j' \right|  \notag \\
& \leq 2 \, (1+\lambda_j)^{1/2} \, \left[ ([u]_j')^2 +
(1+\lambda_j) \, [u]_j^2 \right]  +   [Vu]_{j}^2 +  ([u]_j' )^2 \,
.   \label{e:4p54}
\end{align}
Summing this   up to $(m-1)$ and bounding the $(Vu)$ terms as in
\eqr{e:bdvuterms} gives
\begin{equation}    \label{e:bdlpri}
    \left|  \Low_m ' (t) \right| \leq  \left[ 2 \, \lambda_{m-1}^{1/2} + C \right]
\, \Low_m  (t) + C  \HH_m (t)  \, ,
\end{equation}
where   $C$ depends only on $||V||_{C^{0,1}}$. Using the   bound
\eqr{e:fupp1}, we get
\begin{align}
     \left| \Low_{m}' (t) \right| \leq  c_1 \Low_m  (t) +  C \,
     \,\e^{\alpha t}\,
     I(0) \, ,
\end{align}
where we set
\begin{equation}   \label{e:defc1}
    c_1 = 2 \, \left[  \lambda_{m-1}^{1/2} + C \right]
\end{equation}
 to simplify notation. In particular,   the function
\begin{equation}
    \Low_m (t) \, \e^{-
  c_1 \, t} +  \, \frac{ C}{c_1 -\alpha} \,
     I(0) \, \e^{(\alpha- c_1) \, t}
\end{equation}
 is non-increasing on $[0,T]$; we
     conclude that
\begin{equation}    \label{e:fuppy}
    \Low_{m}(t) \leq \e^{ c_1 \, t} \,
\Low_{m}(0) + \frac{ C}{c_1 -\alpha} \,
     I(0) \,\left(\e^{ c_1 \, t} -
\e^{\alpha t} \right) \leq c_2  \, \e^{ c_1 \, t} \, I(0)
  \,   ,
\end{equation}
where we have set $c_2 = \left( 1 + \frac{ C}{c_1 -|\alpha|}
\right) \geq 1$.   Substituting  \eqr{e:fuppy} into \eqr{e:fupp1}
gives a bound for $I(t) = \HH_m(t) + \Low_m(t)$
\begin{equation}   \label{e:boundll}
    I(t)
    \leq 2\, \Low_{m}(t) + \e^{\alpha t} \, I(0)
    \leq 2 \, c_2 \,
   \e^{ c_1 \, t} \, I(0)   +
   \e^{\alpha t} \, I(0) \leq (2 \, c_2 + 1) \,  \e^{ c_1 \, t} \, I(0)
     \, ,
\end{equation}
where the last inequality used  that $c_1 > |\alpha|$.  The
theorem follows from \eqr{e:boundll}.
\end{proof}

We will next apply the three circles inequality of Theorem
\ref{t:ests0} to prove uniform estimates for the $W^{1,2}$ norm of
an at most exponentially growing solution $u$ on the half-cylinder
$N \times [0,\infty)$, i.e., to prove Corollary
\ref{c:ests}.{\footnote{A similar argument, with Theorem
\ref{t:estsc0} in place of Theorem \ref{t:ests0}, gives  a
corresponding result on spheres and Zoll surfaces even when $V$ is
just bounded.}}   As in the statement of Theorem \ref{t:ests0}, we
will let $I(s)$ denote the $W^{1,2}$ of $u$ on $N\times \{ s \}$.

\begin{proof}
(of Corollary \ref{c:ests}.)    We will  assume   that $\alpha
> 0$ (we can do this since $H_{\alpha} \subset H_{\bar{\alpha}}$
whenever $\alpha \leq \bar{\alpha}$). We will first use the
definition of $H_{\alpha}$ to bound $I(T)$ for large values of $T$
and then use the three circles inequality to bound $I(t)$ in terms
of $I(0)$ and $I(T)$.

 The interior Schauder estimates (theorem
$6.2$ in \cite{GT}) give a constant $C$ depending only on  the
$C^{\beta}$ norm of $V$, where
 $\beta \in (0,1)$ is fixed, so that for all $t \geq 1$
\begin{equation}
    I(t) = \int_{ N \times \{ t\} } \left( u^2 + |\nabla u|^2 \right)  \,d\theta \leq
    C \, \sup_{N\times [t-1,t+1]}
\, |u|^2  \, .
\end{equation}
If we also bring in the definition of $H_{\alpha}$, i.e.,
\eqr{e:defhal}, then we get that
\begin{equation}    \label{e:upperg}
    \limsup_{T\to \infty} \, \left( \e^{-2 \,\alpha \, (T+1)} \,
    I(T) \right) = 0 \,
    .
\end{equation}
Note that the bound \eqr{e:upperg} applies only in the limit as
$T$ goes to $\infty$ and, hence, does not give the Corollary.
However, it does give a sequence $T_j \to \infty$ with
\begin{equation}
    \frac{\log I(T_j) - \log I(0)}{T_j} \leq 2\alpha \, .
\end{equation}
 Applying the three circles inequality of Theorem
\ref{t:ests0} on $[0,T_j]$
  gives
\begin{equation}
    \log I(t) \leq   C\,(1+t) +  2\alpha \, t +  \log I(0) \, ,
\end{equation}
and exponentiating this gives the corollary.
\end{proof}

Note that $\nu$ in Corollary \ref{c:ests} has to also depend on
the  norm of  $V$ and not just on $\alpha$.  In particular, $\nu$
may have to be chosen positive even when $\alpha$ is zero. This
can easily be seen by the following example for the
one-dimensional Schr\"odinger equation. Suppose that  $\Psi: \RR
\to \RR$ is a smooth monotone non-decreasing function with
\begin{equation}
    \Psi (x) = -1 {\text{ for }} x < -1, \,  \Psi (x) = x {\text{ on }}
    [0,\ell], \, \Psi (x) = \ell +1 {\text{ for }} x > \ell \, .
\end{equation}
  Then $u(x) = \e^{\Psi (x)}$ satisfies the Schr\"odinger
equation $u'' = \left( (\Psi')^2 + \Psi'' \right) \, u = V \, u$
for a bounded potential $V$ with compact support. However, $u$ is
constant on each end, but grows exponentially on $[0,\ell]$.
Similarly, one can easily construct a bounded (but no longer with
compact support) potential so that the corresponding Schr\"odinger
equation has a solution that grows exponentially on $[0,\ell]$,
yet at infinity the solution vanishes.

\subsection{Unique continuation}

Rather than  stating the most general   three circles inequality
possible, we have tailored the statement of Theorem \ref{t:ests0}
to fit our geometric applications.  We will show here how to
modify the proof to get strong unique continuation for the
operator $L$ on $N \times [0, \infty)$ since this is of
independent interest:

\begin{Pro}     \label{p:uc}
If $u$ is a solution on $N \times [0,\infty)$, where $N$ satisfies
\eqr{e:spgaps}, and
\begin{equation}    \label{e:inord}
    \liminf_{T\to \infty} \, \frac{\log I(T)}{T} = - \infty \, ,
\end{equation}
then $u$ is the constant solution $u=0$.
\end{Pro}

\begin{proof}
Observe that Theorem \ref{t:ests0} implies that if $I(0) = 0$,
then $u$ is the constant solution $u=0$ (this also follows from
\cite{A}). Therefore,
  it suffices to show that $I(0) = 0$.  We will argue by
  contradiction, so suppose that $I(0) > 0$.  After replacing $u$
  by $I^{-1/2}(0) \, u$, we can assume that $I(0)=1$.

 Choose some $m_0$ so that $\HH_m (0) < \Low_m (0)$ for all $m
 \geq m_0$.
 Using the spectral gaps of $N$ and Proposition
\ref{l:growth}, we can choose an
 arbitrarily large integer $m > m_0$ with
 \begin{equation} \label{e:advyauc}
    \left( \HH_{m} -  \Low_{m} \right) ''  \geq   (4 \lambda_{m-1} + 1)
\, \left( \HH_{m} -  \Low_{m} \right)
      \, .
\end{equation}
The rapid decay given by \eqr{e:inord} guarantees that we can find
$T> 0$ with
  \begin{equation}  \label{e:vanish2N}
    \frac{\log I(T)}{T} < -4 \, \lambda_{m-1}^{1/2} \, .
  \end{equation}
On the interval $[0,T]$, we define a function $f$ by
\begin{equation}
f(t)=\e^{2 \lambda_{m-1}^{1/2} t} \left( \HH_{m} -  \Low_{m}
\right)(t)\, .
\end{equation}
Using first that $m\geq m_0$ and then using \eqr{e:vanish2N}, we
get that
\begin{equation}    \label{e:endpoints}
    f(0) < 0 {\text{ and }} f(T) \leq \e^{2 \lambda_{m-1}^{1/2} t} \,
    I(T) < \e^{- 2 \lambda_{m-1}^{1/2} T}  \, .
\end{equation}
 Using the maximum principle  as in \eqr{e:advya}--\eqr{e:450},
we see that $f$ cannot have a positive interior maximum and,
hence, that
\begin{equation}    \label{e:bdhml}
 \left( \HH_{m} -  \Low_{m}
\right)(t)  \leq \e^{-2 \lambda_{m-1}^{1/2} (t+T)}   \, .
\end{equation}
Combining this with the bound \eqr{e:bdlpri} for $\Low_m '$ gives
\begin{equation}    \label{e:bdlpriuc}
    \left|  \Low_m ' (t) \right| \leq  \left[ 2 \, \lambda_{m-1}^{1/2} + C \right]
\, \Low_m  (t) + C  \, \e^{-2 \lambda_{m-1}^{1/2} (t+T)}  \, ,
\end{equation}
where   $C$ depends only on $||V||_{C^{0,1}}$.  It follows that
the function
\begin{equation}    \label{e:nondec1}
     \e^{\left[ 2 \, \lambda_{m-1}^{1/2} + C \right] \, t} \,  \Low_m (t)  +   \e^{Ct -2 \lambda_{m-1}^{1/2} T}
\end{equation}
is non-decreasing on $[0,T]$.  Evaluating this function at $0$ and
$T$ gives
\begin{equation}    \label{e:nondec2}
     \Low_m (0) \leq \e^{\left[ 2 \, \lambda_{m-1}^{1/2} + C \right] \, T} \,  \Low_m (T)  +   \e^{CT -2 \lambda_{m-1}^{1/2} T}
        \leq   2 \,  \e^{CT -2 \lambda_{m-1}^{1/2} T} \, .
\end{equation}
However, since $\lambda_{m-1}$ can be arbitrarily large and $C$ is
fixed (i.e., does not depend on $m$), we conclude that $\Low_m
(0)=0$. Finally, this gives the desired contradiction since $I(0)=
\Low_m (0) + \HH_m (0) =1$ and $\Low_m (0)
> \HH_m (0)$.
\end{proof}

\section{A three circles theorem for  bounded potentials}
\label{s:ne}

We will show in this section that Theorem \ref{t:ests} holds even
for potentials that are just bounded, i.e.,   the potential $V$
does not have to be Lipschitz; this is Theorem \ref{t:estsc0}.
This result will require larger spectral gaps than were needed for
the arguments in the Lipschitz case. Throughout this section, $u$
will be a solution of
\begin{equation}    \label{e:equn2}
    \Delta u = - Vu
\end{equation}
 on a product $N \times [0,t]$, where the potential $V$ is bounded,
 but is not assumed to be Lipschitz.

The main place where the Lipschitz bound entered previously was
when we took second derivatives of $|\nabla u|^2$.  To avoid doing
this, we will work with the $L^2$ norm of the spectral projections
of a solution $u$. Namely,  we
 fix a large non-negative integer $m$ and let $\HHH_{m}(t)$ and $\Loww_m (t)$ be the ``high
 frequency'' and ``low frequency'' parts, respectively,
 of $u(t,\theta)$   given by
\begin{equation}
    \HHH_m (t) =     \sum_{j= m}^{\infty}
        [u]_j^2(t)   {\text{ and }}
\Loww_m (t) =     \sum_{j= 0}^{m-1}
        [u]_j^2(t)   \, .
\end{equation}
Note that  $\HHH_{m}^{1/2}$ is the  $L^{2}$ norm of the projection
of   $u$ to the eigenspaces from $m$ to $\infty$.

As in   Section \ref{s:sthree}, we will derive a second order ODE
for $\HHH_{m}(t)$ and use this to control its growth.
Unfortunately,   the   ODE \eqr{e:last0} for the quantity
$[u]_j^2$, and thus also for $\HHH_{m}(t)$, is not as nice as for
$[u_t]_j^2 + \lambda_j \, [u]_j^2$ because of the $[u_t]_j^2$ term
on the right hand side. We will use the next lemma to get around
this.

\begin{Cor}     \label{c:fvarn}
Given $t_1 < t_2$, we get that
\begin{equation}
\left([u_t]_j^2 -  \lambda_j \, [u]_j^2  \right)(t_2) - \left(
[u_t]_j^2 -  \lambda_j \, [u]_j^2  \right)(t_1) = - 2 \,
\int_{t_1}^{t_2} \left( [Vu]_j \, [u_t]_j \right) \, dt \, .
\end{equation}
\end{Cor}

\begin{proof}
 Differentiating $\left(
[u_t]_j^2 -  \lambda_j \, [u]_j^2  \right)$ and then using Lemma
\ref{l:diff} gives
\begin{equation}
    \partial_t \left([u_t]_j^2 -  \lambda_j \, [u]_j^2  \right) =
    2 (\lambda_j \, [u]_j   - [Vu]_j )\,  [u_t]_j
    - 2 \, \lambda_j \, [u]_j \, [u_t]_j
    = -
    2 [Vu]_j \,  [u_t]_j \, .
\end{equation}
The corollary now follows
 from the fundamental
theorem of calculus.
\end{proof}

The next lemma will give the key differential inequality for
$\HHH_{m}(t)$.  To state this, it will be useful to define
$\elltwo (t)$ to be the square of the $L^2$ norm of $u$ on
$N\times \{ t \}$, i.e.,
\begin{equation}    \label{e:ell2}
    \elltwo (t) = \int_{N \times \{ t \} } u^2 \, d\theta \, .
\end{equation}

\begin{Lem} \label{l:maindiffa}
\begin{align}     \label{e:maind5a}
      \HHH_{ m}''(t)  &\geq (4\lambda_m -1) \, \HHH_{  m}(t) -      \int_{N\times \{ t \} }   (Vu)^2   \,
      d\theta - 2 \, \int_{N\times \{ t_0 \} } |\nabla u|^2 \, d\theta -   \elltwo'(t) +   \elltwo'(t_0)
      \notag \\
      &\quad \quad \quad
      -  2   \int_{N \times (t_0,t) } \left[ (V^2 + |V|)\, u^2 \right]   \, .
\end{align}
\end{Lem}

\begin{proof}
  Applying Lemma \ref{l:diff}
and then summing over $j$ gives
\begin{align}   \label{e:maind1a}
        \HHH_{ m}''  &=  2\, \sum_{j= m}^{\infty} \left[
        [u_t]_j^2 +
      \lambda_j \, [u]_j^2  - [u]_j \, [Vu]_j \right]
      \\
    &=  2\, \sum_{j= m}^{\infty} \left[
      2 \,\lambda_j \, [u]_j^2 - [u]_j \, [Vu]_j +
\left(  [u_t]_j^2 - \lambda_j \, [u]_j^2    \right)(t_0) - 2 \,
\int_{t_0}^{t} \left( [Vu]_j \, [u_t]_j \right) \, ds
        \right] \, , \notag
\end{align}
where the second equality used Corollary \ref{c:fvarn}.    We will
now handle each of the three ``error terms'' in the second line.
First, the Cauchy-Schwarz inequality   gives
\begin{equation}     \label{e:maind2a}
      2 \, \left| \sum_{j= m}^{\infty}
          [u]_j \, [Vu]_{j} \right| \leq
          \sum_{j= m}^{\infty} \left[ [u]_j^2 + [Vu]_{j}^2 \right]
          \leq  \HHH_{ m} + \int_{N\times \{ t \} }   (Vu)^2   \, d\theta
        \, ,
\end{equation}
where the second inequality used the standard relation between
  the Fourier coefficients of a function on $N$ and  its
  $L^2$ norm.  The second error term
   is bounded by
\begin{equation}     \label{e:maind3a}
        2\, \left|  \sum_{j= m}^{\infty} \left( [u_t]_j^2 - \lambda_j \, [u]_j^2    \right) \right|(t_0) \leq
            2\, \sum_{j= 0}^{\infty} \left( \lambda_j \, [u]_j^2 + [u_t]_j^2  \right)(t_0)
            = 2\, \int_{N\times \{ t_0 \} } |\nabla u|^2 \, d\theta
        \, ,
\end{equation}
where the  equality used the standard relation between
  the Fourier coefficients of a function on $N$ and those of its
  derivative.
Similarly, for the last error term,    the Cauchy-Schwarz
inequality gives
\begin{equation}     \label{e:maind4a}
      4 \, \left| \sum_{j= m}^{\infty}
        \int_{t_0}^{t} \left( [Vu]_j \, [u_t]_j \right) \, ds \right| \leq
     2\, \int_{t_0}^{t}
        \left( \int_{N \times \{s\} } \left[ (Vu)^2  +  (u_t)^2 \right] \, d\theta \right) \, ds
    \, .
\end{equation}
The first term in \eqr{e:maind4a} is of the right form, but it
will be convenient to get a lower order bound for the $(u_t)^2$
term. To do this,
 we use Stokes' theorem to get
\begin{equation}     \label{e:mainet1b}
          2 \, \int_{N \times (t_0,t) } \left[ |\nabla
         u|^2 - V \, u^2  \right] = \elltwo'(t) - \elltwo'(t_0)
    \, ,
\end{equation}
so we get that
\begin{equation}     \label{e:mainet2}
          2 \, \int_{N \times (t_0,t) }   (u_t)^2 \leq  2 \, \int_{N \times (t_0,t) }   |\nabla
         u|^2 = \elltwo'(t) - \elltwo'(t_0) + 2 \,
     \int_{N \times (t_0,t) }   V \, u^2
    \, .
\end{equation}
Finally, substituting the bounds \eqr{e:maind2a}--\eqr{e:maind4a}
and \eqr{e:mainet2} into \eqr{e:maind1a}  gives the lemma.
\end{proof}

We get the following immediate corollary of Lemma
\ref{l:maindiffa}; note that the square of the $W^{1,2}$ norm of
the projection to the low frequencies, i.e., $\Low_m$, appears in
the bound.

\begin{Cor} \label{c:maindiffa}
There exists a constant $C$ depending only on $\sup |V|$ so that
\begin{align}     \label{e:maindifff}
      \HHH_{ m}''  &\geq (4\lambda_m -C) \, \HHH_{  m} -  \HHH_{  m}' - 3 \, I(t_0) - C \, \Loww_m
    - 2 \, (\Loww_m \, \Low_m )^{1/2}   \notag \\
    &\quad \quad \quad
      - C \,
       \int_{t_0}^t \left[ \HHH_m (s) + \Loww_m (s) \right] \, ds \,
       .
\end{align}
\end{Cor}

\begin{proof}
To bound the first ``error term'' from Lemma \ref{l:maindiffa},
  bound  $V$ by $\sup |V|$ to get
\begin{equation}
     \int_{N\times \{ t \} }   (Vu)^2 \, d\theta \leq \sup |V|^2 \, \int_{N\times \{ t \} }   u^2 \, d\theta
        = \sup |V|^2 \, \left[ \HHH_m (t) +  \Loww_m (t) \right] \,
        .
\end{equation}
 Similarly, the last error term is bounded by
\begin{equation}
    2 \, \int_{N \times (t_0,t) } \left[ (V^2 + |V|)\, u^2 \right] \leq  2\,  (\sup |V| + \sup |V|^2) \,
       \int_{t_0}^t \left[ \HHH_m (s) + \Loww_m (s) \right] \, ds \, .
\end{equation}
The second error term $2 \, \int_{N\times \{ t_0 \} } |\nabla u|^2
\, d\theta$ is trivially bounded by $2\,I(t_0)$.  This leaves only
the two $\elltwo'$ terms.  Use Cauchy-Schwarz to bound the second
of these by
\begin{equation}
    |\elltwo'(t_0)| = 2 \, \left| \int_{N \times \{ t_0 \} } u \, u_t \,
    d\theta \right| \leq \int_{N \times \{ t_0 \} } (u^2 + u_t^2) \,
         d\theta \leq I(t_0)  \, .
\end{equation}
To bound $\elltwo'(t)$, observe first that
\begin{equation}        \label{e:bdlprimea}
    \left| \Loww_m'   \right| = 2 \left| \sum_{j=0}^{m-1} \, [u]_j \, [u_t]_j \right| \leq
    2 \, \left[ \sum_{j=0}^{m-1} \, [u]_j^2  \right]^{1/2} \,
    \left[ \sum_{j=0}^{m-1} \, [u_t]_j^2  \right]^{1/2} \leq
     2 \, (\Loww_m \, \Low_m )^{1/2}    \, ,
\end{equation}
so we get
\begin{equation}    \label{e:useqm}
         \elltwo '   =    \HHH_m'   +    \Loww_m'   \leq
         \HHH_m'   + 2 \, (\Loww_m \, \Low_m )^{1/2}  \,
       .
\end{equation}
The corollary now follows from
 Lemma \ref{l:maindiffa}.
\end{proof}

 The next lemma gives the key differential inequality for
$\Low_m$ that will be used later to get an upper bound for $\Low_m
(t)$ (a similar, but slightly less sharp, bound was given in
\eqr{e:bdlpri}).

\begin{Lem}     \label{l:Qbd}
There exists a constant $C > 0$ depending only on  $\sup |V|$ so
that
\begin{equation} \label{e:Qbd}
    \left| \Low_m '  \right|  \leq   2 (\lambda_{m-1} + 1)^{1/2}
    \Low_m   + C \, \sqrt{\Loww_m   + \HHH_m   } \, \sqrt{\Low_m }
      \, .
\end{equation}
\end{Lem}

\begin{proof}
To get the differential inequality,  recall from \eqr{e:4p54} that
Lemma \ref{l:diff2} gives
\begin{align}
\partial_t \left[  ([u]_j')^2 +
(1+\lambda_j) \, [u]_j^2 \right] &=  (4\, \lambda_j + 2) \, [u]_j
\, [u]_j' - 2   [u]_j'\, [Vu]_j   \notag \\
& \leq 2 \,  (\lambda_j + 1)^{1/2} \, \left[ [u_t]_j^2 +
(\lambda_j + 1) \, [u]_j^2 \right]  + 2 \, \left| [Vu]_{j}   \,
[u_t]_j \right| \, .
\end{align}
Summing this   up to $(m-1)$ and then using the Cauchy-Schwarz
inequality for series gives
\begin{equation}
    \left|  \Low_m '   \right| \leq  2 \,   (\lambda_{m-1} + 1)^{1/2}
\, \Low_m    + C \, \sqrt{\elltwo  } \, \sqrt{\Low_m } \, ,
\end{equation}
where the constant $C$ depends only on  $\sup |V|$.
\end{proof}

\subsection{Exponentially weighted sup bounds for $\HHH_m$, $\Loww_m$ and
$\Low_m$}

We will next record an immediate consequence of Corollary
\ref{c:maindiffa} where the last three terms in \eqr{e:maindifff}
are bounded in terms of the sup norms of $\HHH_m$, $\Loww_m$ and
$\Low_m$ against an exponential weight. To make this precise, for
each constant $\alpha
> 0$, we define the exponentially weighted sup norm bounds
$\bar{h}_{\alpha,m}$, $\bar{\ell}_{\alpha , m}$ and
$\ell_{\alpha,m}$ by
\begin{align}
    \bar{h}_{\alpha,m} &= \max_{ [0,T] } \, \left[ \HHH_m (t) \, \e^{-
    \alpha \, t} \right] \, , \\
    \bar{\ell}_{\alpha,m} &= \max_{ [0,T] } \, \left[ \Loww_m (t) \, \e^{-
    \alpha \, t} \right] \, , \\
     \ell_{\alpha,m} &= \max_{ [0,T] } \, \left[ \Low_m (t) \, \e^{-
    \alpha \, t} \right] \, .
\end{align}
Clearly, by definition, we have that
\begin{equation} \label{e:bdhm}
    \HHH_m (t) \leq \bar{h}_{\alpha,m}   \, \e^{
    \alpha \, t}  \, , \,
     \Loww_m (t) \leq \bar{\ell}_{\alpha,m}   \, \e^{
    \alpha \, t}  \, , \,  
     \Low_m (t) \leq \ell_{\alpha,m}   \, \e^{
    \alpha \, t} \, . 
\end{equation}

Substituting these bounds into the differential inequality for
$\HHH_m$ gives:

\begin{Cor}     \label{c:diweight}
There exists a constant $C$ depending only on $\sup |V|$ so that
for $\alpha \geq 1$
\begin{equation}     \label{e:diw}
      \HHH_{ m}''  \geq (4\lambda_m -C) \, \HHH_{  m} -    \HHH_{  m}' -  3 I(t_0) -
      \left[ C \, \left( \bar{\ell}_{\alpha, m} +  \bar{h}_{\alpha,m} \right) + 2 \,
      \left( \bar{\ell}_{\alpha, m} \, \ell_{\alpha, m} \right)^{1/2} \right] \, \e^{\alpha \,
    t} \, .
\end{equation}
\end{Cor}

\begin{proof}
The corollary will follow directly from Corollary
\ref{c:maindiffa} by using \eqr{e:bdhm}  to bound the last three
terms in \eqr{e:maindifff}.
 The bounds on $\Loww_m$ and $(\Loww_m \, \Low_m )^{1/2}$ follow immediately from
\eqr{e:bdhm}.  Finally, to bound the last term in
\eqr{e:maindifff}, note that
\begin{equation}
       \int_{t_0}^t \left[ \HHH_m (s) + \Loww_m (s) \right] \, ds \leq
        (  \bar{h}_{\alpha,m} + \bar{\ell}_{\alpha,m}) \, \int_{t_0}^t \e^{\alpha
        \, s} \, ds \leq \frac{ \bar{h}_{\alpha,m} + \bar{\ell}_{\alpha,m} \,}
        {\alpha} \, \e^{\alpha \, t}
       \, .
\end{equation}
The corollary now follows from substituting these bounds into
\eqr{e:maindifff}.
\end{proof}

\subsection{Taking advantage of gaps in the spectrum}

 Fix a
constant $\kappa \geq 1$   to be chosen (depending only on $\sup
|V|$) and then choose a constant $\bar{\alpha} \geq \kappa$ with
\begin{equation}    \label{e:defalpha}
    \alpha \equiv \frac{1}{T}\left[\log \frac{ I(T)}
  {I(0)}\right]  \leq \bar{\alpha}  \,
    ,
  \end{equation}
  and so that there exists $m$ with
  \begin{align} \label{e:chom1}
    2 \, (\lambda_{m-1} + 1)^{1/2}  + 1 &\leq \bar{\alpha}  \, , \\
    \bar{\alpha}^2 +   \bar{\alpha}    &\leq 4 \, \lambda_m - \kappa   \, ,
    \label{e:chom2}
\end{align}
We will use the spectral gaps to show that such an $\bar{\alpha}$
always exists when $N$ is a round sphere of a Zoll surface.

\begin{Pro}     \label{p:technicalbound}
If $\bar{\alpha} \geq 1$ satisfies  \eqr{e:defalpha},
\eqr{e:chom1}, and \eqr{e:chom2} for some constant $\kappa \geq 1$
depending only on
  $\sup |V|$, then for all $t \in [0,T]$
\begin{equation}
      \int_{N \times \{
    t \} } u^2   \, d\theta \leq C \,I(0) \, \e^{\bar{\alpha} \, t} \, ,
\end{equation}
where $C$ depends only on $\sup |V|$.
\end{Pro}

The proof of Proposition \ref{p:technicalbound} will be divided
into four steps. First, we bound $\bar{\ell}_{\bar{\alpha} , m}$
in terms of $\ell_{\bar{\alpha} , m}$ and $ I(0)$. Second,  we use
\eqr{e:chom1} to bound $\ell_{\bar{\alpha} , m}$  in terms of
$\bar{\ell}_{\bar{\alpha} , m}$, $\bar{h}_{\bar{\alpha} , m}$ and
$I(0)$. Third, we combine these to bound both
$\bar{\ell}_{\bar{\alpha} , m}$ and  $\ell_{\bar{\alpha} , m}$ in
terms of $\bar{h}_{\bar{\alpha} , m}$ and $I(0)$.  Finally, we
substitute these bounds into the differential inequality for
$\HHH_m''$ to bound $\bar{h}_{\bar{\alpha} , m}$ in terms of
$I(0)$.  In this last step, $\bar{h}_{\bar{\alpha} , m}$ will show
up on both sides of the inequality, but \eqr{e:chom2} will allow
us to absorb the terms on the right hand side.

\begin{proof}
(of Proposition \ref{p:technicalbound}.)

{\underline{Bounding $\bar{\ell}_{\bar{\alpha} , m}$}}.  To bound
$\bar{\ell}_{\bar{\alpha} , m}$, use \eqr{e:bdlprimea} to get
\begin{equation}        \label{e:bdlprime2}
    \left| \Loww_m'   \right|   \leq
     2 \, (\Loww_m \, \Low_m )^{1/2} \leq  2 \, \Loww_m^{1/2} \, \ell_{\bar{\alpha} , m} \, \e^{\bar{\alpha} \, t/2}  \,
     .
\end{equation}
On the interval $[0,T]$, we define a function $
f_1(t)=\e^{-\bar{\alpha} t}   \Loww_{m}(t)$, so that
\begin{align}
f_1(0) &\leq I(0)  {\text{ and }} f_1(T) \leq I(0) \, ,\\
f_1'&=\e^{-\bar{\alpha} t}\left[ \Loww_{m}' -\bar{\alpha} \,
\Loww_{m} \right]\, .
\end{align}
Observe that the maximum of $f_1$ on $[0,T]$ is precisely
$\bar{\ell}_{\bar{\alpha} , m}$.  Hence, if the maximum of $f_1
(t)$ occurs at a point $s$ in the interior $(0,T)$, then we get
\begin{equation}
     \bar{\alpha} \, \bar{\ell}_{\bar{\alpha} , m} \, \e^{\bar{\alpha} s} = \bar{\alpha} \, \Loww_m(s)  = \Loww_m'(s) \leq
        2 \, \bar{\ell}_{\bar{\alpha} , m}^{1/2} \,
    \ell_{\bar{\alpha} , m}^{1/2}   \, \e^{\bar{\alpha} \,s}
        \, .
\end{equation}
Combining this with the fact that $f_1 \leq I(0)$ at both
endpoints gives that
\begin{equation}    \label{e:gotellalp}
       \bar{\ell}_{\bar{\alpha} , m}  \leq
        \frac{4}{\bar{\alpha}^2 } \,
    \ell_{\bar{\alpha} , m}    + I(0)
        \, .
\end{equation}

{\underline{Bounding $\ell_{\bar{\alpha} , m}$}}. Substituting the
bound \eqr{e:bdhm} for $\HHH_m$ into Lemma \ref{l:Qbd} gives
\begin{equation}
    \left| \Low_m'(t) \right|  \leq    2 (\lambda_{m-1} + 1)^{1/2}
    \Low_m (t) + C \, \sqrt{\bar{h}_{\bar{\alpha} , m} + \bar{\ell}_{\bar{\alpha} , m} }\, \e^{\bar{\alpha} \,t/2}
         \, \sqrt{\Low_m (t)}  \, .
\end{equation}
Consequently, if the maximum of $f_2 (t) = \e^{-\bar{\alpha} \, t}
\, \Low_m (t)$ occurs at a point $s$ in the interior $(0,T)$, then
we get
\begin{equation}
     \bar{\alpha} \, \ell_{\bar{\alpha} , m} \, \e^{\bar{\alpha} s} = \bar{\alpha} \, \Low_m(s)  = \Low_m'(s) \leq
    2 (\lambda_{m-1} + 1)^{1/2} \, \ell_{\bar{\alpha} , m} \, \e^{\bar{\alpha}
    s}+ C (\bar{h}_{\bar{\alpha} , m} + \bar{\ell}_{\bar{\alpha} , m} )^{1/2} \,  \ell_{\bar{\alpha} , m}^{1/2} \, \e^{\bar{\alpha} \,s}
        \, ,
\end{equation}
so we would get that
\begin{equation}
     \left( \bar{\alpha} - 2 (\lambda_{m-1} +1)^{1/2} \right) \, \ell_{\bar{\alpha} , m}    \leq
     C (\bar{h}_{\bar{\alpha} , m} + \bar{\ell}_{\bar{\alpha} , m} )^{1/2} \,  \ell_{\bar{\alpha} ,
     m}^{1/2}
        \, .
\end{equation}
 Using \eqr{e:chom1}, we would
 then get that
\begin{equation}    \label{e:qfromh}
    \ell_{\bar{\alpha} , m} \leq C \, \left( \bar{h}_{\bar{\alpha} , m} + \bar{\ell}_{\bar{\alpha} , m} \right) \,
    ,
\end{equation}
where the constant $C$ depends only on $\sup |V|$. Combining this
with the fact that $\e^{-\bar{\alpha} \, t} \, \Low_m (t)  \leq
I(0)$ at the endpoints, we get that
\begin{equation}    \label{e:qfromh2}
    \ell_{\bar{\alpha} , m} \leq  C \, \left( \bar{h}_{\bar{\alpha} , m} + \bar{\ell}_{\bar{\alpha} , m} \right) + I(0) \,
    .
\end{equation}

{\underline{Bounding both $\ell_{\bar{\alpha} , m}$ and
$\bar{\ell}_{\bar{\alpha} , m}$ in terms of $\bar{h}_{\bar{\alpha}
, m}$ and $I(0)$}}.   If we substitute the bound \eqr{e:gotellalp}
into \eqr{e:qfromh2}, then we get
\begin{equation}    \label{e:qfromh3}
    \ell_{\bar{\alpha} , m} \leq C \, I(0) + C \, \left( \bar{h}_{\bar{\alpha} , m} + \frac{4}{\bar{\alpha}^2 } \,
    \ell_{\bar{\alpha} , m}       \right)   \, .
\end{equation}
As long as $\bar{\alpha}^2 \geq 8\, C$, then we can absorb the
$\ell_{\bar{\alpha} , m}$ term on the right to get
\begin{equation}    \label{e:qfromh4}
    \ell_{\bar{\alpha} , m} \leq 2 \, C \, I(0) + 2 \, C \,   \bar{h}_{\bar{\alpha} , m}     \,
    .
\end{equation}
Finally, substituting this back into \eqr{e:gotellalp}  gives
\begin{equation}    \label{e:gotellalp2}
       \bar{\ell}_{\bar{\alpha} , m}  \leq  2 \, I(0) +
    \bar{h}_{\bar{\alpha} , m}
        \, .
\end{equation}

{\underline{Bounding $\bar{h}_{\bar{\alpha} , m}$ in terms of
$I(0)$}}. The starting point is to substitute the bounds
\eqr{e:qfromh4} and \eqr{e:gotellalp2} into Corollary
\ref{c:diweight}, to get
\begin{align}     \label{e:3ci1}
      \HHH_{ m}''  &\geq (4\lambda_m -C) \, \HHH_{  m} -   \HHH_{  m}'
      - 3 I(0) -
     \left[ C \, \left( \bar{\ell}_{\alpha, m} +  \bar{h}_{\alpha,m} \right) + 2 \,
      \left( \bar{\ell}_{\alpha, m} \, \ell_{\alpha, m} \right)^{1/2} \right] \,  \e^{\bar{\alpha} \,
    t}      \\
    &\geq      (4\lambda_m -C) \, \HHH_{  m} -   \HHH_{  m}'
      - C  \,
     \left[   \bar{h}_{\alpha,m} + I(0) \right] \,  \e^{\bar{\alpha} \,
    t}
    \, ,  \notag
\end{align}
where   $C$ depends only on $\sup |V|$ and we absorbed the $3 \,
I(0)$ term into the last term.  Define a function
$f_3(t)=\e^{-\bar{\alpha} t} \HHH_{m}(t)$ on $[0,T]$, so that
$f_3$ is bounded by $I(0)$ at $0$ and $T$ and
\begin{align}
f_3'&=\e^{-\bar{\alpha} t}\left[ \HHH_{m}'
-\bar{\alpha} \, \HHH_{m}  \right]\, ,\label{e:f3'}\\
f_3''&=\e^{-\bar{\alpha} t}\left[  \HHH_{m} '' -2\bar{\alpha}
\HHH_{m}' +\bar{\alpha}^2 \, \HHH_{m} \right]\label{e:f3''}\, .
\end{align}
At an interior maximum $s \in (0,T)$ for $f_3$, we have
$f_3'(s)=0$ and $f_3''(s)\leq 0$.   Hence, by \eqr{e:f3'} and
\eqr{e:f3''}
\begin{align}
    \HHH_m' (s) &= \bar{\alpha} \, \HHH_m (s) =  \bar{\alpha} \, \bar{h}_{\bar{\alpha} , m} \, \e^{\bar{\alpha} \,
    s} \, , \\
 \HHH_{m}''(s) &\leq \bar{\alpha}^2 \,
\HHH_{m}  (s) = \bar{\alpha}^2 \, \bar{h}_{\bar{\alpha} , m} \,
\e^{\bar{\alpha} \, s} \, .
\end{align}
Combining these with \eqr{e:3ci1} and multiplying through by
$\e^{-\bar{\alpha} \, s}$ would give
\begin{equation}
 (4\lambda_m -C) \, \bar{h}_{\bar{\alpha},m}  -   \bar{\alpha} \, \bar{h}_{\bar{\alpha} , m} -
      C \,   ( \bar{h}_{\bar{\alpha},m} + I(0))   \leq \bar{\alpha}^2 \, \bar{h}_{\bar{\alpha} ,
      m} \, .
\end{equation}
If we now substitute \eqr{e:chom2} into this, then we would get
that
\begin{equation}    \label{e:wouldget}
  \bar{h}_{\bar{\alpha},m} \leq \left( 4\lambda_m - 2 C - \bar{\alpha}^2
  -   \bar{\alpha}
  \right) \bar{h}_{\bar{\alpha},m}  \leq
      C \,     I(0)  \, .
\end{equation}
On the other hand, if the maximum of $f_3$ occurs at $0$ or $T$,
  then we would get $\bar{h}_{\bar{\alpha},m} \leq   I(0)$
so we conclude that \eqr{e:wouldget} holds in either case.
Combining all of this gives that
\begin{equation}
    \max_{[0,T]} \, \left( \e^{-\bar{\alpha} \, t} \, \int_{N \times \{
    t \} } u^2  \, d\theta \right) \leq \bar{h}_{\bar{\alpha} , m} + \bar{\ell}_{\bar{\alpha} , m} \leq
    C \, I(0) \, .
\end{equation}
\end{proof}

\subsection{Choosing $\bar{\alpha}$}

We will now show that Proposition \ref{p:technicalbound} implies
Theorem \ref{t:estsc0}. The difference between the bounds in
Proposition \ref{p:technicalbound} and those in Theorem
\ref{t:estsc0} is that the constant $\bar{\alpha}$ in Proposition
\ref{p:technicalbound} depends on the spectral gaps for the
manifold $N$. On the other hand,  when $N = \SS^n$ (or a Zoll
surface) we can use the explicit eigenvalue gaps to bound
$(|\bar{\alpha}| - |\alpha|)$ uniformly. Namely, since the $m$-th
cluster of eigenvalues on $\SS^n$ occurs at
\begin{equation}
    b_m = m^2 + (n-1)\, m \, ,
\end{equation}
we get that
\begin{equation}
    b_{m-1} = m^2 + (n-3) \, m +2 -n {\text{ and }} (b_{m-1} + 1)^{1/2} = m +
    \frac{n-3}{2}
    + O(m^{-1})
    \, ,
\end{equation}
where $O(m^{-1})$ denotes a term that is   bounded by   $ C\,
m^{-1}$ for all $m\ne 0$. This gives
\begin{equation}
    4 \, b_m - \left(2 \, (b_{m-1} + 1)^{1/2}  + 1 \right)^2 -   \left(2 \, (b_{m-1} + 1)^{1/2}  + 1 \right)
        = 2 \, m + O(1) \, .
\end{equation}
The key point is that the coefficient of the leading order term is
positive, so there exists some $m_0$ depending only on $\kappa$
and $n$ so that both \eqr{e:chom1} and \eqr{e:chom2} hold for  all
$m \geq m_0$ with
\begin{equation}    \label{e:barapl}
    \bar{\alpha} = 2 \,
(b_{m-1} + 1)^{1/2} + 1  \, .
\end{equation}
Next, let $m_1$ be the smallest positive integer with
\begin{equation}    \label{e:m1}
    2 \,
(b_{m_1-1} + 1)^{1/2} + 1 \geq \alpha    \, .
\end{equation}
Since $(b_{m-1} + 1)^{1/2}$ grows linearly in $m$, there is a
uniform bound for  $2 \, (b_{m_1-1} + 1)^{1/2} + 1 - \alpha$.
Finally, let $m$ be the maximum of $m_0$ and $m_1$ and define
$\bar{\alpha}$ by \eqr{e:barapl}. It follows that we get a uniform
bound for $|\bar{\alpha}| - |\alpha|$ that depends only on
$\kappa$ and $n$.

   A similar argument
  applies for Zoll surfaces.  In particular, this discussion shows
  that
 Proposition \ref{p:technicalbound} gives
  Theorem \ref{t:estsc0}.

\subsection{The frequency function}

The frequency function often gives an alternative approach to
proving a three circles inequality for second order elliptic
equations. This method is predicated  upon having a function whose
hessian is diagonal, such as $|x|^2$ on $\RR^{n+1}$ or the
function $t$ on $N \times \RR$.  However, we will see that this
method does not yield our three circles inequality, but would
instead require some integrability of  $V$ as in \cite{GL}.  For
simplicity and clarity, we will restrict to the case $\SS^1 \times
\RR$.

The frequency function $U(t)$ measures the logarithmic rate of
growth of a function $u$.  Namely,   if we set $J(s) = \int_{\SS^1
\times \{ s \} }   u^2   \, d\theta$, then the frequency is given
by
\begin{equation}
    U(t) = \partial_t \log  J(t) = \frac{J'(t)}{J(t)} \, .
\end{equation}
This is useful because $U(t)$ is a monotone non-decreasing
function of $t$ if $u$ is  harmonic and{\footnote{Equation
\eqr{e:constatinf} rules out functions like the linear function
$t$ where the frequency is {\emph{not}} monotone.}}
\begin{equation}    \label{e:constatinf}
    \lim_{t\to -\infty} \int_{\SS^1 \times \{ t \} }   |\nabla u|^2 \, d\theta  = 0
        \, .
\end{equation}
To see why $U$ is monotone, first differentiate $J$ to get  that
\begin{align}
   J'(s) &= 2 \, \int_{\SS^1 \times \{ s \} }   \left( u u_t  \right) \, d\theta
      \\
    J''(s) &= 2 \, \int_{\SS^1 \times \{ s\} }   \left( u_t^2 +
    u_{\theta}^2
    \right) \, d\theta \, , \label{e:freqjpp}
\end{align}
where the second equation used that $u_{tt} = - u_{\theta \theta}$
and integration by parts on $\SS^1$.
 To get this in a better form,
observe that since $u_{tt} = - u_{\theta \theta}$
\begin{equation}
    \partial_s \int_{\SS^1 \times \{ s \} } (u_t^2 - u_{\theta}^2)
    \, d\theta
    = -2 \, \int_{\SS^1 \times \{ s \} } \partial_{\theta} \, \left( u_t u_{\theta}
     \right) \, d\theta = 0 \, .
\end{equation}
Since we assumed that $\nabla u$ vanishes at $-\infty$ in
\eqr{e:constatinf},  it follows that
\begin{equation}
    \int_{\SS^1 \times \{ s \} } u_t^2 \, d\theta =  \int_{\SS^1 \times \{ s \} }
    u_{\theta}^2 \, d\theta \, .
\end{equation}
Plugging this into the formula \eqr{e:freqjpp} for $J''$ gives
\begin{equation}
    J''(s)  = 4 \, \int_{\SS^1 \times \{ s\} }    u_t^2 \, d\theta
      \, .
\end{equation}
It now follows from the Cauchy-Schwarz inequality that $(J')^2
\leq J \, J''$, so we conclude that $U' = \left[ J'' \, J - (J')^2
\right] / J^2 \geq 0$ as claimed.

Suppose now that $u$ is no longer harmonic, but instead satisfies
the Schr\"odinger equation $\Delta u = - V u$.  In this case, we
get that
\begin{equation}
    J''(s) = 2 \, \int_{\SS^1 \times \{ s\} }   \left( u_t^2 +
    u_{\theta}^2 - V u^2
    \right) \, d\theta \, ,
\end{equation}
 introducing the  ``error term'' $-\int_{\SS^1 \times \{ s \} }
V u^2 \, d\theta $ in $J''(s)$ and giving an  estimate of the form
\begin{equation}    \label{e:fregi}
    \left( \log J \right) '' = U' \geq - C \, \sup |V| \, .
\end{equation}
However, this  lower bound is not integrable in $t$, so $U$ can
decrease by an arbitrarily large amount over a long enough
stretch.  We will see next that this method does not yield the
three circles inequality of Theorem \ref{t:ests}, i.e.,
\begin{equation}    \label{e:ours}
    \log I(t)\leq  C (1+t)   +     \frac{t}{T} \, \left|\log \frac{I(T)}{I(0)}\right| + \log I(0) \, .
\end{equation}
Namely, integrating \eqr{e:fregi} from $s$ to $T$ gives
\begin{equation}
    U(s) \leq U(T) + C \, (T-s) \, ,
\end{equation}
and integrating this from $0$ to $t$ gives
\begin{equation}    \label{e:theirs}
    \log J(t) = \log J(0) + \int_0^t U (s) \, ds
    \leq \log J(0) + t \, (U(T) + C   \, T ) \, .
\end{equation}
To see why \eqr{e:ours} is sharper than \eqr{e:theirs}, suppose
that $|U(T)|$ and $|\log I(T)|/T$ are uniformly bounded but let
$T$ go to infinity. In this case, the upper bound in \eqr{e:ours}
goes to
   $ \log I(0) + C (1+t)$,
whereas the upper bound in \eqr{e:theirs} goes to infinity.

It is interesting to note that N. Garofalo and F.H. Lin,
\cite{GL},   proved unique continuation in a similar setting by
using the frequency function under the stronger assumption that
\begin{equation}
    \int_{\RR} \left( \sup_{\SS^1 \times \{ s\} } |V| \right) \,
    ds < \infty \, .
\end{equation}

\section{Dimension bounds on a manifold with  cylindrical ends}
\label{s:sfour}

In this section, we consider functions $u$ in $H_0$ that solve the
Schr\"odinger equation $\Delta u = -V  u$ for a general bounded
potential $V$ on a manifold $M$ with finitely many cylindrical
ends, each of which is the product of a half-line with a round
sphere or a Zoll surface.{\footnote{A similar argument applies
when the ends have spectral gaps as in \eqr{e:spgaps} and $V$ is
Lipschitz.}} In particular,  $M$ can be decomposed into a bounded
region $\Omega$ together with a finite collection of ends $E_1 ,
\dots , E_k$ where
\begin{itemize}
\item $\Omega$ has compact closure. \item  Each $E_j$ is isometric
to $N_j \times [0,\infty)$, where $N_j$ is either a sphere or a
Zoll surface.
\end{itemize}

The main result of this section is that $H_{\alpha}(M)$ is finite
dimensional for every $\alpha \in \RR$.

\begin{Thm}      \label{t:mlbound}
The linear space $H_{\alpha}$
 has dimension at most
 \begin{equation}
    d = d \left( \alpha , \sup |V|  , \Omega   \right) \, .
\end{equation}
\end{Thm}

Note that F. Hang and F.H. Lin, \cite{HL}, proved a similar result
under the stronger hypothesis that $\sup |V| < \epsilon$ for some
sufficiently small $\epsilon > 0$.

\subsection{A consequence of unique continuation}

We will need an estimate that relates the $W^{1,2}$ norm of a
solution $u$ on the boundary of $\Omega$ to its $W^{1,2}$ norm
inside $\Omega$. This is given  in the next lemma, where we will
use $T_1(\Omega)$ to denote the tubular neighborhood of radius one
about $\Omega$.

\begin{Lem}     \label{l:uc1}
Given $\alpha \geq 0$, there exists a constant $C$ depending on
$\alpha$, $\Omega$, and  $\sup |V|$ so that if $u \in H_{\alpha}$,
then
\begin{equation}    \label{e:uc1}
    \int_{T_1 (\Omega) } \left( u^2 + |\nabla u|^2
    \right) \leq C  \int_{\partial \Omega} \left( u^2 + |\nabla u|^2
    \right) \,  .
\end{equation}
\end{Lem}

\begin{proof}
We will argue by contradiction, so suppose instead that there is a
sequence of functions $u_j$ with $\Delta u_j = - V_j u_j$ where we
have a uniform  bound for   $\sup |V_j|$, each $u_j$ is in
$H_{\alpha}(\Delta + V_j)$, and
\begin{equation}    \label{e:uc2}
    \int_{T_1(\Omega) } \left( u_j^2 + |\nabla u_j|^2
    \right) >  j \,  \int_{\partial \Omega} \left( u_j^2 + |\nabla u_j|^2
    \right) \, .
\end{equation}
The key to the compactness argument is that Corollary \ref{c:ests}
gives a constant $\nu
> 0$  (independent of $j$) so that
\begin{equation}    \label{e:ests}
    \int_{T_1 (\Omega) \setminus \Omega } \left( u_j^2 + |\nabla u_j|^2
    \right) \leq \nu \,  \int_{\partial \Omega} \left( u_j^2 + |\nabla u_j|^2
    \right) \, .
\end{equation}
Therefore, after renormalizing the $u_j$'s, we get that
\begin{equation}    \label{e:uc3}
    \int_{T_1(\Omega) } \left( u_j^2 + |\nabla u_j|^2
    \right) = 1 {\text{ and }} \int_{T_1(\Omega) \setminus \Omega} \left( u_j^2 + |\nabla u_j|^2
    \right) < \nu /j \, .
\end{equation}
The
 interior $W^{2,p}$ estimates (theorem $9.11$ in \cite{GT}) then
 give a
uniform   $W^{2,2}(T_{3/4}(\Omega))$ bound for the $u_j$'s. By
combining this with the Sobolev inequality (theorem $7.26$ in
\cite{GT}), we get uniform higher $L^p$ bounds on the $u_j$'s, and
hence on $V_ju_j$, and then elliptic theory again gives a higher
$W^{2,p}$ bound on $u_j$'s. After repeating this a finite number
of times (depending on $n$), we will get a uniform $W^{2,p}$ bound
for $p > (n+1)$.  Once we have this,  the Sobolev embedding
(theorem $7.26$ in \cite{GT}) gives a uniform $C^{1,\mu}$ bound
\begin{equation}    \label{e:c1mu}
    ||u_j||_{C^{1,\mu}(T_{1,2}(\Omega))} \leq \tilde{C} \, ,
\end{equation}
where $\mu > 0$ and $\tilde{C}$ does not depend on $j$.  We will
refer to this argument as ``bootstrapping.''

It follows from \eqr{e:c1mu} that
 a subsequence  of the $u_j$'s converges uniformly in  $C^1(T_{1/2}(\Omega))$
  to a function $u$  and thus, by \eqr{e:uc3}, $u$ satisfies
\begin{equation}    \label{e:uc4}
    \int_{ \Omega } \left( u^2 + |\nabla u|^2
    \right) = 1 {\text{ and }} \int_{T_{1/2}(\Omega) \setminus \Omega} \left( u^2 + |\nabla u|^2
    \right) = 0 \, .
\end{equation}
We will see that this violates unique continuation of \cite{A}
since $u$ vanishes on an open set, but is not identically zero.
Namely, since the $u_j$'s satisfy
\begin{equation}
    |\Delta u_j| = |V_j| \, | u_j| \leq  \left( \sup_j \sup_M
    |V_j| \right) \, |u_j| \leq C' \, |u_j|   \, ,
\end{equation}
it follows that
\begin{equation}    \label{e:arons}
    |\Delta u| \leq C' \, |u|   \, .
\end{equation}
Finally, the differential inequality \eqr{e:arons} allows us to
directly apply
 \cite{A}.
\end{proof}

\subsection{The proof of Theorem \ref{t:mlbound}}

We will prove Theorem \ref{t:mlbound} by getting an upper bound
for the number of $W^{1,2}(\partial \Omega)$-orthonormal functions
in  $H_{\alpha}(M)$; cf. \cite{CM2}.

\begin{proof}
(of Theorem \ref{t:mlbound}.) Assume that $u_1 , \dots , u_d$ are
functions in $H_{\alpha}(M)$ that are $W^{1,2}(\partial
\Omega)$-orthonormal, i.e., with
\begin{equation}
    \int_{\partial \Omega} \left( u_i u_j + \langle \nabla u_i ,
    \nabla u_j \rangle \right) = \delta_{ij} \, .
\end{equation}
It follows from Corollary \ref{c:ests} that we can find a set of
such functions for any finite $d$ that is less than  or equal to
$\dim (H_{\alpha})$.  Therefore, the theorem will follow from
proving an upper bound on $d$.

Let $U$ denote the vector space spanned by the $u_j$'s and  define
the projection kernel $K(x,y)$ to $U$ on $\partial \Omega \times
\partial \Omega$ by
\begin{equation}
    K(x,y) = \sum_{j=1}^d  \left( u_j(x) \, u_j (y) + \langle \nabla u_j (x) ,
    \nabla u_j (y) \rangle \right) \, .
\end{equation}
Note that
    $\int_{\partial \Omega} K(x,x) = d$.
 We will also need the following  standard estimate for $K(x,x)$
 \begin{equation}   \label{e:varial}
    K(x,x) \leq (n+1) \, \sup_{u \in U\setminus \{ 0\} } \,  \, \frac{  u^2 (x) + |\nabla u|^2
    (x)} {
      \int_{\partial
    \Omega} \left( u^2   + |\nabla u|^2   \right) } \, .
 \end{equation}
To see this, observe first that $K(x,x)$ can be thought of the
trace of a symmetric quadratic form on $U$ and is therefore
independent of the choice of $W^{1,2}(\partial
\Omega)$-orthonormal basis for $U$.  Since the map taking $u \in
U$ to $(u(x) , \nabla u(x))$ is a linear map from $U$ to
$\RR^{n+1}$, we can choose a new $W^{1,2}(\partial
\Omega)$-orthonormal basis $v_1 , \dots , v_d$ for $U$ so that
$v_j (x)$ and $\nabla v_j (x)$ vanish for every $j
> (n+1)$.  Expressing $K(x,x)$ in this new basis gives
\begin{equation}   \label{e:varialproof}
    K(x,x) = \sum_{j=1}^{n+1} \left( v_j^2 (x) + |\nabla v_j|^2
    (x) \right)   \, ,
 \end{equation}
 and   \eqr{e:varial} follows.

We will use \eqr{e:varial} to prove a pointwise estimate for
$K(x,x)$. Namely, Lemma \ref{l:uc1} implies for any $u \in U
\setminus \{ 0 \}$ that
\begin{equation}    \label{e:uc1a}
    \int_{T_1 (\Omega) } \left( u^2 + |\nabla u|^2
    \right) \leq C  \int_{\partial \Omega} \left( u^2 + |\nabla u|^2
    \right) \, ,
\end{equation}
where $C$ depends only on $\alpha$, $\Omega$, and  $\sup |V|$.
Applying the bootstrapping argument of \eqr{e:c1mu} to $u$, i.e.,
 interior $W^{2,p}$ estimates and the Sobolev embedding (theorems $9.11$ and $7.26$
 in \cite{GT}), we get
\begin{equation}
    u^2 (x) + |\nabla u|^2
    (x) \leq C \, \int_{\partial \Omega} \left( u^2 + |\nabla u|^2
    \right) \, ,
\end{equation}
where the new constant $C$ still depends only on $\alpha$,
$\Omega$, and $\sup |V|$.  Substituting this back into
\eqr{e:varial} and integrating gives
\begin{equation}
    d =  \int_{\partial \Omega} K(x,x)  \leq (n+1) \, C \, {\text{Volume}} \,
    (\partial  \Omega) \, ,
\end{equation}
giving the desired upper bound for $d$.
\end{proof}

\section{Density of potentials with $H_0 = \{ 0 \}$}
\label{s:dense}

As in the previous section, we will consider Schr\"odinger
operators $\Delta + V$ where $V$ is a  bounded potential   on a
fixed manifold $M$ with finitely many cylindrical ends, each of
which is the product of a half-line with a round sphere or a Zoll
surface.{\footnote{A similar argument applies when the ends have
spectral gaps as in \eqr{e:spgaps} and $V$ is Lipschitz.}} In
particular, $M$ can be decomposed into a bounded region $\Omega$
together with a finite collection of ends $E_1 , \dots , E_k$
where
\begin{itemize}
\item $\Omega$ has compact closure. \item  Each $E_j$ is isometric
to $N_j \times [0,\infty)$, where $N_j$ is either a sphere or a
Zoll surface.
\end{itemize}

The main result of this section, Proposition \ref{p:dense}, shows
that the set of potentials $V$ where $H_0 = \{ 0 \}$ is dense.

\begin{Pro} \label{p:dense}
Suppose that $L= \Delta + V$ and $f$ is a  non-negative bounded
function with compact support in $\Omega$ that is positive
somewhere. There exists $\epsilon_0 > 0$ so that for all $\epsilon
\in (0,\epsilon_0)$ we have
\begin{equation}
    H_0 (L + \epsilon \, f) = \{ 0 \} \, .
\end{equation}
\end{Pro}

One of the key properties that we will need in the proof is that
if $u$ and $v$ are solutions of $Lu=Lv=0$, then
\begin{equation}        \label{e:dvz}
    \dv \left( u \nabla v - v \nabla u
    \right)= 0 \, .
\end{equation}
Motivated by this, we define the skew-symmetric bilinear form
$\omega ( \cdot , \cdot )$ on functions that are in $L^2 (\partial
\Omega)$ and whose normal derivatives are in $L^2 (\partial
\Omega)$ by setting
\begin{equation}
    \omega ( u, v) = \int_{\partial \Omega} \left( u \partial_n v - v \partial_n u
    \right) \, .
\end{equation}
The next lemma uses Stokes' theorem and \eqr{e:dvz} to prove that
if $u$ and $v$ are solutions of $Lu=Lv=0$
 on $M
\setminus \Omega$ that vanish at infinity, then  $\omega ( u, v) =
0$.

\begin{Lem}     \label{l:stokesp}
If $u$, $v$ are functions on $M \setminus \Omega$ that satisfy
$Lu=Lv=0$ and go to zero at infinity on each end $E_j$, then
\begin{equation}    \label{e:stokesp}
    \omega ( u, v) = 0 \, .
\end{equation}
\end{Lem}

\begin{proof}
 Since $u$ and $v$ go to zero at infinity, the
bootstrapping argument of \eqr{e:c1mu} implies that
 \begin{equation}
    \lim_{t \to \infty} \, \sum_{j=1}^m \int_{N_j\times \{t\}}
\left(u^2 + v^2 +
    |\nabla u|^2 + |\nabla v|^2 \right)
    \big|_{E_j} \, d\theta = 0 \, .
 \end{equation}
The   lemma now follows since Stokes' theorem and \eqr{e:dvz}
imply that {\underline{for any}} $t \geq 0$ we  have
\begin{equation}
    \omega ( u, v) =  \sum_{j=1}^m \, \int_{N_j\times \{t\}}
\left( u \partial_n v - v \partial_n u
    \right)
    \big|_{E_j}\, d\theta  \, .
\end{equation}
\end{proof}

\begin{proof} (of Proposition \ref{p:dense}.)
 Fix a non-negative bounded function $f$ with compact support in
$\Omega$.  We will prove the existence of such a $\epsilon_0 > 0$
by contradiction.  Suppose therefore that there exists a sequence
$\epsilon_j \to 0$ and functions $u_j$ with
\begin{equation}
    u_j \in H_0 (L +
    \epsilon_j \, f) \setminus \{ 0 \} \, .
\end{equation}
  After
dividing each $u_j$ by its $W^{1,2}$ norm on $\partial \Omega$
(this is non-zero by Lemma \ref{l:uc1} and unique continuation,
\cite{A}), we can assume that
\begin{equation}    \label{e:normalize}
    \int_{\partial \Omega} \left( u_j^2 + |\nabla u_j|^2 \right) =
    1 \, .
\end{equation}
Lemma \ref{l:uc1} then gives a constant $C$ depending on $\Omega$,
    $\sup |V|$  and $\sup |f|$ so that
\begin{equation}    \label{e:uc1bad}
    \int_{T_1 (\Omega) } \left( u_j^2 + |\nabla u_j|^2
    \right) \leq C   \, .
\end{equation}
The  bootstrapping argument (i.e.,  $W^{2,p}$ estimates and
Sobolev embedding; cf. \eqr{e:c1mu})  then gives uniform
$C^{1,\mu}$ estimates for the $u_j$'s on the smaller tubular
neighborhood for some $\mu > 0$. Therefore, there is a subsequence
(which we will still denote $u_j$) so that $u_j$ and $\nabla u_j$
  converge uniformly in $T_{1/2}
(\Omega)$ and the limiting function $u$ satisfies the limiting
equation{\footnote{Initially, we only know that $Lu = 0$ weakly,
but elliptic regularity then implies that $u$ is a strong
solution.}}
\begin{equation}
    L u = 0 \, .
\end{equation}
Since $u_j$ and $\nabla u_j$ converge uniformly on $\partial
\Omega$, we get that
\begin{align}
    \int_{\partial \Omega} \left( u^2 + |\nabla u|^2 \right) &
    = \lim_{j \to \infty} \int_{\partial \Omega} \left( u_j^2 + |\nabla u_j|^2 \right) =
    1  \, , \label{e:normone} \\
    \omega ( u , u_k ) &= \lim_{j\to \infty} \, \omega ( u_j , u_k )  = 0 \,
    ,   \label{e:pairz}
\end{align}
where the last equality follows from Lemma \ref{l:stokesp} since
$Lu_j = L u_k = 0$ outside of $\Omega$  and both vanish at
infinity on $E_1 , \dots , E_k$ (recall that $f$ has compact
support in $\Omega$).  Note that \eqr{e:pairz} would have followed
  from Lemma \ref{l:stokesp} alone if we had known that
$u$ also vanishes at infinity.

To get a contradiction, we note that
\begin{equation}
    \dv \, (u \nabla u_k - u_k \nabla u) = - \epsilon_k \, f \,
    u_k \, u \, ,
\end{equation}
so that \eqr{e:pairz} and Stokes' theorem gives
\begin{equation}
    0 = \omega ( u , u_k )  = \int_{\Omega} \dv \, (u \nabla u_k - u_k \nabla u) =  -\epsilon_k \, \int_{\Omega} f\, u \, u_k   \, .
\end{equation}
In particular, since $\epsilon_k > 0$, we must have
\begin{equation}
     \int_{\Omega} f\, u \, u_k = 0 \, .
\end{equation}
Since $u_k \to u$ uniformly in $\Omega$, these integrals converge
to the integral of $fu^2$, so we get that
\begin{equation}
     \int_{\Omega} f\, u^2 = 0 \, .
\end{equation}
Since $f$ is non-negative, but positive somewhere, we conclude
that $u$ vanishes on an open set and, by unique continuation,
\cite{A}, that $u$ is identically zero.  However, this contradicts
\eqr{e:normone}, so we conclude that no such sequence could have
existed.  The proposition follows.
\end{proof}

We can now sum up what we have proven so far as:

\begin{Thm}     \label{t:main2}
Theorem \ref{t:main} holds when each end is {\emph{isometric}} to
a   half-cylinder.
\end{Thm}

Strictly speaking, we have shown the density of potentials where
$H_0 = \{ 0 \}$, but have not yet addressed the density of metrics
where $H_0 = \{ 0 \}$, namely, the more general case of the
theorem that holds for $n=1$. However, this is an easy consequence
of what we have already shown.  To see this, assume that $M$ is
$2$-dimensional and change the metric $g$ conformally by
$\e^{2f}$, where $f$ is bounded, to get an equivalence between
$H_0(\Delta_g + \e^{2f}\, V)$ and $H_0(\Delta_{\e^{2f}g} +  V)$
(see \eqr{e:confH0}). So long as $V$ is not identically zero, this
allows us to perturb the potential to a nearby potential with $H_0
= \{ 0 \}$.  In the remaining case, where $V \equiv 0$ and the
operator is the Laplacian, it follows from  Stokes' theorem that
$H_0 = \{ 0 \}$. Namely, the gradient estimate implies that $|u| +
|\nabla u| \to 0 $ at infinity, so Stokes' theorem gives that
$\int |\nabla u|^2 = 0$ and $u$ must be constant; since the only
constant that goes to zero at infinity is $0$, we get $H_0 = \{ 0
\}$.

\subsection{The set of potentials where $H_0 = \{ 0 \}$ is NOT
open} \label{s:open}

\begin{Exa}
$H_0=\{0\}$ is not an open condition.  Namely, it is easy to
construct a sequence of potentials $V_i$ on $\RR$ with
$|V_i|_{C^1}\to 0$ as $i\to \infty$ and so that $\dim
H_{0}(V_i)>0$ for each $i$ (note that $H_0(0)=\{ 0 \}$; the
limiting Schr\"odinger equation has potential equal to $0$).

To be precise, as we saw right after Theorem \ref{t:main}, it is
easy to construct a potential $V$ on the line (or the cylinder)
such that there exists a solution   to the corresponding
Schr\"odinger equation that goes exponentially to zero at both
plus and minus infinity. (On the cylinder the potential, as well
as the solution, can be taken to be rotationally symmetric, i.e.,
independent of $\theta$.) In fact, the potential can be taken to
be constant (negative) outside a compact set.  Pick such a
potential and name it $V$. On the line look at the rescalings,
$V_{\epsilon}(t)=\epsilon^2\,  V(\epsilon t)$ (on the cylinder
rescale just in the $t$ direction, everything is rotationally
symmetric). Each of the Schr\"odinger equations $u''+V_{\epsilon}
u=0$ has solutions (namely, $u_{\epsilon}(t)=u(\epsilon t)$, where
$u$ is a solution to $u'' + V u = 0$) that decay exponentially to
zero at plus and minus infinity and clearly
$|V_{\epsilon}|_{C^1}\to 0$ as $\epsilon \to 0$. Note also that in
this example $u_{\epsilon}$ converges to the constant function
$u(0)$ as $\epsilon$ goes to zero, which may be taken to be
non-zero and is, of course, in any case, a solution to the
limiting Schr\"odinger equation $u''=0$.
\end{Exa}

\section{Surfaces with more general ends: The case $n=1$ of Theorem
\ref{t:main}}  \label{s:ssix}

We will show in this section that our results apply to a  more
general class of surfaces, namely surfaces with bounded curvature,
injectivity radius, and finitely many ends that are each
bi-Lipschitz to a half-cylinder.   We will prove this by finding a
bi-Lipschitz conformal change of metric on such a surface that
makes each end isometric to a  (flat) half-cylinder and then apply
our earlier results. For this, we will need the following
proposition.

\begin{Pro}     \label{p:}
Suppose that $E$ is topologically a half-cylinder $\SS^1 \times
[0,\infty)$ with a Riemannian metric satisfying the following
bounds:
\begin{enumerate}
\item[(B1)] The Gauss curvature $K_{E}$ is bounded above and below
by $|K_{E}| \leq 1$ and the injectivity radius of $E$ is bounded
below by $i_0 > 0$. \item[(B2)] There is a   bi-Lipschitz
(bijection)
\begin{equation}    \label{e:defFr}
    F = (\theta , r): E \to \SS^1 \times [0,\infty)
\end{equation}
 with   $|dF|\leq \ell_0$ and $|dF^{-1}|\leq \ell_0$.
\end{enumerate}
Then there exists a conformal map $\Phi: E \to \SS^1 \times
[0,\infty)$ satisfying $|d\Phi|_{C^1} \leq C$ and
$|d\Phi^{-1}|_{C^1} \leq C$ for a constant $C$ depending only on
$i_0$ and $\ell_0$.
\end{Pro}

\begin{Rem}     \label{r:norm}
There are two natural equivalent norms for the differential $dF$,
depending on whether one thinks of $dF(x)$ as a vector in $\RR^4$
(Hilbert-Schmidt norm) or as a linear operator on $\RR^2$
(operator norm). We will use the Hilbert-Schmidt norm, but will
often use that $|dF(x)(v)| \leq |dF(x)| \, |v|$. If $\Phi: (\Sigma
, g) \to (\tilde{\Sigma} , \tilde{g})$ is a conformal map between
surfaces, then this convention gives that
\begin{equation}
    \Phi^{\ast}(\tilde{g}) = \frac{1}{2} \, |d\Phi|^2 \, g {\text{
    and }} |d\Phi^{-1}|^2 \circ \Phi = 4 \, |d\Phi|^{-2} \, .
\end{equation}
\end{Rem}

Assuming Proposition \ref{p:} for the moment, we will now complete
the remaining case of Theorem \ref{t:main}  where $n=1$ and the
end merely has bounded geometry and is bi-Lipschitz to a flat
half-cylinder (as opposed to being isometric to it).  Namely, we
prove the following theorem:

\begin{Thm}     \label{t:nequal1}
Let $\Sigma$ be a complete non-compact surface with finitely many
cylindrical ends, each of which  has bounded geometry and is
bi-Lipschitz to a flat half-cylinder.
\begin{enumerate}
\item[(1)] If $V$ is   bounded  (potential) on $\Sigma$, then
$H_{\alpha} (\Sigma, \Delta_{\Sigma} +V)$ is finite dimensional
for every $\alpha$; the bound for   $\dim H_{\alpha}$ depends only
on $\Sigma$, $\alpha$, and $||V||_{L^{\infty}}$. \item[(2)] For a
dense set of
  bounded potentials $H_{0} = \{ 0 \}$.  For each fixed bounded $V$, there is a dense set of
metrics (with finitely many cylindrical ends) where $H_{0}= \{ 0
\}$.
\end{enumerate}
\end{Thm}

\begin{proof}
Using Proposition \ref{p:} we will first show that there exists
another metric $\tilde g$ on $\Sigma$ for which each end is a flat
cylinder and a conformal diffeomorphism $\Phi :(\Sigma,g)\to
(\Sigma,\tilde g)$ with $|d\Phi|_{C^1}\leq C$ and
$|d\Phi^{-1}|_{C^1} \leq C$.  To do this, we first apply
Proposition \ref{p:} to each end $E_j$ of $\Sigma$ to get
conformal diffeomorphisms
\begin{equation}
    \Phi_j : E_j \to \SS^1 \times \RR \, ,
\end{equation}
with a uniform $C^1$ bound for every $j$, i.e., a constant
$\kappa$ so that away from $\partial E_j$ we have
\begin{equation}
    |d\Phi_j| , \, |\nabla d\Phi_j| , \, |d \Phi_j^{-1}| , \, |\nabla d
    \Phi_j^{-1}| \leq \kappa \, .
\end{equation}
Note that the pullback metric $|d \Phi_j|^2 \, g$  makes the end
$E_j$ into a flat cylinder.  It remains to patch these metrics
together across the compact part $\Omega = \Sigma \setminus
\bigcup_j E_j$ of $\Sigma$.  Let $\phi$ be a smooth function on
$\Sigma$ that is identically one on the tubular neighborhood of
radius one about $\Omega$ and has compact support in $\Sigma$ and
then set
\begin{equation}    \label{e:tildeg}
    \tilde g = \left( \phi + (1-\phi) \, \chi_{E_j} \, |d
    \Phi_j|^2 \right) \, g \, .
\end{equation}
Here $\chi_{E_j}$ is the characteristic function of $E_j$.

To see how the operator $\Delta +V$ changes under a conformal
change of metric, let $\Sigma$ be a surface with a Riemannian
metric $g$ and $f$ a smooth function on $\Sigma$. If $u$ is a
solution of the Schr\"odinger equation $\Delta_{g}u+  V\,u=0$ on
$\Sigma$, then $u$ also solves the equation
$\Delta_{\e^{2f}\,g}u+\e^{-2f}\,  V\,u=0$ for the conformally
changed metric.  Namely,
\begin{equation}    \label{e:confch}
\Delta_{\e^{2f}\,g}u=\text{div}_{\e^{2f}\,g}(\nabla_{\e^{2f}\,g}u)
                    =\text{div}_{\e^{2f}\,g}(\e^{-2f}\nabla_{g}u)
                    =\e^{-2f}\Delta_{g}u\, .
\end{equation}
  From this we see that if $g$ and $\tilde g = \e^{2f} \, g$ are
  conformal metrics on a surface $\Sigma$, then
\begin{equation}    \label{e:confH0}
     H_0(\Delta_g + V) = H_0(\Delta_{\e^{2f}g} + \e^{-2f}\,  V)
       \, .
\end{equation}
In particular, applying \eqr{e:confH0} to the surface $\Sigma$
with cylindrical ends with  $\tilde g$ given by \eqr{e:tildeg}, so
that the ends of $(\Sigma , \tilde g)$ are isometric to the
cylinder, we get that $\dim (H_0(\Delta_g + V))$ is finite and is
equal to zero for a generic potential $V$.  More precisely, the
finite dimensionality follows from applying Theorem
\ref{t:mlbound} to the operator $\Delta_{\e^{2f}g} + \e^{-2f}\,
V$, since the Lipschitz bounds on $|d \Phi|$ and $|d\Phi|^{-1}$
imply that
\begin{equation}    \label{e:8p8}
    ||\,|d\Phi|^{-2} \, V||_{C^1(\tilde g)}\leq
C\,||V||_{C^1(g)} \, .
\end{equation}
 Arguing similarly,   the
zero dimensionality for a dense set of  potentials $V$ follows
from \eqr{e:8p8} and the density  for the metric $\tilde g$ proven
in Theorem \ref{t:main2}.
\end{proof}

It remains to  use the bi-Lipschitz map $F$ to find a bi-Lipschitz
conformal map $\Phi$ from each end   to a flat half-cylinder. This
will be accomplished in the next two subsections. The first
subsection proves the existence of a conformal diffeomorphism
$\Phi$ and proves  an $L^{\infty}$ estimate, bounding the second
component of $\Phi$ above and below in terms of $r$ (see
\eqr{e:defFr}). The second subsection proves the uniform Lipschitz
estimates on the conformal factor $|d\Phi|^2$ and its inverse.

\subsection{Uniformization of a cylindrical end}

 The next lemma constructs a harmonic function $u$ on a cylindrical end $E$
 that is bounded above and below by the Lipschitz function $r$,
 i.e., the second component of the map $F$
 given by \eqr{e:defFr}.

\begin{Lem}     \label{l:technical}
There exist constants $C_0 , C_1 > 0$ depending on $i_0$ and
$\ell_0$ so that if $E$ satisfies (B1) and (B2), then there is a
positive harmonic function $u$ on $E$  that vanishes on $\partial
E$ and satisfies
\begin{align}
      C_0^{-1}  &< \int_{\partial E}  \partial_n u  <
     C_0 \, ,  \label{e:p1}  \\
            C_1^{-1} \, r(x)
            &\leq u(x) \leq C_1 \, r(x) {\text{ for }} r(x) \geq 1 \, ,
    \label{e:p2} \\
    0&<|\nabla u| \, . \label{e:p3}
\end{align}
\end{Lem}

\begin{Rem}     \label{r:technical}
It follows that $u$ and its harmonic conjugate $u^{\ast}$ together
give a proper conformal diffeomorphism
\begin{equation}
    (u^{\ast},u): E \to
 \tau \SS^1 \times [0,\infty) \, ,
 \end{equation}
  where the radius $\tau$  is given by
\begin{equation}    \label{e:deftau}
    2 \pi \, \tau = \int_{\partial E} \partial_n u \, .
\end{equation}
To see this, observe that while
 $u^{\ast}$  is not a globally
well-defined function,   it is well-defined up to multiples of
$2\pi \tau$. Hence,  $u^{\ast}$ is a well-defined map to the
circle of radius $\tau$.  Finally, note that \eqr{e:p1} gives
upper and lower bounds for $\tau$.
\end{Rem}

\begin{proof}
(of Lemma \ref{l:technical}). We will construct $u$ as a limit of
harmonic functions $u_j$ on $\{ r\leq j \} \subset E$ as $j\to
\infty$. Namely, let $u_j$ be the harmonic function on $\{ r\leq j
\}$ with $u_j = 0$ on $\partial E$ and $u_j = j$ on the other
boundary component $\{ r = j \}$ (note that $u_j$ exists and is
$C^{2,\alpha}$ to the boundary by standard elliptic theory; cf.
theorem $6.14$ in \cite{GT}).

We will repeatedly use the following consequences of Stokes'
theorem for any $s$ between $0$ and $j$ (Stokes' theorem is used
in the first and last equality below)
\begin{equation}    \label{e:stokesover}
    s \, \int_{\partial E}  \partial_n u_j = s \, \int_{\{ u_j = s
    \}}
    \partial_n u_j = \int_{\{ u_j = s \}} u \,
    \partial_n u_j = \int_{\{ u_j \leq s \}} |\nabla u_j|^2 \, .
\end{equation}

The first step will be to establish the   bounds in \eqr{e:p1} for
the function $u_j$ for a constant $C_0$ that does not depend on
$j$. To get the upper bound, use Stokes' theorem to get
\begin{equation}    \label{e:fluxu1}
    j \, \int_{\partial E}  \partial_n u_j   = \int_{\{ r\leq j \}} |\nabla u_j|^2
    \leq \int_{\{ r\leq j \}} |\nabla r|^2 \leq \ell_0^2 \, \Area \, (\{ r\leq j \}) \leq
    2\pi \,  \ell_0^4
    \, j \, ,
\end{equation}
where the first inequality above uses that $u_j$ and $r$ have the
same boundary values and harmonic functions minimize energy for
their boundary values. (The last two inequalities in
\eqr{e:fluxu1} used the bi-Lipschitz bound for the map $F$. We
will use this again several times in the proof without comment.)
We conclude that
\begin{equation}    \label{e:fluxu}
    \int_{\partial E}  \partial_n u_j \leq 2 \pi \, \ell_0^4
      \, .
\end{equation}
 To get the lower bound, note  that it follows easily from the maximum
 principle that the level set $\{ u_j = s \}$ cannot be
 contractible, so we must have
\begin{equation}    \label{e:llbou}
      i_0  \leq   \Length ( u_j^{-1} (s))
      \, .
\end{equation}
Here, length means the one-dimensional Hausdorff measure if
$u_j^{-1} (s)$ is not a collection of smooth curves.  However,
  we will integrate \eqr{e:llbou} with respect to $s$ and
Sard's theorem implies that almost every level set is smooth, so
this is not an issue.{\footnote{This is really not an issue here
since the argument below for \eqr{e:p3} also implies that $|\nabla
u_j| \ne 0$, so every level set is smooth.}}  Integrating
\eqr{e:llbou} from $0$ to $j$
  and
 using the coarea formula gives
\begin{equation}
     j \, i_0  \leq  \int_0^j \, \Length ( u_j^{-1} (s))
    \, ds
 = \int_{\{ r\leq j \}}  |\nabla u_j|
      \, .
\end{equation}
 Plugging this into Cauchy-Schwarz gives
\begin{equation}
    j^2 \, i_0^2 \leq \left( \int_{\{ r\leq j \}}  |\nabla u_j| \right)^2
    \leq 2 \pi \, \ell_0^2 \, j \,
\int_{\{ r\leq j \}}  |\nabla u_j|^2 = 2 \pi \, \ell_0^2 \, j^2 \,
\int_{\partial E}  \partial_n u_j
      \, ,
\end{equation}
where the last equality comes from applying Stokes' theorem twice,
first to $\dv (u_j \nabla u_j)$ and then to change the boundary
integral of $\partial_n u_j$ on $\{ r = j \}$ into the boundary
integral on $\partial E$.
 We conclude that
\begin{equation}    \label{e:lowerflux}
    \frac{i_0^2}{2 \pi \, \ell_0^2} \leq \int_{\partial E}  \partial_n u_j
      \, .
\end{equation}
Hence, we have uniform upper and lower bounds for the flux of the
$u_j$'s; this will give \eqr{e:p1} in the limit.

We will next establish the bounds \eqr{e:p2} for  the $u_j$'s for
a constant $C_1$ that does not depend on $j$. These uniform
estimates will allow us to extract a limit $u$ that also satisfies
\eqr{e:p2}. We will first show the lower bound in \eqr{e:p2}.  It
will be convenient to let $m_j(s)$ and $M_j (s)$ denote the
minimum and maximum of $u_j$ on $\{ r= s \}$, i.e.,
\begin{align}
    m_j (s) & = \min_{ \{ r = s \} } \, u_j (s) \, , \\
    M_j (s) & = \max_{ \{ r = s \} } \, u_j (s) \, .
\end{align}
 To get the lower bound,  first use Stokes' theorem and the coarea formula
to get that
\begin{equation}
      s  \,   \int_{\partial E} \partial_n u_j = \int_0^s \left( \int_{ \{ r = t \} }
      \partial_n u_j \right) \, dt
       \leq   \int_0^s \, \left( \int_{ \{ r = t \} }
      |\nabla u_j| \right)  \, dt
      = \int_{ \{ r\leq s \} } |\nabla r| \,
      |\nabla u_j| \, .
\end{equation}
 Next, apply
Cauchy-Schwarz to this and then use the bi-Lipschitz bound on $F$
and the fact that $\{ r \leq s \} \subset \{ u_j \leq M_j(s) \}$
to get
\begin{align}
      s^2 \, \left( \int_{\partial E} \partial_n u_j
      \right)^2 &\leq \int_{\{ r\leq s \}} |\nabla r|^2  \, \int_{\{ r\leq s \}} |\nabla u_j|^2
\leq   2 \pi \,  \ell_0^4 \, s \, \int_{ \{ u_j \leq M_j(s)\} }
|\nabla
 u_j|^2 \notag \\
        &=   2 \pi \, \ell_0^4 \, s \, M_j(s) \, \int_{\partial E} \partial_n
        u_j \, ,
\end{align}
where the last equality follows from \eqr{e:stokesover} with
$M_j(s)$ in place of $s$. Combining this with \eqr{e:lowerflux}
gives the desired linear lower bound for $M_j (s)$
\begin{equation}
     s \,  \frac{i_0^2}{2 \pi \, \ell_0^2} \leq
     s \,  \int_{\partial E} \partial_n u_j \leq 2 \pi \, \ell_0^4 \, M_j
     (s) \, .
\end{equation}
To get the upper bound for $m_j (s)$, use the fact that the level
sets of $u_j$ cannot be contractible, see \eqr{e:llbou}, and the
coarea formula to get
\begin{equation}
      i_0  \, m_j  (s) \leq \int_0^{m_j(s)} \, \Length ( u_j^{-1} (t))
    \, dt = \int_{ \{ u_j \leq m_j (s) \} } |\nabla u_j | \, .
\end{equation}
Applying Cauchy-Schwarz and noting that $\{ u_j \leq m_j (s) \}
\subset \{ r \leq s \}$ gives
\begin{align}
     i_0^2 \, m_j^2 (s) & \leq
   \Area \, (\{ r \leq s \})  \, \int_{ \{ u_j \leq m_j (s) \}  } |\nabla
    u_j|^2
     \leq 2 \pi \,
   \ell_0^2 \, s \, \int_{ \{ u_j \leq m_j (s) \} } |\nabla
    u_j|^2  \notag \\
    &= 2 \pi \, \ell_0^2 \, s \,  m_j(s) \, \int_{\partial
    E} \partial_n u_j \leq 4 \pi^2 \, \ell_0^6 \, s \, m_j (s)
  \, ,
\end{align}
where the last inequality uses \eqr{e:fluxu}.  We conclude that
\begin{equation}
    m_j (s) \leq 4 \pi^2 \, s \,  \frac{\ell_0^6}{i_0^2} \, .
\end{equation}
As long as we stay away from
 the
boundary of $\{ r\leq j \}$, we can apply the Harnack inequality
to the positive harmonic function $u_j$.  In particular, the lower
bound for $M_j(s)$ gives a uniform lower bound for $u_j$ and the
upper bound for $m_j(s)$ gives a uniform upper bound for $u_j$.

We will now use these uniform bounds on the $u_j$'s on each
compact set to extract a limit $u$.  Note first that  the upper
bounds for the $u_j$'s in terms of $r$ and standard elliptic
theory (the boundary Schauder estimates; see theorem $6.6$ in
\cite{GT}) give a $C^{2,\alpha}$ bound for the $u_j$'s on each
compact subset of $E$.{\footnote{This bound grows with $r$.  In
the next subsection, we will come back and prove a Lipschitz bound
for $|\nabla u|$ that does not grow with $r$.}} Therefore,
Arzela-Ascoli gives a subsequence of the $u_j$'s that converges
uniformly in $C^2$ on compact subsets of $E$ to a continuous
non-negative harmonic function $u$.  The uniform convergence
implies that $u$ vanishes on $\partial E$ and $u$ also satisfies
\eqr{e:p1} and \eqr{e:p2}; in particular, $u$ is not identically
zero.

We will prove \eqr{e:p3} by contradiction.  Suppose therefore that
$|\nabla u|$ vanishes at some $x \in E$.  Note that $x$ cannot be
in $\partial E$ because of the Hopf boundary point lemma (see
lemma $3.4$ in \cite{GT}). Note also that \eqr{e:p2} implies that
$u$ is proper, so the nodal set $u^{-1} (u (x))$ must be compact.
It follows from the standard structure of the nodal set of a
harmonic function on a surface (see, e.g., lemma $4.28$ in
\cite{CM3}) and the fact that $E$ is a topological cylinder that
there is at least one connected component of $\{ y \, | \, u (y)
\ne u (x) \}$ that both has compact closure and does not contain
$\partial E$ in its boundary. However, this violates the strong
maximum principle and, hence, we conclude that $|\nabla u| \ne 0$.
\end{proof}

\subsection{A uniform Lipschitz bound on the conformal factor}

In this subsection, we will    prove a uniform Lipschitz bound for
$|d \Phi|$ and $|d \Phi^{-1}|$ for any conformal diffeomorphism
$\Phi$ from $E$ satisfying (B1) and (B2) to the flat half-cylinder
$\SS^1 \times [0,\infty) $.  We will then apply this estimate to
the conformal map constructed in the previous subsection to get
Proposition \ref{p:}.    The desired Lipschitz estimate is given
in the next lemma.

\begin{Lem}     \label{l:c1}
There is a constant $\mu$ depending on   $i_0$ and $\ell_0$, so
that if $E$ satisfies (B1) and (B2) and $\Phi : E \to \SS^1 \times
[0,\infty)$ is a conformal diffeomorphism, then away from the
$i_0$-tubular neighborhood of $\partial E$ we get
\begin{equation}
    |d\Phi|_{C^1} \leq \mu {\text{ and }} |d \Phi^{-1}|_{C^1} \leq
    \mu \, .
\end{equation}
\end{Lem}

We will need two preliminary lemmas in the proof of Lemma
\ref{l:c1}. The first is the following result of Bloch (see
Appendix \ref{s:a1}):

\begin{Lem}     \label{l:bloch}
Given $r_0 >0$ and $\kappa > 0$, there exists a constant $B> 0$ so
that if  $\Sigma$ is a surface with $|K| \leq \kappa$, the ball
$B_{r_0}(p)$ is a topological disk in $\Sigma \setminus
\partial \Sigma$,  and
$f$ is a holomorphic function on $B_{r_0}(p)$,  then the image
$f(B_{r_0}(p))$ covers some disk of radius $B \, |df(p)|$ in
$\CC$.
\end{Lem}

We will also need a standard Bochner type formula:

\begin{Lem} \label{l:boch}
If $f$ is a holomorphic function on a surface $E$ and $|\nabla f|
\ne 0$, then
\begin{equation}
    \Delta \log |\nabla f| = K \, .
\end{equation}
\end{Lem}

\begin{proof}
Let  $u$ and $v$ be the real and imaginary parts of $f$, so that
$f= u + i \, v$. The Cauchy-Riemann equations give $|\nabla f| =
\sqrt{2} \, |\nabla u|$ and, hence,
\begin{equation}
    \Delta
\log |\nabla f| = \Delta \log |\nabla u| \, .
\end{equation}
 The
Bochner formula for the harmonic function $u$ gives
 \begin{equation}
    \Delta \log |\nabla u|^2 = 2\, K      + \frac{2 |\Hess_u|^2}{|\nabla u|^2} - \frac{|\nabla |\nabla
    u|^2|^2}{|\nabla u|^4}  \, .
 \end{equation}
    Fixing a point $x$ and working in geodesic normal coordinates that diagonalize the hessian of $u$ at
 $x$ (so that $u_{11}(x) = \lambda = - u_{22}(x)$ and $u_{12}(x) =
 u_{21}(x) = 0$), we get
 \begin{align}
    2 |\Hess_u|^2 \,  |\nabla u|^2 - |\nabla |\nabla
    u|^2|^2&=   4 \lambda^2\, |\nabla u|^2 -  4 \, u_j u_{jk} u_{\ell} u_{\ell k}\notag \\
        &=   4   \, \left( \lambda^2 |\nabla u|^2  -
        \lambda^2 (u_1^2 + u_2^2) \right) = 0 \, ,
 \end{align}
 giving the lemma.
\end{proof}

We can now prove the Lipschitz estimate, i.e., Lemma \ref{l:c1}.

\begin{proof}
(of Lemma \ref{l:c1}.)  We will first use Lemma \ref{l:bloch} to
get the upper bound for $|d\Phi|$.  Given a point $x$ with
$B_{i_0/2}(x) \subset E \setminus \partial E$, let $\tilde \Phi:
B_{i_0/2}(x) \to \CC$ denote the composition of $\Phi$ with the
covering map from the cylinder to $\CC$.  The fact that
$B_{i_0/2}(x)$ is a disk implies that $\tilde \Phi$ is a
well-defined holomorphic function. Furthermore, the fact that
$\Phi$ is injective implies that the projection of $\tilde \Phi
(B_{i_0/2}(x))$ to the cylinder is also an injection. However,
 Lemma
\ref{l:bloch} implies that $\tilde \Phi (B_{i_0/2}(x))$ covers a
disk in $\CC$ of radius
\begin{equation}
    B \, |\nabla \tilde \Phi
(x)| = B \, |d\Phi (x)| \, ,
\end{equation}
 so we must have
\begin{equation}    \label{e:ub}
    B \, |d\Phi (x)| \leq \pi   \, ,
\end{equation}
as desired.

We will next use the upper bound \eqr{e:ub} together with the fact
that $|d\Phi| \ne 0$ to get a lower bound for $\log |d \Phi|$, and
hence an upper bound for $|d \Phi^{-1}|$.  We will use that the
 map $F = (\theta , r)$ maps $E$ to the cylinder
 with bi-Lipschitz constant $\ell_0$.  The key for getting the
 lower bound for $\log |d \Phi|$
 is that the
function $w = \log |d \Phi| = \log |\nabla \tilde \Phi|$ satisfies
$|\Delta w |= |K| \leq 1$ (by Lemma \ref{l:boch})
 and
\begin{equation}    \label{e:forlb}
    c_1 = \log (1/\ell_0) \leq
    \inf_s \, \max_{r = s} w
    \leq \sup
    w \leq \log (  \pi  / B) = c_2 \, .
\end{equation}
 The first inequality in
\eqr{e:forlb} follows from the fact that the curve $\Phi (\{ r =
s\} )$ wraps around the cylinder and hence has length at least
$2\pi$, so that
\begin{equation}
    2 \pi \leq \Length \, (\Phi (\{ r =
s\} )) = \int_{ \{ r = s \} } |d\Phi| \leq
    2\pi \, \ell_0 \, \max_{ \{ r = s \} } \, |d\Phi| \, .
\end{equation}
The last inequality in \eqr{e:forlb} comes from the upper bound
\eqr{e:ub} for $|d\Phi|$. Applying the Harnack
inequality{\footnote{The Harnack inequality that we use here is
that if $w\geq 0$   on $B_{2R}$, then $\sup_{B_R} w \leq C_1 \,
\inf_{B_R} w + C_2 \, \sup_{B_{2R}} |\Delta w|$ where $C_1$ and
$C_2$ depend  on $R$, $\sup |K|$, and $i_0$. Using standard
estimates for the exponential map, it suffices to prove that if
$L$ is a uniformly elliptic second order operator on $B_2 \subset
\RR^2$ and $w\geq 0$ on $B_2$, then $\sup_{B_1} w \leq C_1 \,
\inf_{B_1} w + C_2 \, \sup_{B_{2}} |L\, w|$ where $C_1$ and $C_2$
depend  only on the bounds for the coefficients of $L$. This
follows by combining theorems $9.20$ and $9.22$ in \cite{GT}.}}
 to the
non-negative function $c_2 - w$ centered on a point where $w \geq
c_1$ gives (away from the $i_0$ tubular neighborhood of $\partial
E$)
\begin{equation}    \label{e:lb1}
    \sup (c_2 - w) \leq k_1 \, (c_2 - c_1) + k_2\, \sup |\Delta w| \leq
        k_1 \, (c_2 - c_1) + k_2
     \, ,
\end{equation}
where the constants $k_1$ and $k_2$ depend only on $i_0$ and
$\ell_0$.  Here we have used that every point in $E$ is a bounded
distance from a point where $w \geq c_1$ (by \eqr{e:forlb}) to
ensure that we apply the Harnack inequality on balls of a fixed
bounded size. Rewriting \eqr{e:lb1} gives the desired lower bound
for $|d \Phi|$
\begin{equation}
     \log |d \Phi| = w \geq   c_2 - k_1 \, (c_2 - c_1) - k_2 \, .
\end{equation}

We have now established uniform   bounds on $|d\Phi|$ and
$|d\Phi^{-1}|$.  To get the Lipschitz bounds, we will first work
on the image $\Phi (E) = \SS^1 \times [0,\infty)$ to estimate
$|\nabla_{\SS^1 \times [0,\infty)} \tilde w|$ where
\begin{equation}
    \tilde w  =  \log \, | d\Phi^{-1}| = \log 2 - w \circ \Phi^{-1}  \, ,
\end{equation}
where the last equality used Remark \ref{r:norm}.  The
$L^{\infty}$ estimates above for $w$ imply that $\tilde w$ is
bounded. In addition, applying the formula for the Laplacian of a
conformally changed metric (see \eqr{e:confch}) to $\Delta w = K$
gives
\begin{equation}
    \Delta_{\SS^1 \times [0,\infty)} \tilde w = - \frac{1}{2} \, |d\Phi^{-1}|^2
    \, K \circ \Phi^{-1} \, .
\end{equation}
(The factor $\frac{1}{2}$ comes from our choice of norm; see
Remark \ref{r:norm}.)
 In particular, both $|\tilde w|$ and $|\Delta_{\SS^1 \times
[0,\infty)} \tilde w|$ are uniformly bounded. Therefore, we can
directly apply the Euclidean Cordes-Nirenberg estimate (see, e.g.,
theorem $12.4$ in \cite{GT}) to get
\begin{equation}
    |\tilde w|_{C^{1,\alpha}} \leq C \, (|\tilde w|_{L^{\infty}} +
    |\Delta_{\SS^1 \times [0,\infty)} \tilde w|_{L^{\infty}}) \leq C' \, .
\end{equation}
This gives
 the desired   bound on $|\nabla_{\SS^1 \times [0,\infty)}  d
 \Phi^{-1}|$ and then,
  using the chain rule, it  also gives   the desired bound on
$|\nabla d \Phi|$.
\end{proof}

  Finally, we can combine Lemma \ref{l:technical}  and Lemma
\ref{l:c1} to prove Proposition \ref{p:}.

\begin{proof}
(of Proposition \ref{p:}.) Let $u$ be the positive harmonic
function given by Lemma \ref{l:technical} and let $u^{\ast}$ be
its harmonic conjugate.  As in Remark \ref{r:technical},  we
conclude that  the map
\begin{equation}
    \Phi =  \tau^{-1} \, (  u^{\ast},   \, u) : E \to
  \SS^1 \times [0,\infty)
 \end{equation}
 is the desired conformal diffeomorphism.  Here $\tau$ is defined in \eqr{e:deftau}
    and bounded in \eqr{e:p1}.   The
Lipschitz bounds on $|d\Phi|$ and $|d \Phi^{-1}|$ follow
immediately from Lemma \ref{l:c1}.
\end{proof}

\appendix

\section{Dimension bounds for rotationally symmetric potentials
on cylinders} \label{s:sone}

In this appendix,  we bound the dimension of the space
$H_{\alpha}$ for a rotationally symmetric potential on a flat
cylinder.  In the rotationally symmetric case, things become
particularly simple, but, never the less, it illustrates some of
the ideas needed for the actual argument. We include some simple
ODE comparison results that will also be used within the body of
the paper.

In this  appendix, we will assume that $M$ is a cylinder $N\times
\RR$ with global coordinates $(\theta,t)$ and that the function
$V$ depends only on $t$, i.e., that $V(\theta,t) = V(t)$, and that
$V(t)$ is bounded.

The first result is that the space of functions that vanish at
infinity in the kernel of $\Delta + V$ is  finite dimensional (we
state and prove this only for $H_0$; arguing similarly gives
dimension bounds for any $H_{\alpha}$, where the bound depends
also on $\alpha$):

\begin{Pro}     \label{l:mlbound}
The linear space $H_0$
 has dimension at most
\begin{equation}
    2 \, \big|  \{ j \, | \, \lambda_j \leq \sup V \} \big| \, .
\end{equation}
In particular, when $N =\SS^1$, the dimension is $0$ if $\sup V <
0$ and is bounded by
 $4\, \sqrt{\sup V(t) } + 2$ otherwise.
\end{Pro}

 The key to prove this proposition is that the Fourier coefficients  $[u]_j(t) $ of a solution
$u$, defined in the previous section, satisfy the ODE
\begin{equation}    \label{e:uk}
    w'' (t) = (\lambda_j - V(t)) \, w (t) \, .
\end{equation}

The proposition will follow by first showing that if $u$ is in
$H_+$ and the $j$-th Fourier coefficient for $\lambda_j>\sup V$ is
non-zero{\footnote{The spaces $H_+$ and $H_-$ were defined right
after Theorem \ref{t:main}.}}, then $u$ grows exponentially at
$-\infty$ and likewise for $H_-$. Thus if $u$ lies in $H_0$, so
that it lies in the intersection of $H_+$ and $H_-$, then all
$j$-th Fourier coefficients must be zero for $\lambda_j >\sup V$
and hence $u$ lies in a finite dimensional space. The exponential
growth will follow from Corollary \ref{c:odecomp} below. This
corollary records a consequence of the standard Riccati comparison
argument in a convenient form that will also be needed within the
body of the paper.  The standard proof is included for
completeness.

\begin{Lem}
If $w$ is a function on $[0,\infty)$ that satisfies the ODE
inequality $w''\geq K^2\,w$, $w(0) > 0$, and $w_K$ is a positive
solution to the ODE $w_K''=K^2\,w_K$ with $(\log w)'(0) \geq (\log
w_K)'(0)$, then $w $ is positive  and for all $ t \geq 0$
\begin{equation}    \label{e:comple}
    (\log w)'(t) \geq (\log w_K)'(t)
    \, .
\end{equation}
\end{Lem}

\begin{proof}
  Fix some $b> 0$ so that $w$ is positive on $[0,b)$.  We will show
that \eqr{e:comple} holds for $t \in [0,b)$.  Once we have shown
this, we can integrate \eqr{e:comple} from $0$ to $t$ to get
\begin{equation}
    \log w(t) \geq  \log w(0) + \int_0^t (\log w_K)'(t) \, dt
        = \log w(0) + \log w_K (t) - \log w_K(0) \, ,
\end{equation}
so that $w(b)=\lim_{t\to b}\exp (\log w(t)) > 0$. It follows that
the set $\{ t \, | \, w(t) > 0 \}$ is both open and closed in
$[0,\infty)$, so that $w(t) > 0$ for all $t\geq 0$. Consequently,
\eqr{e:comple}    holds for all $t\geq 0$.

It remains to show that  \eqr{e:comple} holds for $t \in [0,b)$.
To see this,  set $v=(\log w)'$ and $v_K=(\log w_K)'$, so that $v$
and $v_K$ satisfy the Riccati equations
\begin{equation} \label{e:ricatti}
v'+v^2-K^2\geq 0 \text{ and } v_K'+v_K^2-K^2=0\, .
\end{equation}
The claim now follows from the Riccati comparison argument.
Namely, by \eqr{e:ricatti} the function
\begin{equation}
(v-v_K)\exp \left(\int (v+v_K) \right)
\end{equation}
is monotone non-decreasing.
\end{proof}

\begin{Cor}     \label{c:odecomp}
Let $K$ be a positive constant.   Suppose that $w$ satisfies the
ODE inequality
$w''\geq K^2\,w$ and $w(0) > 0$.\\
1).  If $w'(0)\geq 0$ and $w$ is defined on $[0,\infty)$, then
$w(t)\geq  w(0) \, \cosh (Kt)$ for $t\geq 0$.\\
2).  If $w'(0)\leq 0$ and $w$ is defined on $(-\infty,0]$, then
$w(t)\geq w(0) \, \cosh (Kt)$ for $t\leq 0$.

\noindent Moreover, we also have:\\
3).  If $0 > w'(0)> - K \, w(0)$ and $w$ is defined on
$[0,\infty)$,  then for $t\geq 0$ we have
\begin{equation}
    w(t)\geq \frac{ K \, w(0) + w'(0)}{2K} \, \e^{Kt}
+ \frac{ K \, w(0) - w'(0)}{2K} \,
    \e^{-Kt}    \, .
\end{equation}
\end{Cor}

\begin{proof}
If we set $w_K = \cosh (Kt)$, then $w_K'' = K^2 \, w_K$, $w_K$ is
positive everywhere, and $(\log w_K)'(0) = 0$.  The first claim
now follows from the lemma by integrating \eqr{e:comple}. The
second claim follows from applying the first claim to the
``reflected function'' $w(-t)$.

To get the third claim, define the positive function $w_K$ by
\begin{equation}
    w_K = \frac{ K \, w(0) + w'(0)}{2K} \, \e^{Kt} + \frac{ K \, w(0) - w'(0)}{2K} \,
    \e^{-Kt} \, ,
\end{equation}
so that $w_K'' = K^2 \, w_K$, $w_K (0)  = w(0)$, and $w_K'(0) =
w'(0)$.  The last claim now also follows from the lemma by
integrating \eqr{e:comple}.
\end{proof}

\begin{proof}
(of Proposition \ref{l:mlbound}.)
 Suppose that $w$ is solution of   \eqr{e:uk} on
$\RR$ with $\lambda_j
> \sup V$.  If $w$ is not identically zero, then we can apply
either ``1).'' or ``2).'' in  Corollary \ref{c:odecomp} to get
that $w$  grows exponentially at either $+\infty$ or $-\infty$ (or
both).   In particular,  the $j$-th Fourier coefficient $[u]_j (t)
$   of a solution $u \in H_0$ must be zero for every
$\lambda_j>\sup V$.

Since each Fourier coefficient of $u$ satisfies a linear second
order ODE as a function of $t$, it is determined by its value and
first derivative at one point (say $0$).  It follows that any
function $u \in H_0$ is completely determined by the values, and
first derivatives, at $0$ of its $j$-th Fourier coefficients for
$\lambda_j \leq \sup V$.
\end{proof}

\vskip2mm The next corollary is used in Appendix \ref{s:stwo}, but
not in the proof of our main theorem.

\begin{Cor}     \label{c:lft}
If $w(t)$ is a solution of \eqr{e:uk} on $[0,\infty)$  with
$\lambda_j
> \sup V$, then either:
\begin{enumerate}
\item $w(t)$   grows exponentially at $+ \infty$ at least as fast
as $\e^{t \sqrt{\lambda_j - \sup V}}$. \item  $w(t)$   decays
exponentially at $+ \infty$ at least as fast as $\e^{-t
\sqrt{\lambda_j - \sup V}}$.
\end{enumerate}
\end{Cor}

\begin{proof}
It suffices to prove that (2) must hold whenever (1) does not.
Assume therefore that $w(t)$ does not grow exponentially at
$+\infty$. It follows from the first and third claims in Corollary
\ref{c:odecomp} that at any $t$ with $w(t) > 0$ we must have
\begin{equation}
    w'(t) \leq - \sqrt{\lambda_j
- \sup V} \, w(t) \, ,
\end{equation}
since it would otherwise be forced to grow exponentially from $t$
on.  Integrating this gives (2) as long as we know that $w\ne 0$
from some point on.  (If $w<0$ from some point on, then we would
apply the argument to $-w$.)

To complete the proof, recall that $w$ can have only one zero
unless, of course, $w$ vanishes identically.  This follows from
integrating $(ww')' = (w')^2 + (\lambda_j - V)w^2 \geq (w')^2$
between any two zeros.
\end{proof}

\subsection{A geometric example}   \label{s:geomex}

We will next consider an example that illustrates the previous
results.  Namely, consider the rotationally symmetric potential
$V(t)$ on the $2$-dimensional flat cylinder $\SS^1 \times \RR$
\begin{equation}
    V(t) = 2 \, \cosh^{-2} (t) \, .
\end{equation}
Since the potential $V$ is rotationally symmetric, the space of
solutions $u$ of $\Delta u = - V u$ can be written as linear
combinations of separation of variables solutions $w_0 (t)$, $\sin
(k\theta) \, w_k (t)$ and $\cos (k\theta) \, w_k (t)$, where $w_k$
is in the two-dimensional space of solutions to the ODE \eqr{e:uk}
with $\lambda_j = k^2$.  Furthermore, Corollary \ref{c:odecomp}
implies that every $w_k$ with $k^2
> 2 = \sup |V|$ must grow exponentially at
plus or minus infinity.  Hence, to find the space of bounded
solutions, we need only check the solutions of \eqr{e:uk} for
$k=0$ and $k=1$. When $k=0$, we get
\begin{equation}    \label{e:kequal0}
    \frac{\sinh (t)}{\cosh (t)} {\text{ and }}  1 - t \, \frac{\sinh (t)}{\cosh
    (t)}\, ;
\end{equation}
the first is bounded, while the second grows linearly.  When
$k=1$, we get an exponentially growing solution together with an
exponentially decaying solution
\begin{equation}
    \frac{\sinh (2t) + 2t}{\cosh (t)} {\text{ and }} \frac{1}{\cosh (t)}     \, .
\end{equation}
It follows that the space of bounded solutions is spanned by
\begin{equation}
    N_1  = \frac{ \sin (\theta) }{ \cosh (t)  } \, \, \, , \, \, \,
     N_2   = \frac{ -\cos (\theta) }{ \cosh (t)  } \, \, \, ,  \, \, \,
      N_3  = \frac{ \sinh (t) }{ \cosh (t) } \,
     ,
\end{equation}
while the space $H_0$ is spanned by $N_1$ and $N_2$.

This   Schr\"odinger operator arises geometrically as a multiple
of the Jacobi operator (i.e., the second variational operator) on
the catenoid.  The catenoid is the conformal minimal embedding of
the cylinder into $\RR^3$ given by
\begin{equation}
    (\theta ,t) \to (-\cosh t \, \sin \theta , \cosh t \cos \theta , t ) \, .
\end{equation}
It follows that the unit normal is given by
\begin{equation}    \label{e:firstn2cat}
    \nn
        =  \frac{ ( \sin \theta , - \cos \theta  , \sinh t ) }{ \cosh t }
       = (N_1 , N_2 , N_3) \,  ,
\end{equation}
so that $N_1$, $N_2$, and $N_3$ are the Jacobi fields that come
from the coordinate vector fields.  The other (linearly growing)
 $k=0$ solution is the Jacobi field that
comes from dilation.

The above discussion completely determined all polynomially
growing functions in the kernel of the Schr\"odinger
 operator $L= \Delta + 2 \cosh^{-2}(t)$ on the cylinder.  Since
   the kernel of $L$ is
 the $0$ eigenspace of $L$, this leads naturally
 to ask what the entire spectrum of $L$ is.{\footnote{The spectrum of $L$ is
 the set of $\lambda$'s such that $(L+\lambda):W^{2,2} \to L^2$ does not have a bounded
 inverse (note the sign convention);  the simplest way that this can occur is when $\lambda$
 is an eigenvalue of $L$, i.e, when there exists $u_{\lambda} \in
 W^{2,2} \setminus \{ 0 \}$ with $Lu = - \lambda u$.
We refer to Reed and Simon's {\emph{Methods of Modern Mathematical
Physics}}, volumes I through IV, for the definitions and results
in spectral theory  that we will use here.}} We will
    show
  that the  spectrum{\footnote{Note that this is not the same as the spectrum of
 the Jacobi operator on
 the catenoid since the two operators differ by multiplication by
 a positive function (which is why they have the same kernel).}} of $L$ is
 \begin{equation}   \label{e:thespec}
    \{ -1 \} \cup
 [0,\infty) \, .
\end{equation}
   To see this,
 first use Weyl's theorem to see that the essential spectrum is
 $[0,\infty)$
 since the potential $V$ vanishes exponentially on both ends.
   Furthermore, we saw above that
 the positive function
 $\cosh^{-1}(t)$ is an eigenfunction of $L$ with eigenvalue $-1$; this
 positivity implies that $-1$ is the lowest eigenvalue.  It
 remains to show that there is no discrete spectrum between $-1$ and
 $0$.  This will follow from standard spectral theory once we show
 that the constant function $u=0$ is the only polynomially growing solution $u$ of
 \begin{equation}   \label{e:pgs}
    L u = \lambda \, u \, ,
 \end{equation}
 for $0 < \lambda < 1$.  It follows from Corollary \ref{c:lft}
 that such a $u$ must vanish exponentially at both plus and minus
 infinity.  Consequently, every Fourier coefficient $[u]_j$
 is an  exponentially decaying solution of
\begin{equation}   \label{e:fcgs}
    [u]_j''+ 2 \cosh^{-2} (t) \, [u]_j = (\lambda+ \lambda_j) \, [u]_j \,
    ,
 \end{equation}
 where the $\lambda_j$'s are the eigenvalues of $\SS^1$.  In
 particular, since the $\lambda_j$'s are integers and $\lambda$ is not,
 it follows that $\lambda+ \lambda_j \ne 1$.  A standard
 integration by parts argument{\footnote{The exponential decay guarantees that the
 integrals are well-defined and  the boundary
 terms go to zero.}}
  then shows that $[u]_j$ must be
 $L^2(\RR)$-perpendicular to the positive function
 $\cosh^{-2}(t)$; hence, $[u]_j$ must have a zero.  After possibly
 reflecting in $t$, we can assume that $[u]_j(t_0) = 0$ for some
 $t_0 \geq 0$.
 Since
 $\tanh (t) $   satisfies the ODE \eqr{e:fcgs} with $\lambda +
 \lambda_j= 0$ and vanishes only at $0$, the lowest eigenvalue of the operator
 $\partial_t^2 + 2 \cosh^{-2} (t)$ on any subdomain of the
 half-line $[0,\infty)$ must be non-negative.  We will use two
 consequences of this:
 \begin{itemize}
 \item  First, $[u]_j  (t)$ cannot vanish for $t>
 t_0$ unless it vanishes identically; suppose therefore that $[u]_j (t) \geq 0$ for $t\geq t_0$.
 \item
 Second,
 the solution $w$ of the ODE \eqr{e:fcgs} with $\lambda +
 \lambda_j= 0$ and initial values $w(t_0) = 0$ and $w'(t_0) = 1$
 must be positive for all $t\geq t_0$.
 \end{itemize}
 Note that we have already
 shown in \eqr{e:kequal0} that any such $w$ grows at most linearly in $t$.  Hence,
 since $[u]_j$ vanishes exponentially, we know that
 \begin{equation}
    \lim_{t \to \infty} \left[ w [u]_j' - w' [u]_j \right](t) \to 0
    \,
\end{equation}
Since $\left[ w [u]_j' - w' [u]_j \right]' = (\lambda+ \lambda_j)
\, w [u]_j$, the fundamental theorem of calculus gives that
\begin{equation}
    (\lambda+ \lambda_j) \, \int_{t_0}^{\infty} w(t) [u]_j (t) \,
    dt = 0 \, ,
\end{equation}
where we also used that $w(t_0) = [u]_j (t_0) = 0$.  Since $w>0$
and $[u]_j \geq 0$ on $[t_0 , \infty)$, we conclude that $[u]_j$
vanishes identically as claimed, completing the proof of
\eqr{e:thespec}.

\section{Growth and decay for generic rotationally symmetric and periodic
potentials}   \label{s:stwo}

In this appendix,  we introduce the Poincar\'e map and use it to
make some remarks about decay and growth for a generic
rotationally symmetric potential on a cylinder; these are not
needed elsewhere (nor are the results of Appendix
\ref{s:symplectic}), but are included for completeness. This
shows, in particular, that for an open and dense set of periodic
potentials with positive operators any solution that vanishes at
infinity decays exponentially.

  For a bounded and rotationally symmetric potential on a  cylinder
$N\times\RR$, Appendix \ref{s:sone} applies to all but a finite
number of small eigenvalues of $\Delta_N$.  To understand the
remaining small eigenvalues of $\Delta_N$, we will need to
understand the ``Poincar\'e maps'' associated to the ODE.  We will
define this next.  For simplicity, we will assume throughout both
this appendix and the next that $V$ is smooth.

  Given a non-negative number $\lambda$ and $t_1 \leq t_2$, define the
{\emph{Poincar\'e map}} $P^\lambda_{t_1,t_2}:\RR^2 \to \RR^2$ by
\begin{equation}    \label{e:defpoin}
    P^\lambda_{t_1,t_2} (a,b) = (u (t_2) , u' (t_2)) {\text{ where }}  u''
    = (\lambda- V ) \, u {\text{ and }} u(t_1) = a , \, u'(t_1) = b \, .
\end{equation}
In general if $f$, $g$ are functions on $\RR$, not necessarily
periodic, satisfying  $f''= (\lambda- V)f$ and $g''= (\lambda- V)
g$, then
\begin{equation}
\frac{d}{dt} \, {\text{det}}   \left( \begin{array}{cc}
f & g \\
f' & g' \end{array} \right) = 0 \, .
\end{equation}
It follows from this and the fact that $P_{t,t}$ is the identity
that $P_{t,t+s}$ is in $SL(2,\RR)$ for all $t$ and $s\geq 0$.

We will below combine this with the simple fact that if $A$ is a
matrix in $SL(2,\RR)$, then either
\begin{enumerate}
\item The absolute value of the trace of $A$ is (strictly) greater
than two, so the characteristic polynomial of $A$ has two distinct
real roots, $c \in \RR$ and $1/c$ where $|c|>1$. Such an $A$ is
said to be {\emph{hyperbolic}} and can be diagonalized even over
$\RR$. \item The absolute value of the trace of $A$ is (strictly)
less than two, so the characteristic polynomial of $A$ has two
distinct complex roots, $\e^{i\,\phi}$ and $\e^{-i\,\phi}$ where
$0 < \phi < \pi$. \item The absolute value of the trace of $A$ is
equal to two, so the characteristic polynomial of $A$ has the root
one, or the root minus one,   with multiplicity two. Thus, there
exists an orthonormal basis where $A$ can be represented by (plus
or minus)
\begin{equation}
   \left( \begin{array}{cc}
1 & \lambda \\
0 & 1 \end{array} \right)  \, .
\end{equation}
\end{enumerate}

\begin{Lem} \label{l:return}
  For an open dense set of  potentials $V$ on $[0,\ell]$,  the
absolute value of the trace of the Poincar\'e map $P_{0,\ell}$ is
not equal to two.  In fact, if $V$ is a potential with   $\left|
{\text{Trace}} (P_{0,\ell}) \right| = 2$, then there are
potentials $V_j \to V$ with $V_j(0) = V(0)$ and $ V_j(\ell) =
V(\ell)$ so $P^{V_j}_{0,\ell} $ is hyperbolic.
\end{Lem}

 To prove the lemma observe first that since  trace is continuous on
  $SL(2,\RR)$ and the
Poincar\'e map $P_{0,\ell}$ depends continuously on the potential,
the set of potentials where the absolute value of the trace of
$P_{0,\ell}$ is not two is clearly open. Consequently, to prove
Lemma \ref{l:return}, it is enough to prove density. The density
is an easy consequence of the next lemma that allows us to perturb
the Poincar\'e map.

We will need a few definitions before stating this perturbation
lemma. Namely, given $\ell > 0$ and a function $f$ on $[0,\ell]$
with $f(0)=f(\ell)$, let $P(f,s) = P_{0,\ell} (f,s)$ denote the
Poincar\'e map from $0$ to $\ell$ for the perturbed operator
$\partial_t^2 + (V(t) + s \, f(t))$.{\footnote{Note that the
perturbed potential agrees at $0$ and $\ell$ if the original
potential $V$ does.}}

\begin{Lem}     \label{l:push}
The  linear map from functions on $[0,\ell]$ with $f(0) =
f(\ell)=0$ to $2\times2$
 matrices given by
\begin{equation}
    f \to \frac{d}{ds} \, \big|_{s=0} \, P(f,s)
\end{equation}
is onto the three-dimensional space of matrices $B$ such that
$P^{-1}(f,0)\, B$ is trace free.
\end{Lem}

\begin{proof}
Let $u(s,t)$ and $v(s,t)$ be the solutions of $\partial_t^2 +
(V(t) + s \, f(t))$ with initial conditions $(u, u_t)(s,0) =
(1,0)$ and $(v, v_t)(s,0) = (0,1)$.  It follows that
\begin{equation}
\frac{d}{ds} \, \big|_{s=0} \, P(f,s) =   \left( \begin{array}{cc}
u_s & v_s \\
u_{ts} & v_{ts} \end{array} \right) (0,\ell) \, .
\end{equation}
Note that when $s=0$, we have $u_{tt} = -Vu$, $v_{tt} = - Vv$,
$u_{stt} = - Vu_s - fu$, and  $v_{stt} = - Vv_s - fv$.  It follows
that
\begin{align}
(u_s v_t - v u_{st})(0,\ell) =& \int_0^{\ell} \partial_t \, (u_s
v_t - v u_{st}) (0,t) \, dt
        = \int_0^{\ell} (f uv)(0,t) \, dt \, , \\
   (u u_{st} - u_s u_t  )(0,\ell) =& \int_0^{\ell}
\partial_t \, (u u_{st} - u_s u_t  ) (0,t) \, dt
        = -\int_0^{\ell} (f u^2)(0,t) \, dt \, ,  \\
(v_s v_t - v v_{st})(0,\ell) =& \int_0^{\ell} \partial_t \, (v_s
v_t - v v_{st}) (0,t) \, dt
        =  \int_0^{\ell} (f v^2)(0,t) \, dt \, , \\
(u v_{st} - v_s u_t )(0,\ell) =& \int_0^{\ell} \partial_t \, (u
v_{st} - v_s u_t ) (0,t) \, dt
        = - \int_0^{\ell} (f uv)(0,t) \, dt \, .
\end{align}
These four quantities are the $11$, $12$, $21$, and $22$
coefficients, respectively, in the $2\times2$ matrix obtained by
multiplying $P^{-1}(0,f)$ by $\frac{d}{ds} \, \big|_{s=0} \,
P(f,s)$.  Since $u^2(0, \cdot)$, $v^2(0, \cdot)$, and $(uv)(0,
\cdot)$ are linearly independent{\footnote{To see this, note that
the $3\times3$ matrix whose columns are $(u^2 (0,0) , (u^2)_t(0,0)
, (u^2)_{tt} (0,0))$, and similarly for $v^2$ and $uv$, is
invertible.}}  as functions on $[0,\ell]$ and composition by a
linear map can only decrease the dimension of a vector space, we
conclude that the image of $\frac{d}{ds} \, \big|_{s=0} \, P(f,s)$
must be at least three-dimensional. Finally, since it is contained
in a three-dimensional space of matrices, the mapping must be
onto.
\end{proof}

\begin{proof}
(of Lemma \ref{l:return}.) As noted right after the statement of
Lemma \ref{l:return}, it is enough to show that if the absolute
value of the trace of $P_{0,\ell}$ is two, then there is some
function $f$ with $f(0)=f(\ell)=0$ so that for all $s> 0$
sufficiently small we have
\begin{equation}
     \left| {\text{Trace}}
    P_{0,\ell}^{V+sf} \right| > 2 \, ,
\end{equation}
where $P_{0,\ell}^{V+sf}$ is the Poincar\'e map for the potential
$V + s f$.  This follows immediately from two facts.  First, Lemma
\ref{l:push} says that we can choose $f$ to arbitrarily perturb
$P_{0,\ell}$ in $SL(2,\RR)$.  Second, if $P$ is a matrix in
$SL(2,\RR)$ with $\left| {\text{Trace}}
    P \right| = 2$, then there are
matrices $P_j \in SL(2,\RR)$ converging to $P$ with $\left|
{\text{Trace}}
    P_j \right| > 2$. Namely, if we consider
$SL(2,\RR)$ as the hyper-surface $x_1x_4 - x_2 x_3 = 1$ in
$\RR^4$, then the normal direction and the gradient of trace are
\begin{align}
    N &= (x_4, -x_3 , - x_2 , x_1) \, , \\
    \nabla {\text{Trace}} &= (1, 0 , 0 , 1) \, .
\end{align}
In particular, the  projection of $\nabla {\text{Trace}}$ to the
tangent space of $SL(2,\RR)$ vanishes only at the identity matrix
and minus the identity matrix. It follows that we can perturb the
trace as desired, at least away from (plus or minus) the identity
matrix.  This is all that we need, since it is easy to perturb
(plus or minus) the identity matrix
 to a hyperbolic matrix.
\end{proof}

In the next corollary, we will assume that the potential $V$ is
both periodic at $+\infty$ with period $\ell_+$ and that the
associated operator $\partial_t^2 + V$ is positive at infinity.
That is, we will assume that there exists some $T> 0$ so that:
\begin{itemize} \item For  all $t>T$, we have that
$V(t+\ell_+)=V(t)$. \item The only solution  of $u''= - V u$ with
at least two zeros on $[T,\infty)$ is the constant zero.
\end{itemize}
Note that the second condition is equivalent to the lowest
eigenvalue{\footnote{By convention, $\lambda$ is an eigenvalue of
$\partial_t^2 + V$ on $[a,b]$ if there is a (not identically zero)
solution $u$ of $u'' + V u = - \lambda u$ with $u(a)=u(b) = 0$.}}
 being positive on every compact subinterval of
$[T,\infty)$; this follows from the domain monotonicity of
eigenvalues.

\begin{Cor} \label{c:baire}
  For an open and dense set of $\ell_+$ periodic at $+\infty$
potentials $V$ on $\RR$ that are also positive at infinity, any
solution $u \in H_+$ to the Schr\"odinger equation $\Delta u = -
V(t) \, u$ on the cylinder $N \times \RR$
 must decay exponentially to $0$ at $+\infty$.  Likewise for
$H_-$.
\end{Cor}

\begin{proof}
We will use the positivity of $\partial_t^2 +V$, and hence also of
$\partial_t^2 +V - \lambda_j$, to show first that the eigenvalues
of $P^{\lambda_j}_{T,T+\ell_+}$ must be real for every $j$.
  To see this,
suppose instead that the eigenvalues  are $\e^{i\phi}$ and
$-\e^{i\phi}$ with $0< \phi < \pi$.  In this case, we can choose
some positive integer $n$ to make both of the eigenvalues of
\begin{equation}
    \left(P^{\lambda_j}_{T,T+\ell_+}\right)^n = P^{\lambda_j}_{T,T+n\,\ell_+}
\end{equation}
as close as we want to $-1$.  In particular, the solution of $f''=
  ({\lambda_j} - V ) f$ on $[T,T+2n\, \ell_+]$ with initial values $f(T) =
1$ and $f'(T)= 0$ must be negative at $T+n\ell_+$ and then
positive again at $T+2n\ell_+$.  This contradicts the positivity
of the operator, so we conclude that the eigenvalues must be real.

Applying Lemma \ref{l:return}, we may assume that
$P^{\lambda_j}_{T,T+\ell_+}$ has two distinct real eigenvalues for
${\lambda_j} \leq \sup |V|$ and hence is hyperbolic. To complete
the proof, we will prove the exponential decay of any $u \in H_+$
for such a potential.

 By expanding a solution $u$ into its Fourier series, it suffices
to prove   a uniform rate of exponential decay for bounded
solutions $f$ of $f''= ({\lambda_j}- V)f$ on $[0,\infty)$.  Here
``uniform'' means independent of $j$. Corollary \ref{c:lft} gives
this uniform exponential decay for every $j$ with ${\lambda_j}
> \sup |V|$; this does not use the periodicity at $+\infty$.

Assume now that ${\lambda_j} \leq \sup |V|$. It remains to show
that if $P^{\lambda_j}_{T,T+\ell_+}$ is hyperbolic, $f$ vanishes
at $+\infty$, and $f''= ({\lambda_j}- V)f$ on $[0,\infty)$, then
$f$ decays exponentially to zero at $+\infty$.
 For simplicity, we will assume that $j=0$ and $T=0$. The argument
in the general case follows with obvious modifications. Let $v_1$
and $v_2$ be the two eigenvectors of the Poincar\'e map
$P_{T,T+\ell_+}$ such that $v_1$ corresponds to the eigenvalue
with norm larger than one.  Let $f_1$ and $f_2$ be solutions on
$\RR$ to the equation $f'' = - V(t) \, f$ defined by
$(f_i(0),f'_i(0))=v_i$. It follows, that $f_1$ grows exponentially
at $+\infty$ and $f_2$ decays exponentially to $0$ at $+\infty$.
Moreover, since the space of solutions is two dimensional and
$f_1$ and $f_2$ are clearly linearly independent any solution $f$
can be written as a linear combination of $f_1$ and $f_2$.  Thus
$f=a\,f_1+b\,f_2$ for constants $a$ and $b$. It follows that if
$f$ vanishes at $+\infty$, then $a=0$ and hence $f$ must decay
exponentially at $+\infty$.
\end{proof}

\begin{Exa}
We will compute the Poincar\'e maps in the geometric example from
subsection \ref{s:geomex}, where $V(t)$ is a
  rotationally symmetric potential
 on   $\SS^1 \times \RR$ given by
\begin{equation}
    V(t) = 2 \, \cosh^{-2} (t) \, .
\end{equation}
 Using the solutions in \eqr{e:kequal0}, we
get that
\begin{equation}
P_{0,t} =    \left( \begin{array}{cc}
1 -t  \tanh t  & \tanh t \\
-\tanh t -  t\, \cosh^{-2} t  & \cosh^{-2} t
\end{array} \right) \, .
\end{equation}
It follows that if $s < 0 < t$ are large, then $P_{s,t} = P_{0,t}
\circ P_{0,s}^{-1}$ is approximately given by
\begin{equation}
   \left( \begin{array}{cc}
-1 & 1 + s - t \\
0 & -1
\end{array} \right) \, .
\end{equation}
Here, ``approximately'' means up to terms that decay exponentially
in $s$ or $t$.
\end{Exa}

\section{The symplectic form and the symplectic Poincar\'e maps}
\label{s:symplectic}

Much of the discussion of the previous appendix generalizes to
general   bounded potentials that are no longer assumed to be
rotationally symmetric. To explain this, we will need to recall
some standard definitions. Let $\cH$ be a Hilbert space with inner
product $\langle \cdot,\cdot \rangle$ and let $\omega$ be the
canonical symplectic form on $\cH^2=\cH\times \cH$. That is, if
$(v_1,v_2)$ and $(w_1,w_2)$ are in $\cH^2$, then $\omega
((v_1,v_2),(w_1,w_2))=\langle v_1,w_2\rangle -\langle
v_2,w_1\rangle$. (The skew symmetric $2$-form $\omega$ is
symplectic since it is non-degenerate.) By definition a linear map
from $\cH^2$ to itself is said to be a {\emph{symplectic map}} if
it preserves $\omega$.  A linear subspace of $\cH^2$ is said to be
a {\emph{symplectic subspace}} if $\omega$ restricted to the
subspace is non-degenerate and a linear subspace of a symplectic
subspace is said to be {\emph{isotropic}} if the restriction of
the  symplectic form vanishes on the subspace. An isotropic
subspace is said to be {\emph{Lagrangian}} when it is maximal,
i.e., is not strictly contained in a larger isotropic subspace.
Finally, if $W$ is a finite dimensional symplectic subspace of
dimension $2n$, then it follows from Darboux's theorem that
$\omega^n$ is a volume form on $W$ and thus if $A$ is a symplectic
map from $W$ to $W$, then $A$ also preserves the volume form.

Consider now again   solutions $u$  of $\Delta u = - Vu$ on the
half-cylinder $N \times [0,\infty)$.   The potential $V$ will be
smooth and bounded, and is now also allowed to depend on
$\theta\in N$.

The Hilbert space will be $L^2(N)$ with the usual inner product
whose norm is the $L^2$ norm. This is because the next lemma will
allow us to identify a solution of $\Delta u = - V u$ with its
Cauchy data $(u, \partial_t u)$ on a slice $ N \times \{ t_0\}$.

\begin{Lem}  \label{l:uniqueness}
If  $u (\cdot , t_0) =\partial_t u (\cdot , t_0) = 0$, then $u$ is
identically zero.
\end{Lem}

\begin{proof}
We will show first that $u$ and all of its derivatives vanish on
  $N \times \{ t_0 \}$.  Since   $u$   vanishes on $N \times \{ t_0 \}$
  and $\partial_t$   commutes with $\nabla_N$,  every
 partial derivative with at least one derivative in a direction
tangent to $N$   also
 vanishes.  It remains to check that $\partial^n_{t} u ( \cdot ,t_0)$ vanishes
 for all $n \geq 2$.  To get this for $n=2$, use the equation $\Delta u = - V
 u$ to write
 \begin{equation}
    \partial_t^2 u (\theta, t_0) = - \Delta_{N}  u (
    \theta,t_0) - V ( \theta, t_0) \, u (\theta, t_0) = 0 \, .
 \end{equation}
 Similarly, differentiating the equation gives for $n > 2$ that
 \begin{equation}
    \partial_t^n u (\theta, t_0) = - \partial_t^{n-2} \,
\Delta_N  u (
    \theta,t_0) - \partial_t^{n-2} \,  \left[ V (\theta, t_0)
\, u (\theta, t_0) \right] \, .
 \end{equation}
By induction, the  terms on the right hand side of the equation
all vanish, so we conclude that $\partial^n_{t} u ( \cdot, t_0 )$
also vanishes
 for all $n \geq 2$.

  Finally,  since the potential $V$ is bounded and $u$ vanishes to
infinite order on   $N \times \{ t_0 \}$, it follows from the
theory of unique continuation, \cite{A}, that $u$ must vanish
everywhere.
\end{proof}

We conclude from the lemma that the linear map that takes a
solution $u$ of $\Delta u=-V\, u$ to its Cauchy data
$(u,\partial_t u)$ is injective and hence we can identify $u$ with
its Cauchy data on an arbitrary but fixed slice $N \times \{ t_0
\}$.

Motivated by this, we define the skew-symmetric bilinear form
$\omega ( \cdot , \cdot )$ on solutions by
\begin{equation}
    \omega (u, v) = \int_{N \times \{ t_0\} }
\left( u \partial_n v - v \partial_n u \right) \, .
\end{equation}
Thus under the Cauchy data identification the space of solutions
is identified with (a subspace of) $L^2(N)\times L^2(N)$ and the
skew symmetric bilinear form is the pull back of the canonical
symplectic form on $L^2(N)\times L^2(N)$.

As an immediate consequence of  Stokes' theorem and  that $\dv
\left( u \nabla v - v \nabla u
    \right)= 0$ if $Lu=Lv=0$,
the next lemma shows that the skew symmetric form $\omega$ does
not depend on the choice of slice $N \times \{ t_0\}$.

\begin{Lem}     \label{l:stokespa}
If $u$, $v$ are functions on $N \times [t_0,t_1]$ that satisfy
$Lu=Lv=0$, then
\begin{equation}    \label{e:stokespa}
    \omega ( u, v) \equiv \int_{N \times \{ t_0\}}
\left( u \partial_n v - v \partial_n u \right) = \int_{N \times \{
t_1 \}} \left( u \partial_n v - v \partial_n u \right) \, .
\end{equation}
\end{Lem}

As an immediate consequence of Lemma \ref{l:stokespa}, we get that
$\omega$ vanishes on the space of solutions of $Lu=0$ that vanish
at $+\infty$, i.e., the image of the map from $H_0$ to its Cauchy
data is an isotropic subspace; cf. Lemma \ref{l:stokesp}.

\subsection{The Poincar\'e map}
In the spirit of the previous appendix, we can use solutions of
the equation $Lu=0$ to define a Poincar\'e map which maps the
Cauchy data at one time to the Cauchy data (of the same solution)
at a later time.  Namely,
  given  $t_1 \leq t_2$, define the
{\emph{Poincar\'e map}}
\begin{equation}
    P_{t_1,t_2}: L^2(N)\times L^2(N) \to
L^2(N)\times L^2(N)
\end{equation}
 by
\begin{equation}
    P_{t_1,t_2} (f,g) = (u (\cdot , t_2) ,
\partial_t u (\cdot, t_2)) {\text{ where }}
    Lu=0 , \,
        u(\cdot , t_1) = f , \, \partial_t u(\cdot , t_1) = g \, .
\end{equation}
Lemma \ref{l:stokespa} then says that the linear map $P_{t_1,t_2}$
preserves the skew symmetric form $\omega$:

\begin{Cor}     \label{c:symp}
The linear Poincar\'e map $P_{t_1,t_2}$ is symplectic, i.e.,
\begin{equation}
 \omega (f,g) = \omega (P_{t_1 , t_2}(f) , P_{t_1 , t_2}(g) ) \, .
 \end{equation}
\end{Cor}

\subsection{Perturbing the Poincar\'e map}
We will now consider a one-parameter family of Schr\"odinger
operators
\begin{equation}
    L + s f = \Delta + V(\theta , t) + s \, f(\theta , t) \, ,
\end{equation}
together with the associated one-parameter family of Poincar\'e
maps
\begin{equation}
    P = P_{0, \ell } (L+ s f) : L^2 (N) \times L^2 (N) \to
    L^2 (N) \times L^2 (N) \, .
\end{equation}
The next lemma will allow us to compute the derivative with
respect to $s$ of the Poincar\'e map $P$.  In order to state the
lemma, it will be convenient to define a map
\begin{equation}
    L^{-1}: L^2 (N) \times L^2 (N) \to L^2
(N \times [0,\ell])
\end{equation}
 which takes a pair of
functions $(f_1 , f_2)$ to the solution $u$ of $Lu=0$ with Cauchy
data $u(\cdot , 0 ) = f_1$ and $\partial_t u(\cdot , 0) = f_2$.
(Note that this is not defined for all $f_1$ and $f_2$ since the
Cauchy problem in not solvable in general for elliptic equations.)

\begin{Lem}     \label{l:pdiff}
Given functions $(f_1 , f_2)$ and $(g_1 , g_2)$ in $L^2 (N) \times
L^2 (N)$, we have
\begin{equation}
    \omega \left( P(f_1 , f_2) ,  \frac{d}{ds} \, \big|_{s=0}
    P(g_1 , g_2) \right) = - \int_{N \times (0,\ell) } f \, L^{-1}(f_1 , f_2) ,
        \, L^{-1}(g_1 , g_2) \, .
\end{equation}
\end{Lem}

\begin{proof}
Let $u(s,t,\theta)$ and $v(s,t,\theta)$ be solutions of
\begin{equation}    \label{e:equv}
    (L+sf) u = (L+sf) v = 0 \, ,
\end{equation}
with initial conditions{\footnote{We must {\emph{assume}} that
these exist, since the Cauchy problem in not generally solvable.}}
\begin{equation}
    (u , \partial_t u)(s , 0 , \theta) = (f_1 , f_2) (\theta)
    {\text{ and }}
    (v , \partial_t v)(s , 0 , \theta) = (g_1 , g_2) (\theta) \, .
\end{equation}
It follows that
\begin{equation}
    \omega \left( P(f_1 , f_2) ,  \frac{d}{ds} \, \big|_{s=0}
    P(g_1 , g_2) \right) =  \omega
    \left( \, (u , \partial_t u)(0,\ell , \cdot) , \,
    (v_s , \partial_t v_s) (0,\ell , \cdot) \, \right)  \ \, .
\end{equation}
The equation \eqr{e:equv}, and its $s$-derivative, implies that
$Lu(0,\cdot , \cdot) = 0$ and \begin{equation}
    Lv_s (0,\cdot ,
\cdot) = -f (\cdot , \cdot) \, v (0,\cdot , \cdot) \, .
\end{equation}
 In particular,
\begin{equation}        \label{e:dvzb}
    \dv \left( u \nabla v_s - v_s \nabla u
    \right) (0, \cdot , \cdot) = - (fuv) (0, \cdot , \cdot)  \, .
\end{equation}

 Since $(v_s , \partial_t v_s) (0, 0 , \cdot) = (0,0)$, we can
use Stokes' theorem and \eqr{e:dvzb} to get
\begin{equation}
     \omega
    \left( \, (u , \partial_t u)(0,\ell , \cdot) , \,
    (v_s , \partial_t v_s) (0,\ell , \cdot) \, \right) =
    - \int_{N \times (0,\ell) } (fuv) (0, \cdot , \cdot) \,
     \ \, .
\end{equation}
\end{proof}

Note that since $\omega$ is a non-degenerate form, Lemma
\ref{l:pdiff} completely determines the mapping $\frac{d}{ds} \,
\big|_{s=0}
    P(g_1 , g_2)$. Finally, observe that if the potential $V$ is $\ell$-periodic on $N \times
[0,\infty)$, then the map $P_{\ell,0}$ maps the Cauchy data of
$H_{\alpha}$ into itself.

\section{Bloch's theorem}   \label{s:a1}

The classical Bloch theorem is usually stated for a disk in $\CC$.
We need a version of Bloch's theorem for a topological disk in a
surface with bounded curvature.  Since we were unable to find an
 exact reference for this, we will explain here how to get the needed version.
 The following lemma is an immediate consequence of the
classical Bloch theorem (see \cite{Ah}):

\begin{Lem}     \label{l:bloch0}
There exists a constant $B_0> 0$, so that if $f$ is a holomorphic
function on the unit disk $D_1 (0) \subset \CC$, then the image
$f(D_1(0))$ covers some disk of radius $B_0\, |f'(0)|$.
\end{Lem}

   The case $|f'(0)| =1$ appears in \cite{Ah} and
the general case follows  from applying the case $|f'(0)|=1$ to
the function $g = f/|f'(0)| $. The version of Bloch's theorem that
we used, i.e., Lemma \ref{l:bloch},  follows by combining Lemma
\ref{l:bloch0} with the following uniformization result:

\begin{Lem}     \label{l:uni}
Given $r_0 > 0$ and $\kappa$, there exists a constant $\mu> 0$ so
that if $\Sigma$ is a surface  with $|K| \leq \kappa$ and  the
ball $B_{r_0}(p) \subset \Sigma$ is a topological disk, then there
is a holomorphic diffeomorphism $F: D_1(0) \to B_{r_0}(p)$ with
$F(0) = p$ and $|dF(p)|
> \mu$.
\end{Lem}

\begin{proof}
The existence of the holomorphic diffeomorphism $F$ follows
immediately from the uniformization theorem.  Namely, after
 extending a neighborhood of the disk $B_{r_0}(p) \subset \Sigma$, we can assume that it
 sits inside a closed Riemann surface.  By the uniformization theorem,
the universal cover of the closed Riemann surface is   either the
flat plane, the flat disk, or the round sphere.  Hence, the
topological disk $B_{r_0}(p)$ that sits inside the Riemann surface
must be conformal to a topological disk in $\CC$ (a topological
disk in the round sphere is conformal to one in the plane by
stereographic projection) and, by the Riemann mapping theorem,
also conformal to the unit disk $D_1(0)$.

Therefore, the point is to get the lower bound on $|dF(p)|$. Note
first that the inverse map $F^{-1}$ is a holomorphic function on
$B_{r_0} (p)$ that is bounded by one and vanishes at $p$. Since
$\Sigma$ has curvature bounded below by $-\kappa$ and each
component of a holomorphic function on a surface is automatically
also harmonic, the gradient estimate of \cite{ChY} implies that
\begin{equation}
    |\nabla F^{-1}(p)| \leq C \sup_{B_{r_0}(p)} |F^{-1}| = C \, ,
\end{equation}
for a fixed constant $C$ depending only on $r_0$ and $\kappa$.
This proves the lemma.
\end{proof}


\begin{thebibliography}{999}
\frenchspacing

\bibitem[Ah]{Ah}
L. Ahlfors,   An extension of Schwarz's lemma, {\emph{Trans. Amer.
Math. Soc.}} 43  (1938),  no. 3, 359--364.


\bibitem[A]{A}
N. Aronszajn,  A unique continuation theorem for solutions of
elliptic partial differential equations or inequalities of second
order, {\emph{J. Math. Pures Appl.}} (9) 36 (1957), 235--249.

\bibitem[Bo]{Bo}
J. Bourgain,   Fourier transform restriction phenomena for certain
lattice subsets and applications to nonlinear evolution equations.
I. Schr\"odinger equations,  {\emph{ Geom. Funct. Anal.}}  3
(1993), no. 2, 107--156.

\bibitem[BGT]{BGT}
N. Burq, P. G\'erard, and N. Tzvetkov,   Bilinear eigenfunction
estimates and the nonlinear Schr\"odinger equation on surfaces,
{\emph{Invent. Math.}}  159  (2005),  no. 1, 187--223.

\bibitem[CaL]{CaL}
M. Calle and D. Lee, in preparation.

\bibitem[ChY]{ChY}
S.Y. Cheng and S.T. Yau,  Differential equations on Riemannian
manifolds and their geometric applications, {\emph{Comm. Pure
Appl. Math.}} 28 (1975) 333--354.

\bibitem[CD]{CD}
T.H. Colding and C. De Lellis, Singular limit laminations, Morse
index, and positive scalar curvature,  {\it Topology}, Vol. 44,
Iss. 1 (2005) 25--45.

\bibitem[CH]{CH}T.H. Colding and N. Hingston,  Geodesic laminations with
closed ends on surfaces and Morse index: Kupka-Smale metrics,
{\emph{Comm. Math. Helv.}} 81 (issue 3)  (2006),  495 -- 522,
math.DG/0208133.

\bibitem[CM1]{CM1} T.H. Colding and W.P. Minicozzi II, Harmonic functions
with polynomial growth,  {\emph{ J. Differential Geom.}}  46
(1997),  no. 1, 1--77.

\bibitem[CM2]{CM2} T.H. Colding and W.P. Minicozzi II, Harmonic functions
on manifolds, {\emph{  Ann. of Math.}} (2)  146  (1997),  no. 3,
725--747.

\bibitem[CM3]{CM3}
T.H. Colding and W.P. Minicozzi II, Minimal surfaces, Courant
Lecture Notes in Math., v. 4, 1999.

\bibitem[Cv]{Cv}
Y. Colin de Verdi\`ere,  Sur le spectre des op\'erateurs
elliptiques \`a bicaract\'eristiques toutes p\'eriodiques,
{\emph{Comment. Math. Helv.}} 54 (1979), no. 3, 508--522.


\bibitem[GL]{GL}
N. Garofalo and F.H. Lin,   Unique continuation for elliptic
operators: a geometric-variational approach,  {\emph{Comm. Pure
Appl. Math.}} 40 (1987), no. 3, 347--366.


\bibitem[GT]{GT}
D. Gilbarg and N. Trudinger,   Elliptic partial differential
equations of second order. Second edition. Grundlehren der
Mathematischen Wissenschaften, 224. Springer-Verlag, Berlin, 1983.



\bibitem[Gu]{Gu}
V. Guillemin,  Lectures on spectral theory of elliptic operators,
{\emph{Duke Math. J.}} 44 (1977), no. 3, 485--517.

\bibitem[HL]{HL}
F. Hang and F.H. Lin, Exponential growth solutions of elliptic
equations,  {\emph{Acta Math. Sin.}}  15  (1999), no. 4, 525--534.

\bibitem[HsLa]{HsLa}
W. Hsiang and H.B. Lawson, Minimal submanifolds of low
cohomogeneity, {\em Jour. of Diff.  Geom.\/}, 5 (1971) 1--38.

\bibitem[JK]{JK}
D. Jerison and C. Kenig,   Unique continuation and absence of
positive eigenvalues for Schr\"odinger operators, {\emph{Ann. of
Math.}} (2) 121 (1985), no. 3, 463--494.





\bibitem[PW]{PW}
Y. Pan and T. Wolff,  A remark on unique continuation,  {\emph{J.
Geom. Anal.}} 8 (1998), no. 4, 599--604.

\bibitem[W1]{W1}
B. White,   The space of minimal submanifolds for varying
Riemannian metrics, {\emph{Indiana Univ. Math. J.}} 40 (1991), no.
1, 161--200.

\bibitem[W2]{W2}
B. White,  New applications of mapping degrees to minimal surface
theory, {\emph{J. Differential Geom.}}, 29 (1989), no. 1,
143--162.



\end{thebibliography}
\end{document}